\documentclass[a4paper,11pt,reqno]{amsart}
\def\ifundefined#1{\expandafter\ifx\csname#1\endcsname\relax}

\ifundefined{preambleloaded}
\usepackage[utf8]{inputenc}

\usepackage{amsmath}
\usepackage{amsthm}
\usepackage{amssymb}
\usepackage{amsfonts}
\usepackage{mathtools} 
\usepackage{mathrsfs} 
\usepackage{mleftright} 
\mleftright 

\usepackage[x11names]{xcolor}

\usepackage{tikz-cd}

\usepackage[breaklinks,colorlinks=true,linkcolor=blue,citecolor=green!60!black,backref]{hyperref} 
\hypersetup{pdftex}
\usepackage{hypcap}
\usepackage[capitalize,nameinlink]{cleveref} 
\crefname{equation}{}{}
\usepackage{autonum}

\usepackage[shortlabels]{enumitem}

\usepackage{csquotes}
\MakeOuterQuote{"}

\makeatletter
\def\@cite#1#2{[{#1\if@tempswa ,\penalty400~#2\fi}]}
\makeatother

\usetikzlibrary{arrows}
\usetikzlibrary{babel}
\tikzcdset{every label/.append style = {font = \small}}
\tikzset{
    arrowdecor/.style={anchor=south, rotate=90, inner sep=.5mm}
}
\usepackage{rotating}

\newtheorem{thm}{\iflanguage{ngerman}{Satz}{Theorem}}[section]
\newtheorem*{thm*}{\iflanguage{ngerman}{Satz}{Theorem}}
\Crefname{thm}{\iflanguage{ngerman}{Satz}{Theorem}}{\iflanguage{ngerman}{Sätze}{Theorems}}

\newtheorem*{conj*}{\iflanguage{ngerman}{Vermutung}{Conjecture}}
\Crefname{conj}{\iflanguage{ngerman}{Vermutung}{Conjecture}}{\iflanguage{ngerman}{Vermutungen}{Conjectures}}

\newtheorem{lem}[thm]{Lemma}
\newtheorem*{lem*}{Lemma}
\Crefname{lem}{Lemma}{Lemmas}

\newtheorem{prop}[thm]{Proposition}
\newtheorem*{prop*}{Proposition}
\Crefname{prop}{Proposition}{\iflanguage{ngerman}{Propositionen}{Propositions}}

\newtheorem{cor}[thm]{\iflanguage{ngerman}{Korollar}{Corollary}}
\Crefname{cor}{\iflanguage{ngerman}{Korollar}{Corollary}}{\iflanguage{ngerman}{Korollare}{Corollaries}}

\theoremstyle{definition}
\newtheorem{defn}[thm]{Definition}
\Crefname{defn}{Definition}{\iflanguage{ngerman}{Definitionen}{Definitions}}
\newtheorem*{notation}{Notation}
\newtheorem{defn-prop}[thm]{Definition/Proposition}

\theoremstyle{remark}
\newtheorem{rmk}[thm]{\iflanguage{ngerman}{Anmerkung}{Remark}}
\Crefname{rmk}{\iflanguage{ngerman}{Anmerkung}{Remark}}{\iflanguage{ngerman}{Anmerkungen}{Remarks}}

\newtheorem{ex}[thm]{\iflanguage{ngerman}{Beispiel}{Example}}
\Crefname{ex}{\iflanguage{ngerman}{Beispiel}{Example}}{\iflanguage{ngerman}{Beispiele}{Examples}}

\newtheoremstyle{claim}
  {0.5\topsep} 
  {0\topsep} 
  {} 
  {} 
  {\itshape} 
  {:} 
  {.5em} 
  {} 
\theoremstyle{claim}

\newtheorem*{claim*}{\iflanguage{ngerman}{Behauptung}{Claim}}
\Crefname{claim}{\iflanguage{ngerman}{Behauptung}{Claim}}{\iflanguage{ngerman}{Behauptungen}{Claims}}
\counterwithin*{claim}{thm}

\makeatletter
\newenvironment{pf}{\par
  \pushQED{\qed}%
  \list{}{\settowidth{\leftmargin}{\itshape{Pf~}}
          \settowidth{\labelwidth}{\itshape{Pf~}}
          \parsep=0pt \listparindent=\parindent
          \topsep=0pt
  }
\item[\hskip\labelsep\itshape{Proof:}]
}{%
    \popQED\endlist\@endpefalse\smallskip
    \aftergroup\@afterindentfalse\aftergroup\@afterheading
}
\makeatother

\newtheoremstyle{Normal}{}{}{}{}{\bfseries}{:}{.5em}{}
\theoremstyle{Normal}

\counterwithin*{problem}{section}
\makeatletter
\let\oldendproof\endproof
\def\endproof{\oldendproof\aftergroup\@afterindentfalse\aftergroup\@afterheading}
\makeatother
\Crefname{cd}{\iflanguage{ngerman}{Diagramm}{Diagram}}{\iflanguage{ngerman}{Diagramme}{Diagrams}}
\creflabelformat{cd}{(#2\thesection.#1#3)}
\usepackage{aliascnt}
\newaliascnt{cd}{equation}

\makeatletter
\newenvironment{cd}{%
    \refstepcounter{cd}%
    $$
        \begin{tikzcd}[sep=large]
}{%
    \end{tikzcd}
    \eqno \hbox{\@eqnnum}
    $$\@ignoretrue
}

\newenvironment{cd*}{%
    $$
    \begin{tikzcd}[sep=large]
}{%
    \end{tikzcd}
    $$\@ignoretrue
}
\makeatother
\makeatletter
\g@addto@macro\bfseries{\boldmath}
\makeatother

\newcommand{\noop}[1]{}

\Crefformat{item}{(#2#1#3)}

\newcommand{\surj}{\twoheadrightarrow}
\newcommand{\inj}{\hookrightarrow}
\newcommand{\weakto}{\rightharpoonup}
\newcommand{\iso}{\xrightarrow{\sim}}

\AtBeginDocument{
    \renewcommand{\phi}{\varphi}
    \renewcommand{\epsilon}{\varepsilon}
}

\renewcommand{\setminus}{\smallsetminus}

\DeclareMathSymbol{\colonrel}{\mathrel}{operators}{"3A}
\DeclareMathSymbol{:}{\mathpunct}{operators}{"3A}
\expandafter\let\expandafter\originall\csname\encodingdefault\string\l\endcsname
\DeclareRobustCommand*\l
    {\ifmmode\ell\else\expandafter\originall\fi}

\newcommand{\iprod}{\mathbin{\lrcorner}}
\newcommand{\acts}{\raisebox{-0.155em}{\scalebox{0.88}{ \rotatebox[origin=lb]{90}{$\curvearrowleft$}}}\,}

\usepackage{accents}
\DeclareFontFamily{U}{matha}{\hyphenchar\font45}
\DeclareFontShape{U}{matha}{m}{n}{
    <5> <6> <7> <8> <9> <10> gen * matha
    <10.95> matha10 <12> <14.4> <17.28> <20.74> <24.88> matha12
}{}
\DeclareSymbolFont{matha}{U}{matha}{m}{n}
\DeclareFontSubstitution{U}{matha}{m}{n}
\DeclareMathSymbol{\tinsubset}{3}{matha}{"80}
\let\oldbigwedge\bigwedge
\renewcommand{\bigwedge}{\textstyle\oldbigwedge}

\newcommand{\qtext}[1]{\quad\text{#1}\quad}

\newcommand{\Ric}{\mathrm{Ric}}
\newcommand{\Rm}{\mathrm{Rm}}

\newcommand{\SO}{\mathrm{SO}}
\newcommand{\GL}{\mathrm{GL}}

\DeclareMathOperator{\Ad}{Ad}
\DeclareMathOperator{\Aut}{Aut}

\DeclareMathOperator{\_div}{div}
\renewcommand{\div}{\_div}

\DeclareMathOperator{\grad}{grad}
\DeclareMathOperator{\Hom}{Hom}

\DeclareMathOperator{\Spec}{Spec}

\DeclareMathOperator{\id}{id}

\DeclareMathOperator{\sHom}{\mathscr{H}\mkern -4.5mu \mathit{om}}
\DeclareMathOperator{\supp}{supp}
\DeclareMathOperator{\tr}{tr}
\renewcommand{\Re}{\operatorname{Re}}
\renewcommand{\Im}{\operatorname{Im}}

\expandafter\let\expandafter\originald\csname\encodingdefault\string\d\endcsname
\DeclareRobustCommand*\d
    {\ifmmode\mathop{}\!d\else\expandafter\originald\fi}

\newcommand{\p}{\mathrlap{~~\text{.}}}
\expandafter\let\expandafter\originalc\csname\encodingdefault\string\c\endcsname
\DeclareRobustCommand*\c
    {\ifmmode\mathrlap{~~\text{,}}\else\expandafter\originalc\fi}

\newcommand{\demph}{\textbf}
\newcommand{\smashop}[1]{\;\;\smashoperator{#1}\;\;}

\makeatletter
\newcommand{\ovl}[1]{%
  \overline{\raisebox{0pt}[\dimexpr\height+0.6pt\relax]{\m@th$#1$}}%
}
\makeatother

\makeatletter
\let\save@mathaccent\mathaccent
\newcommand*\if@single[3]{%
    \setbox0\hbox{${\mathaccent"0362{#1}}^H$}%
    \setbox2\hbox{${\mathaccent"0362{\kern0pt#1}}^H$}%
    \ifdim\ht0=\ht2 #3\else #2\fi
}
\newcommand*\rel@kern[1]{\kern#1\dimexpr\macc@kerna}
\newcommand*\widebar[1]{\@ifnextchar^{{\wide@bar{#1}{0}}}{\wide@bar{#1}{1}}}
\newcommand*\wide@bar[2]{\if@single{#1}{\wide@bar@{#1}{#2}{1}}{\wide@bar@{#1}{#2}{2}}}
\newcommand*\wide@bar@[3]{%
    \begingroup
    \def\mathaccent##1##2{%
        \let\mathaccent\save@mathaccent
        \if#32 \let\macc@nucleus\first@char \fi
        \setbox\z@\hbox{$\macc@style{\macc@nucleus}_{}$}%
        \setbox\tw@\hbox{$\macc@style{\macc@nucleus}{}_{}$}%
        \dimen@\wd\tw@
        \advance\dimen@-\wd\z@
        \divide\dimen@ 3
        \@tempdima\wd\tw@
        \advance\@tempdima-\scriptspace
        \divide\@tempdima 10
        \advance\dimen@-\@tempdima
        \ifdim\dimen@>\z@ \dimen@0pt\fi
        \rel@kern{0.6}\kern-\dimen@
        \if#31
        \overline{\rel@kern{-0.6}\kern\dimen@\macc@nucleus\rel@kern{0.4}\kern\dimen@}%
        \advance\dimen@0.4\dimexpr\macc@kerna
        \let\final@kern#2%
        \ifdim\dimen@<\z@ \let\final@kern1\fi
        \if\final@kern1 \kern-\dimen@\fi
        \else
        \overline{\rel@kern{-0.6}\kern\dimen@#1}%
        \fi
    }%
    \macc@depth\@ne
    \let\math@bgroup\@empty \let\math@egroup\macc@set@skewchar
    \mathsurround\z@ \frozen@everymath{\mathgroup\macc@group\relax}%
    \macc@set@skewchar\relax
    \let\mathaccentV\macc@nested@a
    \if#31
    \macc@nested@a\relax111{#1}%
    \else
    \def\gobble@till@marker##1\endmarker{}%
    \futurelet\first@char\gobble@till@marker#1\endmarker
    \ifcat\noexpand\first@char A\else
    \def\first@char{}%
    \fi
    \macc@nested@a\relax111{\first@char}%
    \fi
    \endgroup
}
\makeatother
\def\semicolon{;}
\def\applytolist#1{
    \expandafter\def\csname multi#1\endcsname##1{
        \def\multiack{##1}\ifx\multiack\semicolon
            \def\next{\relax}
        \else
            \csname #1\endcsname{##1}
            \def\next{\csname multi#1\endcsname}
        \fi
        \next}
    \csname multi#1\endcsname}

\def\calchar#1{\expandafter\newcommand\csname c#1\endcsname{{\mathcal #1}}}
\applytolist{calchar}QWERTYUIOPLKJHGFDSAZXCVBNM;
\def\bbchar#1{\expandafter\newcommand\csname bb#1\endcsname{{\mathbb #1}}}
\applytolist{bbchar}QWERTYUIOPLKJHGFDSAZXCVBNM;
\def\bfchar#1{\expandafter\newcommand\csname bf#1\endcsname{{\mathbf #1}}}
\applytolist{bfchar}qwertyuiopasdfghjklzxcvbnmQWERTYUIOPLKJHGFDSAZXCVBNM;
\def\sfc#1{\expandafter\newcommand\csname s#1\endcsname{{\sf #1}}}
\applytolist{sfc}QWERTYUIOPLKJHGFDSAZXCVBNM;
\def\fchar#1{\expandafter\newcommand\csname f#1\endcsname{{\mathfrak #1}}}
\applytolist{fchar}qwertyuopasdfghjkzxcvbnmQWERTYUIOPLKJHGFDSAZXCVBNM;
\def\scchar#1{\expandafter\newcommand\csname sc#1\endcsname{{\mathscr #1}}}
\applytolist{scchar}QWERTYUIOPLKJHGFDSAZXCVBNM;

\else
\fi

\usepackage{geometry}
\usepackage{microtype}
\AtBeginDocument{
    \setlength{\abovedisplayshortskip}{-4pt plus 5pt minus 0pt}
}
\newenvironment{nicealign}{\begin{equation}\begin{aligned}}{\end{aligned}\end{equation}}

\renewcommand{\demph}[1]{\textbf{\textcolor{DarkOrange2}{#1}}}

\DeclareMathOperator{\Lie}{Lie}

\DeclareMathOperator{\Exp}{Exp}
\DeclareMathOperator{\Dom}{Dom}

\newcommand{\aut}{\mathfrak{aut}}

\DeclareMathOperator{\vol}{vol}

\DeclareMathOperator{\diam}{diam}
\newcommand{\PSH}{\mathrm{PSh}}
\DeclareMathOperator{\xxCap}{Cap}
\renewcommand{\Cap}{\xxCap} 

\DeclareMathOperator{\xxMA}{MA}
\DeclareMathOperator{\Isom}{Isom}
\DeclareMathSymbol{\colonrel}{\mathrel}{operators}{"3A}
\DeclareMathSymbol{:}{\mathpunct}{operators}{"3A}

\DeclareRobustCommand*{\MA}[1]{\ifmmode\xxMA#1\else{Monge--Ampère}\fi}
\def\DH/{Duistermaat--Heckman}

\newcommand{\MAX}{\MA_X}

\newcommand{\II}{I\mkern-2.5mu I}
\newcommand{\III}{I\mkern-2.5mu I\mkern-2.5mu I}
\newcommand{\IV}{IV}
\newcommand{\V}{V}
\newcommand{\VI}{V\mkern-2mu I}

\usepackage{bbm}
\newcommand{\ind}{\mathbf{1}} 

\newcommand{\pr}{\mathrm{pr}}

\renewcommand{\hat}{\widehat}
\renewcommand{\bar}{\widebar}
\renewcommand{\tilde}{\widetilde}

\newcommand{\Oalg}{\mathcal{O}^{\mathrm{alg}}}

\newcommand{\iddball}{(i \partial\bar{\partial})_{\mkern-2mu M, t}}
\newcommand{\ddat}[2]{\left.\frac{d}{d #1}\right|_{#2}}
\newcommand{\ddc}{\mathop{}\!dd^c}

\counterwithin{equation}{section}

\def\hlcolor{orange!70!black}
\newcommand{\hl}[1]{\textcolor{\hlcolor}{#1}}

\newcommand{\colorA}[1]{\textcolor{blue}{#1}}
\newcommand{\colorB}[1]{\textcolor{green!55!black}{#1}}

\newcounter{autoItemCounter}
\newcommand{\point}[1]{\textit{\labelcref{item:\theautoItemCounter#1}}}

\newlist{xxthmlist}{enumerate}{1}
\setlist[xxthmlist]{
    label=(\roman{xxthmlisti})\protect\label{item:\theautoItemCounter\thexxthmlisti},
    ref=(\roman{xxthmlisti})
}
\newenvironment{thmlist}{
    \refstepcounter{autoItemCounter}
    \begin{xxthmlist}
}{
    \end{xxthmlist}
}

\usepackage{thmtools}
\usepackage{thm-restate}

\newcommand{\thminflow}{\vspace{-\topsep}}

\newcommand{\myqedhere}{\tag*{\popQED}\smallskip}

\title[Uniqueness of AC shrinking gradient Kähler--Ricci solitons]{Uniqueness of asymptotically conical shrinking gradient Kähler--Ricci solitons}
\author{Carlos Esparza}
\address{Department of Mathematics, University of California, Berkeley, CA 94720, USA}
\email{esparzac@berkeley.edu}
\date{\today}

\begin{document}

\maketitle
\begin{abstract}
    We show that, up to biholomorphism, a given noncompact complex manifold only admits one shrinking gradient Kähler--Ricci soliton with Ricci curvature tending to zero at infinity. Our result does not require fixing the asymptotic data of the metric, nor fixing the soliton vector field.
    The method used to prove the uniqueness of the soliton vector field can be applied more widely, for example to show that conical Calabi--Yau metrics on a given complex manifold are unique up to biholomorphism.
    We also use it to prove that if two polarized Fano fibrations, as introduced by Sun--Zhang, are biholomorphic, then they are isomorphic as algebraic varieties.
\end{abstract}

\section{Introduction}%
\label{sec:introduction}

A \demph{Kähler--Ricci soliton} is a triple $(M, \omega, X)$ consisting of a complete Kähler manifold $(M, \omega)$ equipped with a complete real \emph{holomorphic}
vector field $X$ satisfying the equation
\begin{equation} \label{eq:KRS}
    \Ric(\omega) + \frac12 \cL_X \omega = \lambda \omega
\end{equation}
for some $\lambda \in \bbR$. Here $\Ric(\omega) \in \cA^{1,1}(M)$ denotes the Ricci curvature $(1,1)$-form of $\omega$.
By rescaling $X$ and $\omega$, the constant $\lambda$ may be assumed to be either $1, 0$ or $-1$, in which case the soliton is called \demph{shrinking}, \demph{steady} or \demph{expanding}. In this paper we will only consider the shrinking case and we assume the normalization $\lambda = 1$.

A Kähler--Ricci soliton is called \demph{gradient} if $X$ is a gradient vector field, i.e.\ if there exists an $f \in C^\infty(M)$, called the \demph{soliton potential} such that $X = \grad^\omega\! f$. In this case $\cL_X \omega = 2 i \partial\bar{\partial} f$ and \cref{eq:KRS} becomes
\[
    \Ric(\omega) + i \partial\bar{\partial} f = \lambda \omega
    \p
\]
For brevity we will simply refer to shrinking gradient Kähler--Ricci solitons as \demph{shrinkers}.

Our main result is
\thminflow
\begin{restatable}{thm}{uniq} \label{thm:uniq}
    Let $M^n$ be a complex manifold and let $(X_i, \omega_i), i=1,2$ be two shrinkers on $M$ such that $\lvert \Ric(\omega_i) \rvert_{\omega_i} \to 0$ at infinity.
    Then there exists a biholomorphism $\psi: M \to M$ such that $\omega_1 = \psi^* \omega_2$ and $\psi_* X_1 = X_2$.
\end{restatable}

This gives a positive answer to Question 2 in \cite[\S 7.2]{CDS} and generalizes \cite{TianZhu00,TianZhu02} to the noncompact case, at least in the case of shrinkers with Ricci curvature vanishing at infinity.

It was established in \cite[Thm.~A]{CDS} that if a complex manifold $M$ admits a shrinker $(\omega, X)$ with curvature as above then $M$ is a resolution of a Kähler cone $C_0$ and $\omega$ is asymptotic to a Kähler cone metric on the smooth part of $C_0$. As a consequence of the theorem above we also conclude that this metric on $C_0$ is unique up to automorphism.
\begin{cor} \label{thm:coneuniq}
    Let $M$ be a complex manifold admitting an asymptotically conical shrinker $(\omega, X)$. Then up to holomorphic isometry, the asymptotic Kähler cone of $(M, \omega, X)$ is uniquely determined by $M$.
\end{cor}

While this result shows that abstractly there is a special kind of Kähler cone metric that can arise as the limit of a soliton on the resolution of a cone, there is no known intrinsic characterization of these metrics.

As in the compact case, the proof of \cref{thm:uniq} consists of two independent parts. In \crefrange{sec:geodesics}{sec:uniq_along} we prove \cref{thm:uniq} in the special case $X_1 = X_2$. Then in \cref{sec:vf_uniq} we establish:
\begin{restatable}{thm}{uniqVF} \label{thm:uniqVF}
    Let $M$ be a complex manifold and let $(X_i, \omega_i), i=1,2$ be two shrinkers on $M$ with bounded Ricci curvature.
    Then there exists a biholomorphism $\psi: M \to M$ such that $\psi_* X_1 = X_2$.
\end{restatable}

This is just the volume minimization principle \cite[Thm.~D]{CDS}, except that we remove the hypothesis that $X_1$ and $X_2$ lie in the Lie algebra of a single torus acting on $M$ (or equivalently, that $X_1$ and $X_2$ commute).
The method developed in \cref{sec:vf_uniq} can be used to remove that same hypothesis from the volume minimization principle for Sasakian manifolds in \cite{MSYSasaki}, giving us

\begin{restatable}[uniqueness of CY cone metrics]{thm}{uniqCY} \label{thm:CYconeuniq}
    If $\omega_1$ and $\omega_2$ are Ricci-flat cone metrics on a complex manifold $Y$ with Reeb vector fields $X_1$ and $X_2$, then there exists $\psi \in \Aut(Y)$ such that $\psi_* \omega_1 = \omega_2$ and $\psi_* X_1 = X_2$.
\end{restatable}

\subsection{Motivation}%
\label{sub:motivation}

Shrinkers are interesting from the perspective of (Kähler--)Ricci flow.
A shrinker is equivalent to a self-similar Kähler--Ricci flow defined for times in $(-\infty, 1]$, where the metric evolves by pullback along the flow of the soliton vector field, while simultaneously rescaling. More precisely, the metric
\[
    \omega(t) = (1-t) \psi_t^* \omega, \quad t \in (-\infty, 1]
\]
satisfies the Kähler--Ricci flow equation $\partial_t \omega = - \Ric(\omega)$,
where $\psi_t$ is the flow of the time-dependent vector field $X / 2 (1-t)$.
By \cite{NaberShrinkers}, if a Kähler--Ricci flow develops a \emph{Type I} singularity (meaning that $|\Rm(\omega(t))|_{\omega(t)} = O((T - t)^{-1})$ as $t$ approaches the singular time $T$), the parabolic rescaling of the flow around a point where the curvature blows up will converge to the flow of a non-flat \cite{EMTtypeI} shrinker (in our sense) with bounded curvature.

Furthermore, shrinkers are interesting from the point of view of finding canonical metrics on Kähler manifolds.
Indeed, our main theorem asserts that shrinker metrics with suitably decaying curvature are canonical (and so are their soliton vector fields).
Shrinkers are a generalization of positive Kähler--Einstein metrics: A Kähler--Einstein metric on a Fano manifold is just a shrinker with zero soliton vector field, and while a noncompact manifold cannot admit a positive Kähler--Einstein metric by Myer's theorem, it can admit a shrinking soliton. The existence of a shrinking soliton on $M$ still implies that $M$ has positive first Chern class, at least in the sense that there exists a smooth hermitian metric on $K_M^{-1}$ with positive curvature. These ideas are explored further in \cite{SunZhang}, where a noncompact version of the Yau--Tian--Donaldson for Fano manifolds is proposed.
This work makes progress on the uniqueness part of that conjecture.

\subsection{Related work}%
\label{sub:related_work}

Examples of asymptotically conical shrinkers have been constructed using the Calabi Ansatz on the total spaces of the line bundles $\cO_{\bbP^{n-1}}(-k)$ for $1 \leq k < n$ \cite{FIK} and, more generally, over certain holomorphic vector bundles of rank 1 or 2 on Kähler--Einstein Fano manifolds, see \cite{CharlieExamples} and the references therein.

In the toric setting, the uniqueness of shrinkers has already been shown in \cite{CharlieThesis}, without any assumptions on the curvature. That work however does not obtain the full uniqueness of the vector field, it only shows uniqueness in the Lie algebra of a fixed torus.

Our result is similar, but distinct from Kotschwar--Wang's rigidity theorem \cite{KWrigidity}, which asserts that if two (Riemannian) shrinking gradient Ricci solitons on possibly different manifolds have isometric asymptotic cones, then they must be isometric near infinity.
However, Kotschwar--Wang's result leaves open the possibility that even on the same manifold there could be two shrinkers with non-isometric asymptotic cones. In the Kähler case this is excluded by our main theorem, leading to \cref{thm:coneuniq}.

In this respect our situation also differs from that of expanding asymptotically conical Kähler--Ricci solitons. Combining \cite[Thm.~A]{CDExpanders} and \cite[Thm.~A]{CDS} we see that if a manifold $M$ admits an asymptotically conical \emph{expanding} gradient Kähler--Ricci soliton then it will have an asymptotic cone $C$ which is uniquely determined as a complex manifold. However, for every Kähler cone metric $\omega_c$ on $C$ there will be an expanding gradient Kähler--Ricci soliton asymptotic to $\omega_c$, showing that for expanders \cref{thm:uniq} fails and \cref{thm:coneuniq} in some sense fails to the maximum possible degree.

\subsection{Methods and Outline}%
\label{sub:methods}

The bulk of this paper is concerned with the case where the soliton vector fields $X_1$ and $X_2$ agree.
Our proof of this case uses the techniques that \cite{BoBM} used to re-prove the uniqueness result of \cite{TianZhu00} in the compact case: In \cref{sec:geodesics}, we show that two shrinker metrics $\omega_0, \omega_1$ on $(M, X)$ can be connected by a very weak notion of geodesic (in the sense of Semmes \cite{SemmesGeod} and Donaldson \cite{DonaldsonGeod}) in a space of Kähler metrics that have quadratic growth in some sense (see \cref{def:Hs}). In \cref{sec:functionals} we define $X$-weighted energy and Ding functionals and in \cref{sec:uniq_along} we exploit their convexity properties, most notably a noncompact version of Berndtsson's convexity result \cite[Thm.~3.1]{BoBM}, to conclude that the geodesic is the pullback of one of the endpoints along the flow $F_t$ of a time-dependent vector field. In particular $\omega_1 = F_1^* \omega_0$, proving the special case of the theorem.
As in the compact case, this uniqueness result allows us to re-prove the Matsushima theorem for shrinkers, see \cref{thm:Matsushima}.

Compared to the compact case, additional difficulties arise when constructing geodesics by Perron's method since there is no easy choice of barrier function. This is because the difference of the Kähler potentials might not be bounded and the asymptotic cone metrics of $\omega_1$ and $\omega_2$ might even differ.
Since our construction of geodesics happens on a noncompact manifold, the standard theory of the \MA/-Dirichlet problem does not guarantee continuity, let alone $C^{1,1}$-regularity of the solution. This complicates the definition of the $X$-weighted \MA/ measure and $X$-weighted \MA/ energy which usually requires taking derivatives of the Kähler potential.
Here we draw on the methods of \cite[\S 2]{BWN} to define weighted \MA/ measures and -energies for locally bounded metrics, adapting their results to the asymptotically conical case.

The question of the uniqueness of the soliton vector field was already partially addressed in \cite{CDS}. Using the strategy of \cite{TianZhu02} and the \DH/ localization formula they prove that soliton vector fields are unique among all vector fields obtained from any maximal torus inside $\Aut(M)$. In \cref{sec:vf_uniq} we address the remaining question of showing that two shrinkers, whose soliton vector fields lie in two different maximal tori $T_1, T_2$ of automorphisms, are also equal up to automorphism. We tackle this problem by restricting the action of the tori to a non-reduced compact subscheme supported on the maximal compact analytic subset $E$ of $M$ (cf.\ \cref{def:AC}(iii)). Iwasawa's theorem then shows that $T_1$ and $T_2$ are conjugate up to terms vanishing on $E$ to arbitrarily high order.
Our final step is to show that if the difference between two Reeb vector fields (in the sense of \cite{SunZhang}; the soliton vector field is such a Reeb vector field) vanishes to sufficiently high order on $E$, then they are related by a biholomorphism.

This method can be applied more broadly in the context of polarized Fano fibrations (cf.\ \cite[\S 2]{SunZhang} or \cref{sub:polarized_fano_fibrations}).
One additional application is proving the uniqueness of Reeb vector fields for Calabi--Yau cone metrics, which gives us \cref{thm:CYconeuniq}.
By observing that the choice of a Reeb vector field on a polarized Fano fibration $M$ determines the sheaf of algebraic functions inside that of holomorphic functions, we can show, as a third application, that biholomorphic polarized Fano fibrations are isomorphic as algebraic varieties (\cref{thm:pffalg}).

\subsection*{Acknowledgements}
I would like to thank Song Sun for suggesting this problem, as well as for his constant support, encouragement and guidance.
I want to thank Junsheng Zhang for sending me a draft version of \cite{SunZhang} and for helpful discussions of its content.
I wish to thank Charles Cifarelli for useful discussions and many helpful comments on the draft version of this article,
and Ronan Conlon for feedback on the draft, and for suggesting \cref{thm:CYconeuniq}.
I also want to thank John Nolan and Qiuyu Ren for help with the proof of \cref{thm:Euniq}, and Garrett Brown, Connor Halleck-Dubé and Zachary Stier for helpful discussions.
Finally I wish to thank the Institute for Advanced Study in Mathematics at Zhejiang University for their hospitality during my visit in the Summer of 2024, during which part of this work was completed.
This work was partially funded by the NSF grant DMS-2304692.

\section{Geodesics in the space of asymptotically conical Kähler metrics}%
\label{sec:geodesics}

\subsection{Setup and Definitions}%
\label{sub:defs}

\begin{defn}
    On a product manifold $C \cong L \times \bbR_{>0}$, where $L$ is a compact $n-1$-manifold called the \demph{link}, we can have cone metrics. These are metrics of the form
    \[
        g_C = \d r^2 + r^2 g_L
    \]
    for some Riemannian metric $g_L$ on $L$ and $r: L \times \bbR_{>0} \to \bbR_{>0}$ given by $r(\ell, t) = t$. The pair $(M, g)$ is called a Riemannian cone.

    A \demph{Kähler cone} is a Riemannian cone equipped with a complex structure that makes the metric Kähler.
\end{defn}

It turns out that the Kähler $(1, 1)$-form of a Kähler cone satisfies $\omega_C = \frac{i}{2} \partial\bar{\partial} r^2$ \cite{MSYSasaki}.
A Kähler cone has a real-holomorphic \demph{radial vector field} $X = r \partial_r = r \grad^{\omega_C} r$.
In this paper we will always be considering Kähler cone metrics with a prescribed radial vector field. We then say that $\omega_C$ is a Kähler cone metric \emph{on} $(C, X)$.

The cone potentials of any two such metrics will be related by $r'^2 = e^{-\psi} r^2$ for some smooth function $\psi$ on the link.
Since $\psi$ is bounded, this means that any two cone potentials are equivalent, i.e.\ their ratio is globally bounded.

We can add a single vertex point $o$ to any Kähler cone, turning it into a normal affine variety with a single singularity at $o$ \cite[Thm.~3.1]{vanCoevering}.
We write $C_0 = C \cup \{o\}$.

In this paper, we will consider an $n$-dimensional Kähler manifold $(M, \omega, X)$ with a holomorphic vector field that is asymptotic to a Kähler cone $(C, \omega_C, r \partial_r)$ in the following (rather strong) sense:
\begin{defn}[cf.\ {\cite[Thm.~A]{CDS}}] \label{def:AC}
    We say $M$ is \demph{asymptotically conical} if
    \begin{enumerate}[(i)]
        \item There exists a holomorphic \demph{resolution map} $\pi: M \to C_0$ to a cone with vertex $C_0$, which is a biholomorphism onto the smooth part $C$ of the cone.
        \item $\pi_* X = r \partial_r$ away from the \demph{exceptional set} $E = \pi^{-1}(\{o\})$.
        \item The metric of $M$ has the asymptotic behavior
            \begin{equation}
                \big\lvert (\nabla^{\omega_C})^k( \pi_* \omega - \omega_C ) \big\rvert_{\omega_C} = O(r^{-2-k})
                \p
            \end{equation}
        \item $J X$ is $\omega$-Killing.
    \end{enumerate}
    Via $\pi$, we can treat the cone metric $\omega_C$ as a metric on $M \setminus E$ and $r$ as a continuous function on $M$, smooth on $M \setminus E$. We will use it as a background metric.
\end{defn}
From this definition we see that $X$ and $JX$ are complete because the flows of $\pi_* X$ and $J \pi_* X$ exist for all times.
Since $J X$ acts (faithfully) on the $\SO$-frame bundle of $M$ restricted to $\{r = 1\}$, the closure of its flow in $\Isom(X, \omega)$ is compact and abelian, so it is a torus that acts faithfully on $(M, \omega)$ by isometries. We call this the \demph{torus generated by $JX$}.
We will now pick and fix a torus $T$ in the holomorphic isometry group that contains the torus generated by $JX$.

\cref{def:AC} is satisfied when $M$ admits a shrinker with quadratic curvature decay \cite[Thm.~A]{CDS}, i.e.
$|\Rm(x)| \leq C d(p, x)^2$
for some $C > 0$ and $p \in M$.
This is essentially the only situation of interest in this paper. However, for future work, we do not directly assume that $M$ admits a shrinker.
Munteanu--Wang have shown that for shrinking gradient Ricci solitons, quadratic curvature decay is implied by vanishing Ricci curvature at infinity (as in the hypothesis of \cref{thm:uniq}) \cite[Thm.~2.2]{MWconical}, and that in turn is implied by having bounded Riemann tensor and scalar curvature decay at infinity \cite[Thm.~2.3]{MWconical}.

\begin{notation}
    \begin{enumerate}[(i)]
        \item If $f$ and $g$ are functions on $M$, we write $\{f \leq g\}$ for the set $\{p \in M \mid f(p) \leq g(p)\}$.
        \item We define $B_R \coloneqq \{r < R\}$ to be the open "ball" of radius $r$.
    \end{enumerate}
\end{notation}

\subsubsection{Equivariant charts}%
\label{sec:equivariant_charts}

We always use $z_k = x_k + i y_k, k=1, \dots, n$ as coordinates on (open subsets of) $\bbC^n$.

\begin{lem}[equivariant charts]
    There exist finitely many holomorphic charts
    \begin{equation} \label{eq:defEqCh}
        U_\alpha \cong \bbR \times i (-1, 1) \times U'_\alpha \quad\text{where}\quad U'_\alpha \subseteq \bbC^{n-1} \text{ is open}
    \end{equation}
    covering $M \setminus E$, such that $X|_{U_\alpha} = \frac{\partial }{\partial x_1^\alpha}$. In these charts $r = G_\alpha(z^\alpha_2, \dots, z^\alpha_n) \exp(x^\alpha_1)$ for some smooth functions $G_\alpha$.
    We refer to this sort of charts as \demph{equivariant charts}.
\end{lem}
\begin{proof}
    Since $X$ and $JX$ commute and are complete, we can use their flows to obtain holomorphic charts of the form $\bbR \times i(-1, 1) \times U_\alpha$. By compactness, finitely many suffice to cover the link of $C \cong M \setminus E$, and these actually cover all of $M \setminus E$.
    Regarding the last assertion: On the cone, $X = r \frac{\partial}{\partial r} = \frac{\partial }{\partial x^\alpha_1}$ implies $\frac{d r}{d x^\alpha_1} = r$, so $r = C(y^\alpha_1, z^\alpha_2, \cdots) \exp(x^\alpha_1)$. Since $J X = \frac{\partial }{\partial y^\alpha_1}$ is Killing, it needs to preserve $r$ and the claim follows.
\end{proof}

\newcommand{\charts}{\cA_M}
We now fix a finite atlas $\charts$ for $M$ consisting of charts that are equivariant or defined on a precompact set in $M$ (the latter are used to cover $E$).
By slightly shrinking $U'_\alpha$ and replacing $(-1, 1)$ by $(-1 + \epsilon, 1 - \epsilon)$ in \cref{eq:defEqCh} (and also shrinking the precompact charts), we can ensure that the Jacobians of the transition functions are bounded. This is essentially because our "shrunken" charts are precompact subsets of the original charts. That statement is not actually true, because of the $\bbR$-factor in \cref{eq:defEqCh}, but since all relevant quantities are invariant in that direction the conclusion still holds.

\begin{notation}
    Let $L \to M$ be a holomorphic line bundle.
    We will use "additive notation" for hermitian metrics: Locally, there is a correspondence between functions $\phi \in C^\infty(M)$ and hermitian metrics $e^{-\phi} \in C^\infty(M; L^* \otimes \bar{L}^*)$. Globally, $\phi$ cannot be considered to be a function, but only a section of some affine bundle.
    Even though calling $\phi$ a potential would be more accurate, we will refer to it as a metric. The curvature of the metric $\phi$ is $\partial\bar{\partial} \phi$, which is globally well-defined.
    We will mostly be concerned with \emph{singular} hermitian metrics, for which $\phi$ is not locally $C^\infty$, but only plurisubharmonic (PSh). We will write $\PSH(M; L)$ for the set of non-negative singular hermitian metrics on $L$.
    An alternative point of view would be to choose a smooth background metric $e^{-\psi}$ on $L$, and then identify $\PSH(M; L)$ with the set of $\theta$-PSh functions for $\theta = i \partial\bar{\partial} \psi$.
\end{notation}

\newcommand{\Hsing}{\bar{\cH}_2}

\begin{notation} We will use asymptotic notation in the following way:
    \begin{enumerate}
        \item By default, all asymptotic notation in this paper is "as $r \to \infty$."

        \item "Big $\Theta$" notation $f = \Theta(g(r))$ means that $\exists C > 0$ such that $C^{-1} g(r) \leq f \leq C g(r)$ for $r \gg 1$.

        \item When $\phi$ is a metric (in "additive notation") we write $\phi = O(g(r))$ to mean that $\phi_\alpha = O(g(r))$ for the trivializations $\phi_\alpha$ of $\phi$ in all charts of $\charts$ (of course the precompact charts in $\charts$ do not matter for asymptotic statements like this).
            More precisely, it means
            \[
                - \log\left(\frac{e^{-\phi}}{i^{n^2} \Omega_\alpha \wedge \bar{\Omega}_\alpha}\right) = O(g(r)) \quad\text{as $r \to \infty$}
            \]
            for every chart $z^\alpha_1, \cdots z^\alpha_n$ in $\charts$, where $\Omega_\alpha = \d z^\alpha_1 \wedge \dots \wedge \d z^\alpha_n$. Equivalently, we can also say
            \begin{equation} \label{eq:bigOalt}
                -\log \left( \frac{e^{-\phi}}{r^{-2n} \omega_C^n} \right) = O(g(r)) \quad\text{as $r \to \infty$}
            \end{equation}
            since the $X$-invariant volume form $\Upsilon \coloneq r^{-2n} \omega_C^n$, defined on $\{r \neq 0\}$ is uniformly equivalent to each $\Omega_\alpha \wedge \overline{\Omega_\alpha}$ on $U_\alpha$ (this is because by construction $U_\alpha$ is a precompact subset of a larger chart).

    \end{enumerate}
\end{notation}
\begin{rmk}
    Since all cone metrics (with fixed radial vector field) on $C_0$ are equivalent, we can use them interchangeably in asymptotic notation, e.g.\ $f = O(\log r_1) \iff f = O(\log r_2)$ for any two cone metrics $r_1, r_2$. Thus we will just write $O(\log r)$ for the background metric from \cref{def:AC}.

    Also note that in the important special case $g(r) = r^2$ we can rephrase $\psi = O(r^2)$ as $\psi_\alpha = O(e^{2x^1_\alpha})$ in every chart in $\charts$, since $r$ is equivalent to $x^1_\alpha$ in every equivariant chart.
\end{rmk}

\begin{defn} \label{def:Hs}
    Let
    \begin{nicealign}
        \bar{\cH} &= \PSH(M; K_M^{-1}) \cap L^\infty_\mathrm{loc} \\
        \cH &= \left\{ \psi \in \PSH(M; K_M^{-1}) \cap C^\infty ~\middle|~ i \partial\bar{\partial} \psi > 0 \right\}
    \end{nicealign}
    and
    \begin{equation} \label{eq:defHs}
    \begin{aligned}
        \cH_2 &= \left\{ \psi \in \cH ~\middle|~ \psi = \Theta(r^2) \right\}
        \\
        \cH^*_2 &=
        \left\{ \psi \in \cH_2 ~\middle|~
            \begin{array}{l}
                \exists \text{ cone metric } \omega_C = \frac{i}{2} \partial\bar{\partial} \tilde{r}^2 \text{ on $(C, r \partial_{r})$, $\rho \coloneqq \psi - \tilde{r}^2/2$:} \\
                \rho = O\left((\log r)^2\right), \quad
                |\partial \rho |_{\omega_C} = O(r^{-1} \log r), \quad
                |i \partial\bar{\partial} \rho |_{\omega_C} = O(r^{-2})
            \end{array}
        \right\}
        \\
        \Hsing &= \left\{ \psi \in \bar{\cH} ~\middle|~ \psi = \Theta(r^2) \right\}
        \p
    \end{aligned}
    \end{equation}
    Since we will need more quantitative bounds later we define $\cH_2(a, b)$, $\Hsing(a,b)$, etc.\ to consist of those $\psi$ for which
    \[
        a r^2 - a^{-1} < \psi < b r^2 + b
    \]
    in all charts of $\charts$.

    Finally $\cH^T$, $(\cH_2^*)^T$, etc.\ denote the subsets of $T$-invariant metrics.
\end{defn}

The decay conditions in the definition of $\cH^*_2$ come from
\thminflow
\begin{lem} \label{thm:KRSasymptotics}
    If $\omega_S$ is a shrinker with quadratic curvature decay then $\omega_S = i \partial\bar{\partial} \phi_S$ for some $\phi_S \in (\cH_2^*)^T$.
\end{lem}
\begin{proof}
    We can rewrite the shrinker equation \cref{eq:KRS} as
    \[
        -i \partial\bar{\partial} \log(\omega_S^n) + i \partial\bar{\partial} f = \omega_S
    \]
    so we start by setting $\phi_S = -\log(\omega_S^n) + f$ and $\omega_S = i \partial\bar{\partial} \phi_S$.
    On every equivariant chart $V_\alpha$, we define $\rho_\alpha$ such that
    \[
        e^{-\phi_S} = e^{-r_C^2 /2 -\rho_\alpha} i^{n^2} \Omega_\alpha \wedge \bar{\Omega}_\alpha
        \p
    \]
    From \cite[Thm.~A(c)]{CDS} we know that
    $\lvert i \partial\bar{\partial} \rho_\alpha \rvert_{\omega_C} = \left\lvert \omega_S - \omega_C \right\rvert_{\omega_C} = O(r^{-2})$. Since $\omega_C$ has $X$-weight 2 this is equivalent to saying that the coefficients of $i \partial\bar{\partial} \rho_\alpha$ are $O(1)$ in equivariant coordinates.
    Since $\cL_{JV} \partial \rho = 0$, we have
    \begin{align}
        \cL_X \partial \rho_\alpha &= \cL_X \partial \rho_\alpha = \cL_{X^{0,1}} \partial \rho_\alpha \\
                    &= X^{0,1} \iprod \bar{\partial} \partial \rho_\alpha
        \c
    \end{align}
    so the coefficients of $\cL_X \partial \rho_\alpha$ are also $O(1)$ and integrating ($\cL_X = \partial_1$ in equivariant coordinates) we conclude that $\partial \rho_\alpha = O(x)$ in equivariant coordinates. This is equivalent to $| \partial \rho_\alpha |_{\omega_C} = O(e^{-x} x) = O(r^{-1} \log r)$.
    Conjugating we also obtain $d \rho_\alpha = O(x)$, so another integration gives $\rho_\alpha = O(x^2) = O( (\log r)^2 )$
\end{proof}
\begin{rmk} \label{rmk:higherAsymp}
    In every equivariant chart $U_\alpha$ we have $r^2 = e^{2x^\alpha} f_\alpha(z_2^\alpha, \dots, z_n^\alpha)$ for some $f_\alpha$, thus a direct computation shows that the Christoffel symbols of the cone metric $\omega_C = i \partial\bar{\partial} r^2/2$ in these coordinates are independent of $z_1^\alpha$. This implies that their components, and all their coordinate derivatives are bounded.
    Thus the higher order asymptotics $|\nabla^k (\omega - \omega_C)|_{\omega_C} = O(r^{-2-k})$ from \cite[Thm.~A(c)]{CDS} allow us to conclude that the coefficients of all (coordinate) derivatives of $i \partial\bar{\partial} \rho_\alpha$ are bounded in equivariant charts. So
    \[
         \frac{\partial \rho}{\partial z^\alpha \partial z^\beta}  = O(1)
    \]
    for any multi-indices $\alpha, \beta \neq 0$.
    Integrating as above we then obtain for $\alpha \neq 0$
    \[
         \frac{\partial \rho}{\partial z^\alpha} = O(x) = O(\log r)
         \p
    \]
\end{rmk}

\subsection{Barrier functions and existence of geodesics}%
\label{sub:barrier_functions_and_existence_of_geodesics}

\begin{notation}
    We write $I = [0, 1]$ and $\Omega = (0, 1) \times i \bbR \subseteq \bbC$. We will use $t$ as a coordinate on $I$ and $\tau = t + i s$ on $\Omega$.

    When dealing with functions/metrics $f$ defined on products like $\Omega \times M$ we write $f_t \coloneqq f|_{\{t\} \times M}$.
    Furthermore, we write $d f, \partial f, \bar{\partial} f$ to denote only the part of the differential of $f$ in the "$M$-directions" and $d_{t, M} f = d f + \frac{\d f}{\d t} \d t$ to denote the full differential.
\end{notation}

\begin{defn}
    A \demph{geodesic} connecting metrics $\phi_0, \phi_1 \in (\cH^*_2)^T$ is a locally bounded PSh metric $\Phi \in \PSH(\Omega \times M; \pr_2^* K_M^{-1} )^T$, constant in the $i \bbR$-direction, such that there exist $a, b$ for which $\Phi_t \in \Hsing(a, b)$ for all $t$, and which solves the boundary value problem
    \[
        \begin{cases}
            \left(\iddball \Phi\right)^{n+1} = 0 \\
            \lim_{t \to i} \Phi_t = \phi_i & i = 0, 1
        \end{cases}
        \p
    \]
    The limit is in the pointwise sense. A \demph{subgeodesic} has the same definition, except we only require $(\iddball \Phi)^{n+1} \geq 0$.
    Since the geodesics provided by the next theorem satisfy $\Phi_t - \Phi_s = O(r^2) |t - s|$, we include this property in our definition of (sub)geodesics.
\end{defn}

The first step in proving \cref{thm:uniq} is establishing the existence of geodesics:
\begin{thm} \label{thm:geod_exist}
    Let $\phi_0, \phi_1 \in (\cH^*_2)^T$. Then there exists a $T$-invariant \emph{geodesic} $\Phi$ connecting $\phi_0$ to $\phi_1$.
\end{thm}

\newcommand{\proofstep}[1]{\subsubsection*{#1}}

\begin{proof}
    Consider the subset $\cB \subseteq \PSH\left(\Omega \times M; \pr_2^* K_M^{-1}\right) \cap L^\infty_\mathrm{loc}$ consisting of those PSh metrics which are independent of $\Im t$
    .
    We claim that the upper-semicontinuous regularization
    \begin{equation}
        \label{eq:perron}
        \Phi(t, x) = \left(\sup \left\{\theta(t, x) ~\middle|~ \theta \in \cB, \quad \limsup_{t \to 0,1} \theta_t \leq \phi_{0, 1} \right\} \right)^*
    \end{equation}
    defines a locally bounded PSh metric (in particular we are asserting that $\Phi(t, x)$ is finite at every point).

    \proofstep{Step 1: Constructing a barrier function.}%
    Since $\phi_0, \phi_1$ lie in $(\cH_2^*)^T$ there exist (potentials for) cone metrics $r_i^2 = e^{\psi_i} r^2$ such that $\phi_i - r_i^2/2 = O((\log r)^2)$.

    Let $g: \bbR \to [0, 1]$ be a smooth function such that $g|_{[0,1]} = 0$ and $g|_{\bbR_{\geq2}} = 1$.
    Now $g(r) r_i^2$ are well-defined functions on $M$, and we can define $\rho_i = \phi_i - g(r) r_i^2/2$, effectively decomposing the metrics into a "conical part" and a part along the exceptional set $E$.

    For $A, B > 0$, define a metric on $\pr_2^* K_M^{-1} \to \bar{\Omega} \times M$, which is independent of the imaginary direction in $\bar{\Omega}$, and for $t \in I$ is given by
    \begin{align} \label{eq:barrier0}
        \tilde{\chi}_t = g(r) \frac{r_t^2}{2} e^{-A t (1-t)}
            + (1 - t) \rho_0 + t \rho_1
            + 2 B t (t-1)
    \end{align}
    where
    \begin{nicealign}
        r_t^2 &= e^{\psi_t} r^2 = (1-t) r_0^2 + t r_1^2 \\
        \psi_t &= \log \left( (1-t) e^{\psi_0} + t e^{\psi_1} \right)
        \p
    \end{nicealign}
    We will specify the parameters $A, B > 0$ later.
    The distance functions $r$ and $r_t$ for $t \in [0, 1]$ are uniformly equivalent, with the constant of equivalence depending only on $\phi_0$ and $\phi_1$.
    Similarly, all cone metrics $\omega_t \coloneqq \frac{i}{2} \partial\bar{\partial} r_t^2$ are uniformly equivalent on $M \setminus E$ because they all have $X$-weight 2 (i.e.\ $\cL_X \omega_t = 2 \omega_t$) and they are uniformly equivalent when restricted to the (compact) link.
    From \cref{eq:barrier0} we immediately see that $\tilde{\chi}_0 = \phi_0, \tilde{\chi}_1 = \phi_1$ and on $B_1$ (recall we write $B_1 \coloneq \{r \leq 1\}$)
    \begin{align}
        \tilde{\chi}_t &\equiv (1-t) \phi_0 + t \phi_1 + 2 Bt (t-1)
        \label{eq:chinearE}
        \p
    \end{align}
    For $r \geq 2$ (so $g(r) = 1$) we compute
    \begin{align}
        \label{eq:ddbchi}
        \beta \coloneq i \partial\bar{\partial} \tilde{\chi}_t
                &= e^{-A t (1-t)} i \partial\bar{\partial} \frac{r_t^2}{2} + (1-t) i \partial\bar{\partial} \rho_0 + t i \partial\bar{\partial} \rho_1 \\
        \label{eq:dtchi}
        \frac{d}{d \bar{\tau}} \tilde{\chi}_t &= \frac12 \frac{d}{d t} \tilde{\chi}_t = \frac{r_t^2}{4} e^{-A t (1-t)} \left( \dot{\psi}_t + A (2t-1) \right) + \frac{\rho_1 - \rho_0}{2} + B (2t -1)\\
        \label{eq:dttchi}
        f \coloneq \frac{d^2}{d \tau d \bar{\tau}} \tilde{\chi}_t &= \frac{r_t^2}{8} e^{-A t (1-t)} \left[ ( \dot{\psi}_t + A (2t-1) )^2 + 2A + \ddot{\psi}_t \right] + B \\
        \label{eq:dmixchi}
        \alpha \coloneq \partial \frac{d}{d \tau} \tilde{\chi}_t &= \frac14 e^{-A t (1-t)} \left( (\dot{\psi}_t + A (2t-1)) 2 r_t \partial r_t + r_t^2 \partial \dot{\psi}_t \right) + \frac{\partial(\rho_1 - \rho_0)}{2}
        \p
    \end{align}
    Note that $d^2 \tilde{\chi}_t / d\tau d\bar{\tau}$ is the only second derivative that depends on $B$.
    The real $(1,1)$-form
    \[
        \iddball \tilde{\chi} = i \beta + i \alpha \wedge \d \bar{\tau} - i \bar{\alpha} \wedge \d \tau + i f \d \tau \wedge \d \bar{\tau}
    \]
    is positive if{}f $f > 0, \beta > 0$ and $i \alpha \wedge \bar{\alpha} < f \beta$. The condition $f > 0$ is easy to guarantee if we pick $B$ sufficiently large.
    Recall that%
    \footnote{We use the norm on differential forms where if $\epsilon_i$ are orthonormal then $\epsilon_1 \wedge \dots \wedge \epsilon_k$ has unit length.}
    for any Hermitian metric with $(1,1)$-form $\omega$ we have $i \alpha \wedge \bar{\alpha} \leq |\alpha|^2_\omega \omega$,
    so to establish the third condition it suffices to check that $|\alpha|^2_{\omega_t} \leq f \beta$. If we write out this equation, multiplying both sides by $16 e^{At(1-t)}$ and writing $\Gamma = \dot{\psi}_t + A(2t-1)$ for brevity we obtain
    \begin{align}
        & \left\lvert 2 \Gamma r_t \partial r_t + r_t^2 \partial\dot{\psi}_t + 2 e^{At(1-t)} \partial(\rho_1 - \rho_0) \right\rvert^2_{\omega_t} \omega_t \\
            &\qquad
            < 2 \left(
                r_t^2 ( \Gamma^2 + 2A + \ddot{\psi}_t )
                + e^{A t (1-t)} B
              \right)
              \left(
                \omega_t +
                e^{A t (1-t)} \left[ (1-t) i \partial\bar{\partial}\rho_0 + t i \partial\bar{\partial} \rho_1 \right]
              \right)
        \label{eq:testPos}
        \p
    \end{align}
    \proofstep{Step 1.1. Positivity on $\{r \geq R\}$ for $R \gg 1$.}

    We will allow the implicit constants in big-O notation to depend on the boundary values $\phi_0, \phi_1$ (but not on the parameters $A$ and $B$).

    The forms $\partial r_t$ and $\partial\dot{\psi}_t$ are $\omega_t$-perpendicular and furthermore $2 |\partial r_t|_{\omega_t} = \lvert \d r_t \rvert_{\omega_t} = 1$, so
    \begin{equation} \label{eq:wedgeEst}
        \left| 2\Gamma r_t \partial r_t + r_t^2 \partial\dot{\psi}_t \right|^2_{\omega_t} \omega_t = \left( 2 r_t^2 \Gamma^2 + r_t^4 \lvert \partial\dot{\psi}_t\rvert_{\omega_t}^2 \right) \omega_t
        \p
    \end{equation}
    Thus we can estimate the LHS of \cref{eq:testPos} by
    \begin{align}
        \cdots &\leq \Big(
            2 r_t^2 \Gamma^2 
                + r_t^4 \lvert \partial \dot{\psi}_t\rvert_{\omega_t}^2
                + 4 e^{A t(1-t)} \left\langle 2 \Gamma r_t \partial r_t + r_t^2 \partial(\psi_0 - \psi_1), \partial(\rho_1 - \rho_0) \right\rangle_{\omega_t}
            \\& \qquad
                + 4 e^{2 A t(1-t)} |\partial(\rho_1 - \rho_1)|_{\omega_t}^2
        \Big) \omega_t
        \p
    \end{align}
    Now we use \cref{thm:KRSasymptotics} together with the fact that $|r_t \partial r_t|_{\omega_t}$ and $|r_t^2 \partial(\psi_0 - \psi_1)|_{\omega_t}$ have $X$-weight 1 and hence are $O(r)$, yielding
    \begin{align} \label{eq:UB}
        \cdots &\leq \Big(
            2 r_t^2 \Gamma^2
            +C r_t^2 + e^A (1 + \Gamma) O(\log r)
            + e^{2A} O\left(r^{-2}(\log r)^2\right)
        \Big) \omega_t
        \end{align}
    where
    \[
        C = \sup_{t \in I} \sup_L | \partial \dot{\psi}_t|^2_{\omega_t} r_t^2
          = \sup_{t \in I} \sup_{M \setminus E} | \partial \dot{\psi}_t|^2_{\omega_t} r_t^2
      \p
    \]
    Next, observe that
    the RHS of of \cref{eq:testPos} is bounded from below by
    \begin{equation} \label{eq:LB}
        2 \left(\Gamma^2 + 2A + \ddot{\psi}_t\right) r_t^2 \left(1 - e^A O(r^{-2})\right) \omega_t
        \p
    \end{equation}
    Now we fix some $A$ such that $2\left(2A - \sup_{I \times L} |\ddot{\psi}|\right) > C$. Then the coefficient of the (leading order) quadratic term in \cref{eq:LB} is greater than that of \cref{eq:UB}, so
    there exists an $R > 0$ such that the inequality \cref{eq:testPos} is satisfied for any $r \geq R, t \in [0, 1]$ and any $B > 0$.

    We will also pick $R$ large enough to ensure that $i \partial\bar{\partial} \tilde{\chi}_t > \kappa \omega_C$ for some $\kappa > 0$. This is possible by \cref{eq:ddbchi} since $|i \partial\bar{\partial} \rho_i|_{\omega_C} = O(r^{-2})$.

    \proofstep{Step 1.2. Fixing spatial positivity for small $r$.}

    Pick a smooth function $\sigma: \bbR \to [0, 1]$, identical to $1$ on $(-\infty, 0]$ and to $0$ on $[1, \infty)$. Put
    \[
        \sigma_{R, L}(x) \coloneqq \sigma\left(\frac{x - R}{L}\right)
    \]
    (where $R$ is the threshold determined in the last step and $L$ is to be determined).
    Further, pick a function $q: [0, \infty) \to [0, \infty)$, depending on a parameter $\epsilon$, with the properties
    \[
        \begin{cases}
            q(r) = 0 & \text{for } r \leq 1/3 \\
            q(r) = \frac12 r^2 & \text{for } 1 \leq r \leq R \\
            R \leq q'(r) \leq 2R  & \text{for } R \leq r \\
            q'' \geq 0
        \end{cases}
        \p
    \]
    Note that these conditions together imply $q(r) \leq 2 R r$ and $q' \geq 0$.
    One way to obtain such a function would be to integrate $\sigma_{R, 2R}$ twice, and then modify it to vanish around zero.

    For any $L, D > 0$ we now define a family of metrics
    \begin{equation} \label{eq:chi}
        \chi_t = \tilde{\chi}_t + D t (1-t) q(r) \sigma_{R, L}(r)
        \p
    \end{equation}
    Observe that $\chi$ agrees with $\tilde{\chi}$ for $r > R + L$ and hence is positive (jointly in $M$- \emph{and} $t$-directions) on $\{r > R + L\}$.

    Further note that
    \[
        i \partial\bar{\partial} (q \circ r) 
            \equiv q'(r) i \partial\bar{\partial} r + q''(r) i \partial r \wedge \bar{\partial} r \geq 0
        \c
    \]
    so on $B_R$, where we can ignore the $\sigma_{R, L}$ in \cref{eq:chi}, we get $i \partial\bar{\partial} \chi_t \geq i \partial\bar{\partial} \tilde{\chi}_t$.
    In particular, by \cref{eq:chinearE} we have $i \partial\bar{\partial} \chi_t \geq i \partial\bar{\partial} \tilde{\chi}_t > 0$ on $B_1$.

    Next, on $\{1 \leq r \leq R\}$ the metrics $\chi_t$ and $\tilde{\chi}_t$ differ by $D t (1-t) \frac{r^2}{2}$.
    On $\{0, 1\} \times M$ we have $i \partial\bar{\partial} \tilde{\chi}_t = i \partial\bar{\partial} \phi_{0, 1} > 0$, so by compactness of $B_R \times \bar{\Omega}$ there exists $\delta > 0$ such that $i \partial\bar{\partial} \tilde{\chi}_t > 0$ on $\{r \leq R\} \times ([0, \delta] \cup [1-\delta, 1])$.
    So if we pick $D$ sufficiently large, we can arrange that on $\{2/3 < r < R\}$
    \[
        i \partial\bar{\partial} \chi_t = i \partial\bar{\partial} \tilde{\chi}_t + D t (1-t) i \partial\bar{\partial} \frac{r^2}{2}> 0
    \]
    ($>0$ here only means positivity in $M$-directions).

    \newcommand{\csig}{\lVert\sigma\rVert_{C^2}}

    Finally consider the "roll-off region" $\{R \leq r \leq R + L\}$.
    First compute
    \begin{align}
        i \partial\bar{\partial} \chi_t &\equiv i \partial\bar{\partial} \tilde{\chi}_t
            + D t (1-t) \Big[
                \big(
                    \colorA{\sigma_{R, L}'' q} +
                    \colorB{2\sigma_{R, L}' q'} +
                    \overbrace{\sigma_{R, L} q''}^{\geq 0}
                \big) i \partial r \wedge \bar{\partial} r
                + \frac{
                    \colorB{\sigma_{R, L}' q} +
                    \overbrace{\sigma_{R, L} q'}^{\geq 0}
                }{r} i r \partial\bar{\partial} r
            \Big]
        \p
    \end{align}
    Using the assumptions on $q$ together with
    \[
        \sigma'_{R, L} \leq \csig L^{-1},\quad \sigma''_{R, L} \leq \csig L^{-2}
        \quad\text{and}\quad i \partial r \wedge \bar{\partial} r, i r \partial\bar{\partial} r \leq \omega_C
    \]
    we obtain a lower bound
    \begin{align}
        i \partial\bar{\partial} \chi_t
        &\geq i \partial\bar{\partial} \tilde{\chi}_t
            - \csig D t(1-t) \left(
                \colorA{2 R L^{-2} r} +
                \colorB{5 R L^{-1}}
            \right) \omega_C \\
        &\geq \kappa \omega_C - \frac54 \csig D \left( R^2 L^{-2} + R L^{-1}\right) \omega_C
        \p
    \end{align}
    So picking $L$ sufficiently large, we can ensure $i \partial\bar{\partial}\chi_t > 0$ on the roll-off region, for every $t \in [0, 1]$.

    \proofstep{Step 1.3: Positivity everywhere.}%
    
    So far, we have shown that for any choice of $B$ and any $t \in [0,1]$:
    \begin{itemize}
        \item $\iddball \chi_t > 0$ on $\{R + L < r\}$ for some large $R$ and $L$, \emph{independent of $B$}, and
        \item $i \partial\bar{\partial} \chi_t > 0$
    \end{itemize}
    Let $K = B_{R + L + 1} \times [0,1]$.
    By compactness, $\exists c > 0$ such that $i \partial\bar{\partial} \chi_t|_K > c i \partial\bar{\partial} \phi_0|_K$ for any $t$, so
    \[
         \frac{d^2}{d \tau^2} \chi_t \big|_K  - \left\lvert \partial \frac{d}{d \tau} \chi_t \big|_K \right\rvert_{i \partial\bar{\partial} \chi_t}
            \geq  \frac{d^2}{d \tau^2} \chi_t \big|_K  - c^{-1} \left\lvert \partial \frac{d}{d \tau} \chi_t \big|_K \right\rvert_{i \partial\bar{\partial} \phi_0}
        \p
    \]
    Looking at \cref{eq:ddbchi,eq:dmixchi} we see that the first term on the RHS increases linearly in $B$ while the second one is independent of $B$, so this expression is positive for $B$ sufficiently large.
    But as we remarked before, that is equivalent to $\iddball \chi_t > 0$ on $K$, so we have $\iddball \chi_t > 0$ everywhere.

    \proofstep{Step 2: The envelope has the required growth.}%

    From \cref{eq:barrier0} it is clear hat $\tilde{\chi}_t$ and hence also $\chi_t$ are asymptotic to the cone metric $\frac{i}{2} i \partial\bar{\partial} r_t^2$, so $\chi \in \cB$. Furthermore they also both agree with $\phi_0, \phi_1$ at $t = 0, 1$, so $\chi$ appears in the set $\cC$ from \cref{eq:perron}.

    This immediately shows that the $\sup$ in \cref{eq:perron} will not be $-\infty$ anywhere.
    Since the max of two quasi-PSh functions is again quasi-PSh, it suffices to consider $\theta \geq \chi$ in \cref{eq:perron}. But such a $\theta$ needs to be convex in $t$ (for a fixed $p \in M$) and has to agree with $\phi_{0,1}$ at the boundary so
    \begin{equation} \label{eq:bounds}
        \partial_t^+ \chi_0 \leq \partial_t^+ \theta_0 \leq \partial^\pm_t \theta_t \leq \partial_t^- \theta_1 \leq \partial_t^- \chi_1
    \end{equation}
    Since $\chi$ is smooth, this shows that $\{\theta \in \cC \mid \theta \geq \chi\}$ is equi-locally Lipschitz%
    \footnote{On a compact set in $M$, we have the same Lipschitz constant for all $\psi \in \cC$}
    in $t$.
    Hence $\Phi$, as defined in \cref{eq:perron}, is the upper-semicontinuous regularization of the supremum of a locally bounded family of PSh functions, and therefore is PSh. As it is obtained from metrics that do not depend on the imaginary direction of $\Omega$, the metric $\Phi$ will also be constant along the imaginary direction of $\Omega$. Hence $\Phi$ is convex in $t$. Since $\Phi \in \cC$ the inequalities \cref{eq:bounds} also apply to $\psi = \Phi$, so $\Phi$ is locally Lipschitz in $t$.
    We can explicitly bound the Lipschitz constant for large $r$ by looking at \cref{eq:barrier0} (recall that $\chi_t = \tilde{\chi}_t$ for large $r$):
    \[
        \partial_t^\pm \Phi_t \leq \partial_t \tilde{\chi}_1 \leq \frac{r_1^2}{2} (A + \dot{\theta}_t) + (\rho_1 - \rho_0) + 4B = O(r^2)
    \]
    and similarly for the lower bound, so $|\partial_t^\pm \Phi| = O(r^2)$. Since $\phi_{0, 1} = O(r^2)$ this shows that $\Phi = O(r^2)$.
    Finally $\Phi \geq \chi = \Theta(r^2)$ gives us the necessary lower bound to conclude $\Phi = \Theta(r^2)$.

    \proofstep{Step 3: $\Phi$ solves the homogeneous \MA/ equation.}%
    This part is standard. For completeness, we will reproduce the classical argument, taken out of \cite[Thm.~5.16]{GZ}.
    The claim is a local statement. Let $p \in \Omega \times M$, let $B$ be a ball around $p$ contained in $\Omega \times M$ and define $g = \Phi|_{\partial B}$, which is an upper-semicontinuous function.
    We can "localize" the definition of $\Phi$ in the following way: We claim that $\Phi|_B = v$, where
    \begin{equation} \label{eq:locsup}
        v = \left(\sup \left\{ w ~\middle\vert~ w \in \PSH(B),\: \limsup_{z \to \partial B} w(z) \leq g \right\}\right)^*
        \p
    \end{equation}
    Clearly $\Phi|_B$ is one of the competitors in the sup, so $v \geq \Phi|_B$ and in particular $\lim_{z \to \partial B} v(z) = g$. This implies that the metric
    \[
        \Psi(z) = \begin{cases}
            \max\{v(p), \Phi(p)\} & p \in B \\
            \Phi(p) & p \not\in B
        \end{cases}
    \]
    still is upper-semicontinuous. It also satisfies the mean-value inequality at every point for sufficiently small radii, so it is PSh.
    Since $\Psi$ still has the same values on $\partial \Omega \times M$ as $\Phi$, we conclude $\Phi = \Psi$ and hence $\Phi|_B = v$ by the maximality of $\Phi$. So it suffices to show that $v$ solves the homogeneous \MA/ equation in $B$, which is just a statement about pluripotential theory on the unit ball:

    The boundary data $g \in L^\infty(\partial B)$ can be approximated pointwise (using convolutions) by a decreasing sequence of $C^2$ functions $g_j$. By \cite[Thm.~5.13, 5.14]{GZ} the maximal functions
    \[
        u_j = \left(\sup \left\{ w ~\middle\vert~ w \in \PSH(B),\: \limsup_{z \to \partial B} w(z) \leq g_j \right\}\right)^*
    \]
    lie in $C^0\left(\bar{B}\right) \cap \PSH(B) \cap C^{1,1}(B)$ and satisfy $(i \partial\bar{\partial} u_j)^{n+1} = 0$ and $u_j|_{\partial B} = g_j$. Clearly $v \leq u_j \leq u_{j-1}$ and actually $u_j \searrow v$, because the (decreasing) limit of $u_j$ is PSh \cite[I.5.4]{agbook} and its $\limsup$ at the boundary is $\leq g_j$ for all $j$, so it is a competitor in \cref{eq:locsup}. Since the \MA/ operator is continuous along decreasing limits we obtain $(i \partial\bar{\partial} v)^{n+1} = 0$.

    This concludes the proof of \cref{thm:geod_exist}.
\end{proof}

\section{Energy and Ding functionals for singular AC metrics}%
\label{sec:functionals}

The goal of this section is to define the weighted Energy functional of a metric in $\Hsing^T$. This is straightforward for smooth, or even just $C^{1,1}$, metrics.
But since
we do not have any regularity in \cref{thm:geod_exist}
we need to do some extra work.
The main technical ingredient is defining the $X$-weighted \MA-measure of a metric $\phi \in \Hsing^T$. For smooth $\phi$ this is just
\[
    \MA_X(\phi) \coloneq e^{-h^\phi} (i \partial\bar{\partial} \phi)^n
    \c
\]
where $h^\phi$ is a suitably normalized Hamiltonian potential for $X$ w.r.t.\ the symplectic form $i \partial\bar{\partial} \phi$, see \cref{thm:hproperties} below.
This direct approach does not work for singular $\phi$ because $h^\phi$ cannot be defined unambiguously on a set of full $(i \partial\bar{\partial} \phi)^n$-measure. Instead we will follow the approach of \cite[Sec.\ 2.2]{BWN} and define certain projection operators $P_\lambda : \Hsing^T \to \PSH(M; K_M^{-1})^T$. For smooth $\phi$ these have the property that
\[
    (i \partial\bar{\partial} P_\lambda \phi)^n = \ind_{\{h^\phi < \lambda\}} (i \partial\bar{\partial} \phi)^n
\]
(see \cref{thm:mavsmooth}), allowing us to define $\MA_X(\phi)$ as
\[
    \MA_X(\phi) = \int_{-\infty}^\infty e^{-\lambda} (i \partial\bar{\partial} P_\lambda \phi)^n \d \lambda
    \p
\]
Since the RHS does not explicitly depend on $h^\phi$ this definition will work for singular $\phi \in \Hsing^T$.
To avoid technicalities, we will actually define $\MA_X(\phi)$ using an ad-hoc Riemann integral (cf.\ \cref{def:MAX}), but we could also use a Bochner integral as sketched here.

\begin{lem} \label{thm:hproperties}
    Let $\phi \in \cH_2^T$. The function $h^\phi = - \frac12 e^\phi \cL_X e^{-\phi}$ has the following properties:%
    \footnote{Note that since $e^{-\phi}$ is $T$-, so in particular $JX$-invariant, $h^\phi = - e^\phi \cL_{X^{1,0}} e^{-\phi}$, which is Berndtsson's definition}

    \begin{thmlist}
        \item If $\phi$ is smooth then $h^\phi$ is a Hamiltonian for $J X$, which can be expressed by any of the following equivalent equations:
            \begin{equation} \label{eq:Hamilton}
            \begin{aligned}
                - JX \iprod i \partial\bar{\partial} \phi &= \d h^\phi \\
                X^{1,0} \iprod \partial\bar{\partial} \phi &= \bar{\partial} h^\phi \\
                X \iprod i \partial\bar{\partial} \phi &= \d^c h^\phi \\[-0.3em]
                X &= \grad^{i \partial\bar{\partial} \phi} h^\phi
            \end{aligned}
            \end{equation}
            where $d^c = -i (\partial - \bar{\partial})$.
            As a consequence, $h^\phi$ is Morse--Bott and if $(M, X, i \partial\bar{\partial} \phi)$ is a shrinker, $h^\phi$ is the soliton potential (up to a constant).
        \item $h^{\phi + u} = h^\phi + \frac12 X(u)$.
        \item If $X$ vanishes at $z \in M$ then $h^\phi(z) = - \frac12 \tr(D_z X)$.%
    \end{thmlist}
\end{lem}
\begin{notation}
    Here, and everywhere else, $\tr(D_z X)$ denotes the "real" trace of the \emph{real} linear map $D_z X: T_z M \to T_z M$.
\end{notation}
\begin{proof}
    \point1
    Let $\Omega$ be a local holomorphic volume form and set $f = - \log (e^{-\phi} / i^{n^2} \Omega \wedge \bar{\Omega})$. Note that
    $i \partial\bar{\partial} \phi = i \partial\bar{\partial} f$. 
    Now compute
    \[
        h^\phi = -e^\phi \cL_{X^{1,0}} \left(e^{-f} i^{n^2} \Omega \wedge \bar{\Omega}\right) = X^{1,0} f - i^{n^2} \frac{\partial X^{1,0} \iprod \Omega}{\Omega}
        \p
    \]
    Using that the second term on the RHS is holomorphic we get
    \[
        \bar{\partial} h^\phi = \bar{\partial} (X^{1,0} f) = X^{1,0} \iprod \partial\bar{\partial} f
        \p
    \]
    Using that $h^\phi$ is real by definition we now get
    \begin{align}
        d h^\phi
    &= X^{1,0} \iprod \partial\bar{\partial} f - X^{0,1} \iprod \partial\bar{\partial} f
    = - i JX \iprod \partial\bar{\partial} f \\
    d^c h^\phi
    &= i(X^{1,0} \iprod \partial\bar{\partial} f + X^{0,1} \iprod \partial\bar{\partial} f)
    = i X \iprod \partial\bar{\partial} f
    \end{align}

    The computations above show that $\xi \in \Lie T \mapsto \mu_\xi \coloneq -e^\phi \cL_{-J \xi} e^{-\phi}$ is a moment map for the action of $T$ on $(M, \omega)$. Thus $h^\phi = \mu_{J X}$ is Morse--Bott since the moment maps of torus actions on symplectic manifolds are always (vector-valued) Morse--Bott functions (see for example \cite[Thm.~3.37]{NicolaescuMorse}).

    \point2 Follows immediately from the product rule.

    \point3
    Writing $e^{-\phi} = e^{-f} \Upsilon$ for an arbitrary smooth volume form $\Upsilon$, we have
    \[
        (\cL_X e^{-\phi})(z) = e^{-f(z)} (\cL_X \Upsilon)(z)
    \]
    because $(X f)(z) = 0$. This shows
    \[
        h^\phi(z) = - \Upsilon(z)^{-1} (\cL_X \Upsilon)(z) = -(\div_\Upsilon X)(z) = -\tr (D_z X)
        \p
        \myqedhere
    \]
\end{proof}

Next we shall recall some basic properties about the flows of "quasi-periodic" holomorphic vector fields (meaning the closure of their flow is a torus).
First recall that compact group actions can be linearized around fixed points (cf.\ \cite[proof of Prop.\ I]{CScstar}):
\begin{lem} \label{thm:linCoordBase}
    Let $M$ be a complex manifold with a holomorphic action $G \acts M$ by a compact Lie group that has a fixed point $p \in M$. Then around $p$ there exists a holomorphic coordinate system in which the action of $G$ is linear. In particular, the fixed point set of $T$ is a disjoint union of complex submanifolds.
\end{lem}
\begin{proof}
    Pick an arbitrary holomorphic coordinate system $z^1, \dots, z^n: U \to \bbC$ centered at $p$. Shrinking $U$ if necessary we can assume it is $G$-invariant. Since $p$ is a fixed point, $G$ acts linearly on $T_p^* M$. So changing coordinates by some $(\Lambda_{ij}) \in \GL(n; \bbC)$ we get $\xi_i = \sum_j \Lambda_{ij} \d_p z_j \in T^*_p M$ which are weight vectors for $G$, i.e.\ $g \cdot \xi_i = \alpha_i(g) \xi_i$ for some homomorphisms $\alpha_i: G \to U(1)$. Now define
    \[
        w_i(x) \coloneq \sum_j \int_G \alpha_i(g) \Lambda_{ij} z_j(g \cdot x) \d g = \sum_j \int_G \alpha_i(g) \Lambda_{ij} (g^{-1} \cdot z_j)(x) \d g
        \p
    \]
    From the definition of $\alpha_i$ we see $\d_p w_i = \xi_i$, so the $w_i$ are a coordinate system around $p$ since their differentials are linearly independent at $p$.
    Furthermore $g \cdot w_i = \alpha_i(g) w_i$, so the $G$-action on $M$ is linear and diagonal in the coordinates $w_i$.
\end{proof}

Applying this to the torus generated by $J X$ we obtain \cite[Prop.\ 6]{BryantSoliton}:
\begin{prop}[linearizing coordinates] \label{thm:linCoord}
    If $p \in M$ is a zero of $X$ then there is a holomorphic coordinate system around $p$ in which $X: M \to TM$ is linear. In particular the fixed point set of $X$ is a disjoint union of connected complex submanifolds.
\end{prop}

\begin{notation}
    We write $\gamma_t: M \to M$ for the time $t$ flow of the vector field $\frac12 X$.
\end{notation}

\def\BB/{Bia\l{}ynicki-Birula}
We have the "\BB/" decomposition of $M $:
\begin{prop}
    Let $Z_0, \dots,  Z_k \subseteq M$ be the connected components of the vanishing set of $X$, or equivalently, the fixed point set of $T$. For any point $p \in M$ the limit $\lim_{t \to \infty} \gamma_{-t}(p)$ exists and lies in the fixed point set of $X$. We define
    \[
        Z_\alpha^- \coloneq \left\{p \in M ~\middle|~ \lim_{t \to \infty} \gamma_{-t}(p) \in Z_\alpha\right\}
    \]
    to be the \demph{attracting set} of $Z_\alpha$.
    Then $M = \bigcup_\alpha Z_\alpha^-$.
\end{prop}

The existence of an accumulation point of $\gamma_{-t}(p)$ as $t \to -\infty$ is easy to see because $r$ is decreasing along the flow of $X$ and $E$ is compact. The uniqueness of the limit can then be seen by considering linearizing coordinates around an accumulation point.

\begin{defn} \label{def:lambda_a}
    Pick an arbitrary $\phi \in \cH_2$.
    On each of the connected sets $Z_\alpha$ the function $h^\phi|_{Z_\alpha} = - \frac{1}{2} \tr(D X|_{Z_\alpha})$ is a constant, which we call $\lambda_{Z_\alpha}$. We call $\Sigma \coloneq \{\lambda_{Z_\alpha}\}_\alpha$ the set of \demph{critical values}.
    We assume w.l.o.g.\ $\lambda_{Z_0} = \min \Sigma$
    and write $\lambda_0 \coloneqq \lambda_{Z_0}$.
    Note that
    \[
        \lambda_0 \coloneqq \min_E h^\phi = \min_M h^\phi
        \p
    \]
    Finally we define
    \begin{equation} \label{eq:defAl}
        A_\lambda \coloneqq \bigcup_{\lambda_{Z_\alpha} \geq \lambda} Z_\alpha^-
        \p
    \end{equation}
\end{defn}

For future reference we note
\thminflow
\begin{lem}
    \begin{thmlist}
        \item $E$ is connected.
        \item If $p$ maximizes $h^\phi$ on $E$ for some $\phi \in \cH_2$, then $X(p) = 0$. This implies $\max \Sigma = \max_E h^\phi$.
    \end{thmlist}
\end{lem}
\begin{proof}
    \point1 See \cite[Claim~3.13]{CDS}.

    \point2 The flow of $X$ preserves $E$, so $h^\phi(\gamma_t(p)) \leq h^\phi(p)$. But $X$ is the gradient flow of $h^\phi$, so $h^\phi(\gamma_t(p))$ increases as $t \to \infty$ unless $X(p) = 0$.
\end{proof}

The sets $Z_\alpha$ might not necessarily form a cell-decomposition, meaning that the closure $Z_\alpha$ might not contain every $Z_\beta$ it intersects. We can however use a Hamiltonian function $h^\phi$ to restrict which attracting sets it might intersect:

\begin{lem} \label{thm:Zupwards}
    We have $\displaystyle{\overline{Z_\alpha^-}  \subseteq \smashop{\bigcup_{\lambda_{Z_\beta} \geq \lambda_{Z_\alpha}}} Z_\beta^-}$. In particular, the sets $A_\lambda$ are closed.
\end{lem}
\begin{proof}
    Pick $\phi \in \cH_2$. Since $h^\phi$ is decreasing along the flow of $X$, we have $h^\phi(Z_\alpha^-) \subseteq [\lambda_{Z_\alpha}, \infty)$ and by continuity
    $h^\phi\Big(\overline{Z_\alpha^-}\Big) \subseteq [\lambda_{Z_\alpha}, \infty)$.
    Since $\overline{Z_\alpha^-}$ is $T$-invariant and hence a union of flow lines of $-X$, the limit of $h^\phi$ along all those flow lines must be $\geq \lambda_{Z_\alpha}$. The claim follows.
\end{proof}

\subsection{The "partial Legendre transform"}%
\label{sub:the_partial_legendre_transform_}

As sketched at the beginning of \cref{sec:functionals} we now define the projection operators $P_\lambda$ that we will use to define $\MA_X$ for singular hermitian metrics.
\begin{defn} \label{def:Plambda}
    For $\lambda > \lambda_0$ and $\phi \in \bar{\cH}^T_2$ we define $P_\lambda \phi$ as the pointwise infimum
    \begin{equation} \label{eq:defPl}
        P_\lambda \phi = \inf_{t \geq 0} \left(\gamma_{-t}^* \phi + \lambda t\right)
    \end{equation}
    or, in multiplicative notation
    \[
        e^{-P_\lambda \phi} = \sup_{t \geq 0} e^{-\lambda t} \gamma_{-t}^* e^{-\phi}
        \p
    \]
\end{defn}

\begin{prop} \label{thm:Plambda}
    The $P_\lambda \phi$ defined above is a non-negative singular hermitian metric and $T$-invariant, so $P_\lambda : \Hsing^T \to  \PSH(M; K_M^{-1})^T$.
    Furthermore $P_\lambda \phi$ is monotone in $\phi$ and $\lambda$ and $P_\lambda \phi \leq \phi$.
\end{prop}

Before we prove the assertion that $P_\lambda \phi$ is PSh, we want to establish
\thminflow
\begin{lem} \label{thm:flow}
    Fix a point $p \in M$.
    \begin{thmlist}
        \item The function $t \mapsto (\gamma_{-t}^* \phi) (p)$ is convex.
        \item The slope at $t \to \infty$ of $t \mapsto (\gamma_{-t}^* \phi)(p)$ is equal to $\frac12 \tr (D_z X)$, where $z = \lim_{t \to \infty} \gamma_{-t}(p)$.
        \item If $\phi$ is smooth then
            \begin{equation} \label{eq:hagain}
                \frac{d}{d t} \left(\gamma_{-t}^* \phi\right)(p) = - h^\phi(\gamma_{-t}(p))
                \p
            \end{equation}
    \end{thmlist}
\end{lem}
\begin{proof}
    \point1.
    Consider the pullback $a^* \phi \in \PSH(\bbC \times M; \pr_2^* K_M^{-1})$ along the complex flow
    \begin{equation} \label{eq:actionMap}
        a: \bbC \times M \to M,\quad (\tau, p) \mapsto \gamma_{-\tau}(p) \coloneq \gamma_{-\Re \tau}(p)
        \p
    \end{equation}
    Since $\phi$ is $JX$-invariant, $a^* \phi$ is a PSh function that does not depend on the imaginary coordinate of $\tau$, hence
    \[
        t \mapsto (a^* \phi)(-t, p) = (\gamma_{-t}^* \phi)(p)
    \]
    is convex.

    \point2.
    Let $U$ be a neighborhood of $z$ in which we can write $e^{-\phi} = e^{-f}\Upsilon$ with some volume form $\Upsilon$ and pick $T \gg 1$ so that $\gamma_{-T}(p) \in U$. For $t > T$ we have
    \[
        \frac{d}{d t}(\gamma_{-t}^* \phi)(p)
            = -\frac{\cL_{-X/2} \Upsilon}{\Upsilon}(\gamma_{-t}(p)) - \frac12 X(f) (\gamma_{-t}(p))
            \xrightarrow{t \to \infty} \frac{\cL_X \Upsilon}{2 \Upsilon}(z) - \frac12 (X f)(z)
        \p
    \]
    The first summand of the limit is equal to $(\div_\Upsilon X) (z) = \frac12\tr(D_z X)$. Since $X(z) = 0$, the second summand of the limit vanishes.

    \point3 is the result of a direct computation:
    \[
        \frac{d}{d t} (\gamma_{-t}^* \phi) = - \frac{\frac{d}{d t} \gamma_{-t}^* e^{-\phi}}{\gamma_{-t}^* e^{-\phi}} = \frac{\gamma_{-t}^* \cL_{X/2} e^{-\phi}}{\gamma_{-t}^* e^{-\phi}} = - \gamma_{-t}^* h^\phi
        \p
        \myqedhere
    \]
\end{proof}

\begin{proof}[Proof of \cref{thm:Plambda}]
    The assertions in the second sentence are immediate from the definition.

    Showing that $P_\lambda \phi$ is PSh is a straightforward application of Kiselman's minimum principle \cite[0.7.B]{agbook}. Consider the pullback $a^* \phi \in \PSH(\bbC \times M; \pr_2^* K_M^{-1})$ along the complex flow from \cref{eq:actionMap} and set
    \[
        \psi = a^* \phi + \lambda \Re \tau
        \c
    \]
    where $\tau$ is the coordinate on $\bbC$.
    If we have charts $U_\alpha$ for $M$ then $\bbC \times U_\alpha$ are charts for for $\bbC \times M$, in which we can consider $\phi$ and hence $\psi$ to be functions.
    Since $\phi$ is invariant under the flow of $J X$, the metric $\psi$ is independent of $\Im \tau$ and Kiselman's minimum principle applies%
    \footnote{Note that whether we take the infimum over $\Re \tau > 0$ or $\Re \tau \geq 0$ does not matter because the function in question is upper-semicontinuous.}
    to the PSh functions $\psi|_{(\bbR_{>0} \times i \bbR) \times U_\alpha}$, showing that $P_\lambda \phi$ is PSh or identically $-\infty$.

    For any $p \in A_0$ (as defined in \cref{eq:defAl}) the slope at $t \to \infty$ of $\psi(p, \infty)$ is $\lambda - \lambda_0 > 0$ by \cref{thm:flow}, so $P_\lambda \phi(p) > -\infty$. In conclusion $P_\lambda \phi \neq -\infty$ somewhere, so $P_\lambda \phi$ is PSh.
\end{proof}

\begin{prop} \label{thm:Plaffine}
    For any $a, b > 0$ there exists $R_0 > 0$ with the following property: If we define $R(\lambda) = R_0 \sqrt{\max\{\lambda, 1\}}$ , then for any $\lambda > \lambda_0$ and $\phi \in \Hsing^T(a, b)$, the metric $(P_\lambda \phi)|_{M \setminus \bar{B}_{R(\lambda)}}$ is affine in $x_1$, with slope $2 \lambda$ in equivariant charts.

    \label{thm:dr}
    Moreover, if $\phi$ is smooth then there exists $\epsilon = \epsilon(a, b) > 0$ such that $\epsilon r^2 \leq h^\phi \leq \epsilon^{-1} r^2$ for $r > 1$. Since $\inf_M h^\phi = \lambda_0$ is independent of $\phi$, we can even assume $h^\phi \geq \epsilon r^2 - \epsilon^{-1}$.
\end{prop}
\begin{proof}
    We work in an equivariant coordinate system with canonical holomorphic volume form $\Omega = \d z_1 \wedge \dots \wedge \d z_n$ and write
    \[
        e^{-\phi} = e^{-f} i^{n^2} \Omega \wedge \bar{\Omega}
        \c
    \]
    so that $f$ is convex in $x \coloneqq x_1$.
    In the smooth case we have $h^\phi = \frac12 X f = \frac12 \frac{\partial f}{\partial x_1}$, since $\cL_X \Omega = 0$.

    The hypothesis $\phi \in \Hsing$ translates to
    \begin{equation} \label{eq:sqmain}
        c e^{2x_1} \leq f \leq C e^{2x_1}
    \end{equation}
    for $x_1 \gg_{a, b} 1$ and some $c = c(a, b), C = C(a, b) > 0$.
    Using convexity this implies that for any $t > 0$,
    \[
        c e^{2x_1} + t \partial_1^+ f \leq C e^{2(x_1 + t)}
        \p
    \]
    Specializing to $t = 1$ we get an upper bound $\partial_1^+ f \leq e^{2x_1} (C e^2 - c)$. To get an analogous lower bound we use that
    \[
        (-t) \partial_1^- f + c e^{2x_1} \leq C e^{2(x_1-t)}
    \]
    provided $x_1 - t \gg 1$. This gives us $\partial_1^- f \geq e^{2x_1} t^{-1} (c - C e^{-2t})$.
    Picking $t > \frac12 \log (2 C / c)$ we get
    \begin{equation} \label{eq:lowerquad}
        \partial_1^- f \geq \frac{c}{\log(2C /c)} e^{2x_1}
    \end{equation}
    for $x_1 \gg_{a, b} 1$.
    This proves the claim about $h^\phi$ in the smooth case.

    Furthermore, the lower bound \cref{eq:lowerquad} shows that for $\lambda > 0$ we have $\partial_1^- f > \lambda$ when
    \[
        \log r \simeq x_1 > \log \sqrt{\lambda} + \frac12 \log \frac{ \log(2 C / c)}{c}
    \]
    and $x_1 \gg_{a,b} 1$.
    After accounting for the ratio between $\log r$ and $x_1$ (which might be different on different charts, but is bounded by some constant) we obtain $R(\lambda) \sim \sqrt{\max\{\lambda, 1\}}$ such that
    \[
        \frac{d^+}{d t} (\gamma^*_{-t} \phi) \big|_{M \setminus B_{R(\lambda)}} = - \frac12 \partial_1^- f \big|_{M \setminus B_{R(\lambda)}} < - \lambda
        \p
    \]

    If $\frac{d^+}{d t} (\gamma^*_{-t} \phi + \lambda t)(p) < 0$ then by convexity
    \[
        P_\lambda \phi(p) = \inf_{t \in \bbR} \left( \gamma_{-t}^* \phi(p) + \lambda t \right)
    \]
    because values of $t < 0$ don't contribute to the $\inf$.
    For $s > 0$, the same is true when replacing $p$ by $\gamma_s(p)$, so
    \[
        P_\lambda \phi(\gamma_s(p)) = \lambda s + \inf_{t \in \bbR} \left( f(\gamma_{-(t-s)}(p)) + \lambda (t-s) \right) = \lambda s + P_\lambda \phi(p)
        \p
    \]
    This shows that $P_\lambda \phi$ is affine outside $B_{R(\lambda)}$.
\end{proof}

\begin{ex} \label{thm:exCn}
    Consider the Gaussian soliton on $\bbC^n$, so $e^{-\phi} = e^{-r^2/2} \Upsilon$ with the euclidean volume form $\Upsilon$. Then $X = r \partial_r$ is the Euler vector field and 
    \begin{align}
        \gamma_{-t}^* e^{-\phi} &\equiv_\Upsilon e^{-e^{-t} r^2/2} e^{-nt} \qquad \text{or equivalently} \\
        \gamma_{-t}^* \phi &\equiv_\Upsilon n t + e^{-t} r^2/2
        \c
    \end{align}
    where we use the notation $\equiv_\Upsilon$ to indicate that we are trivializing using $\Upsilon$.
    Using \cref{eq:hagain} we compute
    \[
        h^\phi = r^2/2 - n
        \p
    \]
    We also see that the $t$-derivative of the expression inside the first inf in \cref{eq:defPl} is
    \[
        \lambda + n - e^{-t} r^2/2
        \p
    \]
    Solving for when this is zero we obtain $t_* = 2 \log r - \log (2 \lambda + 2 n)$, so
    \[
        P_\lambda \phi \equiv_\Upsilon \begin{cases}
            r^2/2 & r^2/2 - n \leq \lambda \\
            (\lambda + n) (1 + 2 \log r - \log(2 \lambda + 2n))& r^2/2 - n \geq \lambda
        \end{cases}
        \p
    \]
    Since we did not work in equivariant coordinates, we do not get slope $2 \lambda$ in $\log r$.

    We now re-do the computation using the equivariant coordinates $z_1, \dots, z_n$ implicitly defined by
    \[
        w_1 = e^{z_1}, \quad w_i = e^{z_1} z_i,\; i \geq 2
        \c
    \]
    where $w_1, \dots, w_n$ are the standard coordinates on $\bbC^n$. The canonical volume form $\Omega_z \wedge \overline{\Omega_z}$ in these coordinates is determined by
    \[
        2^n i^{-n^2} \Upsilon = \Omega_w \wedge \overline{\Omega_w} = \d w_1 \wedge \dots \wedge \d w_n \wedge \overline{\d w_1 \wedge \dots \wedge \d w_n} = e^{2n \Re z_1} \Omega_z \wedge \overline{\Omega_z}
        \p
    \]
    Note that $r^2 = e^{2 \Re z_1} (1 + |z_2|^2 + \cdots |z_n|^2)$, so we can write
    \[
        i^{n^2} \Omega_z \wedge \overline{\Omega_z} = r^{-2n} \Upsilon f(z_1, \dots,  z_n)
    \]
    for an appropriate function $f$ and both $\Omega_z \wedge \overline{\Omega_z}$ and $r^{-2n} \Upsilon$ are $X$-invariant.

    In the trivialization of $K_M^{-1}$ given by the coordinates $z_i$ the metric $\phi$ defined before corresponds to the function $\frac{r^2}{2} - 2n \log r + \log f$ and
    \[
        \gamma_{-t}^* \phi \equiv_{\Omega_z} n t + e^{-t} r^2/2 - 2n \log r + \log f
        \p
    \]
    (Note that in this trivialization we can actually compute $\gamma_{-t}^* \phi$ by substituting $r \to e^{-t/2} r$ in the expression for $\phi$). The minimum of $\gamma_{-t}^* \phi$ still occurs at $t_* = 2 \log r - \log(2\lambda + 2n)$ hence
    \[
        P_\lambda \phi \equiv_{\Omega_z} \begin{cases}
            r^2/2 - 2n \log r + \log f                   & r^2/2 - n \leq \lambda \\
        (\lambda + n) (1 - \log(2\lambda + 2n)) + 2\lambda \log r + \log f & r^2/2 - n \geq \lambda
        \end{cases}
    \]
    and we see that the slope in $\log r$ is $2\lambda$.
\end{ex}

\def\BT/{Bedford--Taylor}

\begin{notation}
    We write $\MA$ for the \BT/-\MA/ operator for locally bounded PSh metrics. As is usual we will use the symbol "$\weakto$" to denote weak convergence of measures.
\end{notation}

The metrics $P_\lambda \phi$ are very helpful in our non-compact situation because \cref{thm:Plaffine} and the next lemma together show that $\MA(P_\lambda \phi)$ has compact support.
\begin{lem} \label{thm:affinePSH}
    Let $u$ be a PSh function on an open subset $V \times U \subseteq \bbC \times \bbC^{n-1}$ that is affine in $x_1$ and independent of $y_1$. Then $u$ is equal to $a x_1 + v(z')$ for some $a \in \bbR$ and $v \in \PSH(U)$.
    If $u \in L^\infty_\mathrm{loc}$, the \MA/ measure of $u$ is zero.
\end{lem}

The proof is very straightforward and probably known in pluripotential theory, but we will include it here since we were not able to find it in the literature.
\begin{proof}
    We can think of $i \partial\bar{\partial} u$ as a matrix with measure entries. By Fubini's theorem we have for any test function $\chi$
    \[
        \int_{V \times U} u(z_1, z') \frac{\partial^2 \chi}{\partial z_1 \bar{\partial} z_1} \d z_1 \d z'
            = \int_U \int_V u(z_1, z') \frac{\partial^2 \chi}{\partial z_1 \bar{\partial} z_1} \d z_1 \d z'
        \p
    \]
    Since $u(\cdot, z')$ is affine we can integrate the inner integral by parts twice, obtaining that the $(1,\bar{1})$-entry of the matrix $i \partial\bar{\partial} u$ is zero.

    Now suppose the $(1,\bar{\jmath})$-entry is nonzero, so there exists a function $f \in C^\infty_c(U, \bbC)$ such that
    \[
        c \coloneqq \int u \frac{\partial^2 f}{\partial z_1 \bar{\partial} z_j} > 0
        \p
    \]
    Let $K \subseteq V \times U$ be a compact set such that $\supp f \subseteq K^\circ$ and observe that
    \[
        \int_K (i \partial\bar{\partial} u) \left(\epsilon^{-1} f \partial_1 - \epsilon \partial_j, \epsilon^{-1} \bar{f} \bar{\partial}_1 - \epsilon \bar{\partial}_j\right)
                = -2 \epsilon c + \epsilon^2 i \int_K \frac{\partial^2 u}{\partial z_j \bar{\partial} z_j} < 0
    \]
    for $\epsilon \ll 1$. This contradicts the fact that $u$ is PSh, so $\partial_1 \partial_{\bar{\jmath}}
    u = 0$.
    Consequently $\bar{\partial} \frac{\partial u}{\partial z_1} = 0$, so $\frac{\partial u}{\partial z_1}$ is a real-valued holomorphic function, i.e.\ a real constant $a$. In other words,
    \[
        \frac{\partial}{\partial x_1}(u - a x_1) = \frac{\partial}{\partial y_1}(u - a x_1) = 0
        \p
    \]
    Now a standard result
    from distribution theory tells us that $u - a x_1 = \pr_U^* v$ for some $v \in \cD'(U)$, proving the first claim ($\pr_U : V \times U \to U$ denotes the projection function).

    If we expand the \BT/-definition
    \[
        \MA(\pr_U^* v + a x_1) = i \partial\bar{\partial} \left[ (\pr_U^* v + a x_1) i \partial\bar{\partial} \left((\pr_U^* v + a x_1) i \partial\bar{\partial} (\cdots) \right) \right ]
    \]
    by multilinearity, we can first note that all terms containing $i \partial\bar{\partial} (a x_1 (\cdots))$ vanish by the product rule for distributions because $i \partial\bar{\partial} x_1 = 0$.
    The only remaining term is $\pr_U^* (i \partial\bar{\partial} v)^n$, which vanishes because $\dim_{\bbC} U < n$.
\end{proof}

The metrics $P_\lambda \phi$ will be unbounded (unless the zero-set of $X$ only has a single connected component) and this is not a pathological phenomenon, but essential to the argument. However we know that singularities will only appear for some values of $\lambda$:
\begin{lem} \label{thm:PlLocBdd}
    Let $\lambda > \lambda_0$, $\phi \in \Hsing^T$. Then $P_\lambda \phi|_{M \setminus A_\lambda}$ is locally bounded. In particular, $P_\lambda \phi$ is locally bounded (on all of $M$) when $\lambda > \max \Sigma$, i.e.\ when $\lambda$ is larger than all critical values.
\end{lem}
\begin{proof}
    First observe that for any $\phi, \phi' \in \Hsing^T$,
    \[
        \gamma_{-t}^*(\phi - \phi') \leq \gamma_{-t}^* \phi + \lambda t + \sup_{s \geq 0} (-\gamma_{-s}^* \phi' - \lambda s) = (\gamma_{-t}^* \phi + \lambda t) - P_\lambda \phi'
    \]
    so $\inf_{t \geq 0} \gamma_{-t}^* (\phi - \phi') \leq P_\lambda \phi - P_\lambda \phi'$ and similarly $P_\lambda - P_\lambda \phi' \leq \sup_{t \geq 0} \gamma_{-t}^* (\phi - \phi')$.
    Since $\phi$ and $\phi'$ are locally bounded and $\gamma_{-t}$ preserves the compact sets $B_r$, this shows that $P_\lambda \phi - P_\lambda \phi'$ is always locally bounded. Therefore we may assume that $\phi$ is smooth for the remainder of the proof.

    The condition $p \not\in A_\lambda$ is equivalent to $\lim_{t \to \infty} h^\phi(\gamma_{-t}(p)) < \lambda$. So for such $p$ we can find $T > 0$ and, by continuity of $h^\phi$, a neighborhood $U$ of $p$ such that
    \[
        \lambda - h^\phi(\gamma_{-t}(U)) > 0
    \]
    for $t = T$ and hence by monotonicity also for $t > T$. By \cref{thm:flow} 
    we have for $q \in U$
    \begin{align}
        \gamma_{-t}^* \phi(q) + \lambda t &\geq \gamma_{-T}^* \phi(q) + \lambda T + (t - T) (\lambda - h^\phi (\gamma_{-T}(q)))  \\
                            &> \gamma_{-T}^* \phi(q) + \lambda T
        \p
    \end{align}
    This gives us a lower bound on $P_\lambda \phi(q)$ that holds for all $q$ near $p$.
\end{proof}

On the other hand, when $p \in Z_\alpha^-$ with $\lambda_{Z_\alpha} > \lambda$ then
\[
    P_\lambda \phi (p) = -\infty
\]
because the slope at $t \to \infty$ of $\gamma_{-t} \phi(p) + \lambda t$ is equal to $-\lambda_{Z_\alpha} + \lambda$.
When $\lambda \not\in \Sigma$ (cf.\ \cref{def:lambda_a}),
\[
    \bigcup_{\lambda_{Z_\alpha} > \lambda} Z_\alpha^- = \bigcup_{\lambda_{Z_\alpha} \geq \lambda} Z_\alpha^- = A_\lambda
\]
so we have a perfect dichotomy where $P_\lambda \phi$ is bounded on $M \setminus A_\lambda$ and $-\infty$ on $A_\lambda$. In other words, the unbounded locus of $P_\lambda \phi$ coincides with $\{P_\lambda \phi = -\infty\}$. In particular we see that the sets $A_\lambda$ are pluripolar.

As mentioned before, we wish to use the \MA/ measures of the metrics $P_\lambda \phi$. Since these are not locally bounded when $\lambda < \max \Sigma$, we cannot use the \BT/ definition of $\MA$. Instead we define
\begin{defn}
    $\MA(P_\lambda \phi) = i_* \MA(P_\lambda \phi |_{M \setminus A_\lambda})$ where $i_*$ denotes the pushforward of measures form the open set $M \setminus A_\lambda$.
\end{defn}
Once we show that these measures have (locally) finite mass, we can conclude that the non-pluripolar product $(i \partial\bar{\partial} P_\lambda \phi)^n$ as defined in \cite[Def.\ 1.1]{BEGZ} exists and agrees with our definition of $\MA(P_\lambda \phi)$. We do this in the next theorem, which is the key technical result needed to define the weighted Energy functional for singular metrics.

\begin{thm} \label{thm:MAlmonotone}
    If $\nu \geq \lambda$
    then $\int_M \MA(P_\nu \phi) \geq \int_M \MA(P_\lambda \phi')$ for any $\phi, \phi' \in \Hsing^T$.
    In particular
    \begin{equation} \label{eq:MAlindep}
        \int_M \MA(P_\lambda \phi) = \int_M \MA(P_\lambda \phi') < \infty
        \p
    \end{equation}
\end{thm}
\begin{proof}
    \textit{Case 1: $\lambda \not\in \Sigma$ and $\nu > \lambda$.~}
    By \cref{thm:PlLocBdd} and the discussion above, $P_\lambda \phi$ and $P_\nu \phi'$ are locally bounded away from $A_\lambda$, which is Zariski-closed by \cref{thm:Zupwards} and $X$-invariant.

    Consider an arbitrary smooth volume form $\Upsilon$ on $M$ and let $B \coloneqq B_R$ for some $R > R(\nu) > R(\lambda)$. Define a function $J: [0, \infty) \to \bbR_{> 0}$ by
    \[
        J(t) = \sup_{\bar{B}} \log \frac{\gamma_{-t}^* \Upsilon}{\Upsilon}
        \p
    \]
    Since $\bar{B}$ is compact, $J$ is continuous.
    Pick $C \gg 1$ such that $e^{-C} \Upsilon \leq e^{-\phi} \leq e^C \Upsilon$ and $|\phi - \phi'| \leq C$ on $B$.
    Thus for any $s \geq 0$
    \begin{align}
        \nu s + \gamma_{-s}^* \phi &\geq \lambda s + \gamma^*_{-s} \phi' + (\nu - \lambda) s - C \\
                         &\geq P_\lambda \phi' + (\nu - \lambda) s - C \\
        \nu s + \gamma_{-s}^* \phi &\geq \phi + \nu s - J(s) - 2C \\
                         &= P_\lambda \phi' + \nu s - J(s) - 3 C + (\phi' - P_\lambda \phi')
    \end{align}
    and we can combine these inequalities to
    \[
        P_\nu \phi - P_\lambda \phi' \geq -3C + \inf_{s \geq 0} \max\{\nu s - J(s) + (\phi' - P_\lambda \phi'), (\nu - \lambda) s\} =: F_{\nu, \lambda}(\phi' - P_\lambda \phi')
    \]
    which holds on all of $B$.
    The function $F_{\nu, \lambda} : \bbR_{\geq 0} \to \bbR$ defined above is increasing, and since $J$ is continuous, $F_{\nu, \lambda}(t) \to \infty$ as $t \to \infty$.

    Since $\phi$ is bounded on $B$ and $P_\lambda \phi'$ is upper-semicontinuous, any sublevel set $\{P_\lambda \phi' - \phi' < -K\} \cap B$ contains an open neighborhood of $A_\lambda \cap B$ inside $B$ (here we used that $\lambda \not\in \Sigma$ to deduce that $P_\lambda \phi \equiv -\infty$ on $A_\lambda$). On that sublevel set we have $P_\nu \phi - P_\lambda \phi' \geq F_{\nu, \lambda}(K)$. By inverting $F_{\nu, \lambda}$, we see that for any $k$
    \begin{equation} \label{eq:maxNbh}
        \max(P_\nu \phi - k, P_\lambda \phi') \equiv P_\nu \phi -k
    \end{equation}
    on some neighborhood $U_k$ of $A_\lambda \cap B$ inside $B$.

    Outside of $B_{R(\nu)}$, both $P_\nu \phi - k$ and $P_\lambda \phi'$ are affine in $\log r$ with slopes $\nu$ and $\lambda$ respectively (in equivariant charts). Since $\nu > \lambda$ we conclude that \cref{eq:maxNbh} remains true on the open set
    \[
        \gamma_{(0, \infty)}\left( \partial B_{R(\nu)} \cap U_k \right) \supseteq A_\lambda \setminus B_{R(\nu)}
        \p
    \]
    Since $U_k \supseteq L(P_\lambda \phi') \supseteq L(P_\nu \phi)$, we know that $P_\nu \phi$ is bounded from below on $\partial B_{R - \epsilon} \setminus U_k$, so $\inf_{\partial B_{R - \epsilon} \setminus U_k} (P_\nu \phi - P_\lambda \phi') > - \infty$. Furthermore $P_\nu \phi$ has a larger slope than $P_\lambda \phi'$ along the flow lines of $X$, so $P_\nu \phi - k > P_\lambda \phi'$ on
    \[
        \gamma_{(T, \infty)}(\partial B_{R(\nu)} \setminus U_k) = M \setminus \bar{B}_{e^T R(\nu)} \setminus \gamma_{(T, \infty)}\left(\partial B_{R(\nu)} \cap U_k\right)
        \qtext{for}
        T_k \gg 0
        \p
    \]
    Altogether, we see that $\max(P_\nu \phi - k, P_\lambda \phi') \equiv P_\nu \phi - k$ on the open set
    \[
        \left( M \setminus \bar{B}_{e^{T_k} R(\nu)} \right) \cup \gamma_{(0, \infty)}\left(\partial B_{R(\nu)} \cap U_k\right) \cup (U_k \cap B)
        \c
    \]
    whose complement is compact and a subset of $M \setminus A_\lambda$.

    Thus using that $\int_{A_\lambda} \MA(P_\nu \phi) = 0$ since $\MA(P_\nu \phi)$ is non-pluripolar and Stokes' theorem (for currents) applied on $M \setminus A_\lambda$ we obtain
    \[
        \int_M \MA(P_\nu \phi) = \int_{M \setminus A_\lambda} \MA(P_\nu \phi - k) = \int_{M \setminus A_\lambda} \MA(\max(P_\nu \phi - k, P_\lambda \phi'))
        \p
    \]
    Now the decreasing limit
    \[
        \max\left(P_\nu \phi - k, P_\lambda \phi'\right) \searrow P_\lambda \phi'
    \]
    as $k \to \infty$ implies weak convergence
    \[
        \MA(\max(P_\nu \phi - k, P_\lambda \phi')) \weakto \MA(P_\lambda \phi')
    \]
    by \BT/'s convergence theorem and hence
    \[
        \int_M \MA(P_\nu \phi) \geq \int_{M \setminus A_\lambda} \MA(P_\lambda \phi') \geq \int_M \MA(P_\lambda \phi')
        \p
    \]

    \medskip\noindent
    \textit{Case 2: $\lambda \in \Sigma$ or $\nu = \lambda$.~}
    Pick an increasing sequence of non-critical $\lambda_k \in (\lambda_0, \lambda) \setminus \Sigma$ converging to $\lambda$.

    \begin{claim*}
        $P_{\lambda_k} \phi \nearrow P_\lambda \phi$ pointwise on $M \setminus A_{\lambda_1}$.
    \end{claim*}
    \begin{pf}
        For $p \not\in A_{\lambda_1}$ we know that $(\gamma_{-t}^* \phi)(p)$ has slope at least $-\lambda_1$ as $t \to \infty$. Hence for $k \geq 2$ the infimum in the definition of $P_{\lambda_k} \phi (p)$ in \cref{eq:defPl} is achieved by some $t_*$. Then for $k \geq 2$ if suffices to take the infimum over $0 \leq t \leq t_*$, which shows that $(P_\lambda \phi - P_{\lambda_k}\phi)(p) \leq (\lambda - \lambda_k) t_*$, proving the claim.
    \end{pf}

    By the weak continuity of $\MA$ along increasing sequences \cite[Thm.~3.23]{GZ} we conclude that we have $\MA(P_{\lambda_k} \phi) \weakto \MA(P_\lambda \phi)$ on $M \setminus A_{\lambda_1}$ and hence
    \[
        \int_{M \setminus A_{\lambda_1}} \MA(P_\lambda \phi')
            \leq \liminf_{k \to \infty} \int_{M \setminus A_{\lambda_1}} \MA(P_{\lambda_k} \phi')
            \leq \int_{M \setminus A_{\lambda_1}} \MA(P_\nu \phi)
        \c
    \]
    where we use Case 1 for the second inequality. Since $\MA(P_\lambda \phi')$ and $\MA(P_\nu \phi)$ are non-pluripolar the integrals do not change if we integrate over all of $M$, hence the claim.
\end{proof}

Using the theorem above we can prove continuity of $\MA(P_\lambda -)$ along decreasing sequences.
\begin{lem} \label{thm:Plconv}
    If $\phi_j \in \Hsing^T$ decrease towards $\phi$ then $\MA(P_\lambda \phi_j) \weakto \MA(P_\lambda \phi)$ on $M$.
\end{lem}
\begin{proof}
    By the standard \BT/ convergence theorem, $\MA(P_\lambda \phi_j|_{M \setminus A_\lambda}) \weakto \MA(P_\lambda \phi|_{M \setminus A_\lambda})$ on $M \setminus A_{\lambda}$. By \cref{thm:MAlmonotone} the weak convergence has no "loss of mass" i.e.\
    \[
        \int_{M \setminus A_\lambda} \MA(P_\lambda \phi|_{M \setminus A_\lambda}) = \lim_{j \to \infty} \int_{M \setminus A_\lambda} \MA(P_\lambda \phi_j|_{M \setminus A_\lambda})
        \c
    \]
    so have convergence against bounded continuous functions on $M \setminus A_\lambda$. So the pushed-forward measures $\MA(P_\lambda \phi_j) = i_* \MA(P_\lambda \phi_j|_{M \setminus A_\lambda})$ converge weakly to $\MA(P_\lambda \phi) = i_* \MA(P_\lambda \phi|_{M \setminus A_\lambda})$.
\end{proof}

The following formula for $\MA(P_\lambda \phi)$ in the smooth case motivates the definition of the operators $P_\lambda$.
\begin{prop} \label{thm:mavsmooth}
    For every $\phi \in \cH_2^T$ and $\lambda > \lambda_0$ we have
    \begin{align}
        \{h^\phi \leq \lambda\} &= \{P_\lambda \phi = \phi\} \qquad\text{and} \\
        \MA(P_\lambda \phi) &= (i \partial\bar{\partial} \phi)^n \ind_{\{h^\phi < \lambda\}}
        \p
    \end{align}
\end{prop}
\begin{proof}
    By the convexity result from \cref{thm:flow}, $(P_\lambda \phi)(p) = \phi(p)$ if{}f $\ddat{t}{0} \gamma_{-t}^* \phi \geq -\lambda$. In view of \cref{eq:hagain} this is equivalent to $h^\phi \leq \lambda$, as desired.

    The second equation can be proven just as in \cite[Lem.\ 2.6]{BWN}.
    First note that
    \[
        \MA(\phi) \ind_{\{h^\phi < \lambda\}} = \MA(P_\lambda \phi) \ind_{\{h^\phi < \lambda\}}
        \c
    \]
    because $\{h^\phi < \lambda\}$ is an open set on which $\phi \equiv P_\lambda \phi$. Since $h^\phi$ is Morse--Bott its level sets have vanishing Lebesgue measure, and hence also vanishing $\MA(\phi)$-measure. Consequently
    \begin{equation} \label{eq:MA_sublevel}
        \MA(\phi) \ind_{\{h^\phi \leq \lambda\}} = \MA(\phi) \ind_{\{h^\phi < \lambda\}} \leq \MA(P_\lambda \phi)
    \end{equation}
    On the other hand, when $\lambda \not\in \Sigma$ the comparison principle on $M \setminus A_\lambda$ gives%
    \footnote{
        The comparison principle in \cite[Thm.~3.29]{GZ} does not immediately apply to our situation. But the proof still applies if we take $\Omega' = B_R \setminus A_\lambda$ for some large $R$ (we cannot take $R = \infty$ because we need $\int_{\Omega'} \MA(\phi) < \infty$).
        Since $\phi$ has a quadratic lower bound while $P_\lambda \phi = O(\log r)$, we have $\phi - P_\lambda \phi \to \infty$ as $r \to \infty$, and since $\lambda \not\in \Sigma$ we also have $A_\lambda = \{P_\lambda = -\infty\}$, so $\limsup_{p \to \partial (B_R \setminus A_\lambda)} \phi(p) - P_\lambda \phi(z) > 0$ for $R \gg 0$.
    }
    \begin{align}
        \int_{\{P_\lambda \phi > \phi - \epsilon\}} \MA(P_\lambda \phi)
            &\leq \int_{\{P_\lambda \phi > \phi - \epsilon\}} \MA(\phi)
            \p
    \end{align}
    Letting $\epsilon \to 0$ we
    see
    that \cref{eq:MA_sublevel} must actually be an equality.

    If $\lambda \in \Sigma$ we can argue as in the second case of the proof of \cref{thm:MAlmonotone} to obtain
    \[
        \int_M \MA(P_\lambda \phi) \leq \liminf_{k \to \infty} \int_M \MA(P_{\lambda_k} \phi) = \liminf_{k \to \infty} \int_{\{h^\phi \leq \lambda_k\}} \MA(\phi) = \int_{\{h^\phi \leq \lambda\}} \MA(\phi)
    \]
    which again implies that \cref{eq:MA_sublevel} is an equality.
\end{proof}

We can use the preceding theorem together with smooth approximation to upgrade \cref{thm:MAlmonotone} in the case where $\phi = \phi'$: We do not just have an inequality of total masses, but actually an inequality of measures.
\begin{cor} \label{thm:MAlMeasureMonotone}
    If $\phi \in \Hsing^T$ and $\nu \geq \lambda$ then $\MA(P_\nu \phi) \geq \MA(P_\lambda \phi)$.
\end{cor}
\begin{proof}
    Pick a smooth sequence of metrics $\cH_2^T \ni \phi_j \searrow \phi$ (this is possible by \cref{thm:regularization}) and observe that
    \[
        0 \leq \ind_{\{\lambda \leq h^{\phi_j} < \nu\}} (i \partial\bar{\partial} \phi_j)^n = \MA_\nu(\phi_j) - \MA_\lambda(\phi_j) \weakto \MA_\nu(\phi) - \MA_\lambda(\phi)
        \p
        \myqedhere
    \]
\end{proof}

We can interpret \cref{thm:MAlmonotone} as saying that a metric $\phi \in \Hsing^T$ defines a \DH/ measure for the action of $X$, and these measures are independent of $\phi$:
\begin{cor} \label{thm:DH}
    The \demph{\DH/ measure} $\nu = h^\psi_* (i \partial\bar{\partial} \psi)^n$ on $[\lambda_0, \infty)$ is independent of $\psi \in \cH_2^T$. In fact, its cumulative distribution function satisfies
    \[
        \int_{\lambda_0}^\lambda \nu = \int_M \MA(P_\lambda \phi)
    \]
    for any $\phi \in \Hsing^T$.
\end{cor}
\begin{proof}
    Observe that by \cref{thm:mavsmooth,thm:MAlmonotone}
    \[
        \int_{\lambda_0}^\lambda h^\psi_* (i \partial\bar{\partial} \psi)^n = \int_{h^\psi \leq \lambda} (i \partial\bar{\partial} \psi)^n = \int_M \MA(P_\lambda \psi) = \int_M \MA(P_\lambda \phi)
        \p
        \myqedhere
    \]
\end{proof}

\subsection{The $X$-\MA/ measure}%
\label{sub:the_v_monge_ampere_measure}

Now we are ready to define the $X$-weighted \MA/ measure $\MAX$, which is a generalization of $(i \partial\bar{\partial} \phi)^n e^{-h^\phi}$ to singular $\phi \in \bar{\cH}_2$.
We have seen that for smooth $\phi$,
\[
    \ind_{[\lambda_0, \lambda]}(h^\phi) (i \partial\bar{\partial} \phi)^n = \MA(P_\lambda \phi)
\]
and the RHS makes sense for singular $\phi$.
Following \cite[Section 2]{BWN} we define
\begin{equation} \label{eq:defMAint}
    \MA_{\ind_{[\lambda_0, \lambda]}}(\phi) = \MA(P_\lambda \phi)
    \p
\end{equation}
Extending by linearity we can now define measures $\MA_g(\phi)$ for singular $\phi$ and any step function $g$. For smooth $\phi$ this definition satisfies $\MA_g(P_\lambda \phi) = g(h^\phi) (i \partial\bar{\partial} \phi)^n$.
Our goal is to define $\MA_{e^{-x}}(\phi)$ for singular $\phi$, which we will do by approximating $e^{-x}$ by step functions.
To this end, let $g_k : \bbR_{\geq 0} \to \bbR_{\geq 0}$ be an increasing sequence of positive linear combinations of the form
\[
    g_k = \sum_{i=1}^{N_k} \alpha_{ki} \ind_{[\lambda_0, \lambda_{ki}]} \qquad \alpha_{ki} \geq 0
\]
that approximate $x \mapsto \exp(-x)$ uniformly on $[\lambda_0, \infty)$.

For technical reasons we will be using auxiliary metrics in $(\cH_2^*)^T$, so we will from now on assume that $(\cH_2^*)^T$ is nonempty. Recall that by \cref{thm:KRSasymptotics} this is the case when $(M, X)$ admits a shrinker with quadratic curvature decay.

\begin{defn-prop} \label{def:MAX}
    Let $\phi \in \Hsing^T(a, b)$
    \begin{thmlist}
        \item The limit
            \begin{equation} \label{eq:defMAX}
                \MAX(\phi) \coloneqq \lim_{k \to \infty} \MA_{g_k}(\phi)
            \end{equation}
            exists in the (Banach) space of bounded Radon measures on $M$, and the convergence is uniform in $\phi$.
        \item If $\phi \in \cH_2^T$ (i.e.\ $\phi$ is smooth) then $\MAX(\phi) = e^{-h^\phi} (i \partial\bar{\partial} \phi)^n$.

        \item $\int_M \MAX(\phi)$ is independent of $\phi$.
    \end{thmlist}
\end{defn-prop}
\begin{proof}
    Pick an arbitrary auxiliary $\psi \in (\cH_2^*)^T$.
    By \cref{thm:MAlMeasureMonotone} we know that $0 \leq \MA_{\ind_{[\lambda, \nu]}}(\phi)$ for any $\nu \geq \lambda \geq \lambda_0$. Therefore
    \begin{equation}
        \| \MA_{\ind_{[\lambda, \nu]}}(\phi) \| = \int_M \MA(P_\nu \phi) - \int_M \MA(P_\lambda \phi)
        \p
    \end{equation}
    By \cref{thm:MAlmonotone}, the total mass $\int_M \MA(P_\lambda \phi)$ is independent of $\phi$, so we can use the formula for smooth metrics from \cref{thm:mavsmooth} to obtain
    \begin{equation}
         \int_M \MA(P_\nu \phi) - \int_M \MA(P_\lambda \phi) = \int_M \MA(P_\nu \psi) - \int_M \MA(P_\lambda \psi)
         = \int_M \ind_{[\lambda, \nu]} (h^\psi) \MA(\psi)
         \p
    \end{equation}
    For $m \geq k$ this implies
    \[
        \left\| \MA_{g_m}(\phi) - \MA_{g_k}(\phi) \right\| = \int (g_m - g_k) \circ h^\psi \MA(\psi) 
            \leq \int \left(e^{-h^\psi} - g_k \circ h^\psi\right) \MA(\psi)
            \p
    \]
    The integrand on the right goes to zero pointwise, so if we establish that $\int \exp(-h^\psi) \MA(\psi) < \infty$ we obtain that $\MA_{g_m}(\phi)$ is Cauchy by the dominated convergence theorem, showing \point1.

    By the definition of $\cH_2^*$ in \cref{eq:defHs} there is a cone metric $\tilde{\omega} = i \partial\bar{\partial} \tilde{r}^2/2$ such that
    $| i \partial\bar{\partial} \psi - \frac{i}{2} \partial\bar{\partial} \tilde{r}^2 |_{\tilde{\omega}} = O\left(r^{-2}\right)$. Together with the lower bound $h^\psi \geq -\epsilon^{-1} + \epsilon r^2$ from \cref{thm:dr} this shows
    \begin{equation} \label{eq:EstimateGaussianInt}
        \int_M e^{-h^\psi} \MA(\psi) = \int_M e^{O(\log r)} e^{-\tilde{r}^2} \left(1 + O\left(r^{-2}\right)\right) \left( \frac{i}{2} \partial\bar{\partial} \tilde{r}^2 \right)^n < \infty
        \p
    \end{equation}

    When $\phi$ is smooth, we use \cref{thm:DH} to see that
    \begin{align}
        \left\| e^{-h^\phi} (i \partial\bar{\partial} \phi)^n - \MA_{g_k}(\phi) \right\| &= \int \left(e^{-h^\phi} - g_k(h^\phi)\right) (i \partial\bar{\partial} \phi)^n = \int \left( e^{-x} - g_k(x) \right) \nu(\d x) \\
                                                     &= \int \left(e^{-h^\psi} - g_k(h^\psi)\right) (i \partial\bar{\partial} \psi)^n \to 0
    \end{align}
    as above, thus establishing \point2.

    \point3 is an immediate consequence of \cref{thm:MAlmonotone}.
\end{proof}

Now we re-prove \cite[Thm.~2.7]{BWN} in our setting.
\begin{thm} \label{thm:psiMAphiConverge}
    Let $\psi_j, \phi_j, \theta_j \in \Hsing^T$ be non-increasing sequences converging pointwise to $\psi_\infty, \phi_\infty, \theta_\infty \in \Hsing^T$. Then we have weak convergence
    \[
        (\psi_j - \theta_j) \MAX(\phi_j) \weakto (\psi_\infty - \theta_j) \MAX(\phi_\infty)
    \]
    with no loss of mass, i.e.
    \begin{equation} \label{eq:psiMAphiConverge}
        \int (\psi_j - \theta_j) \MAX(\phi_j) \to \int (\psi_\infty - \theta_j) \MAX(\phi_\infty)
        \p
    \end{equation}
\end{thm}
\begin{proof}
    \begin{claim*}
        For a fixed $\lambda$, we have $(\psi_j - \theta_j) \MA(P_\lambda \phi_j) \weakto (\psi_\infty - \theta_\infty) \MA(P_\lambda \phi_\infty)$.
    \end{claim*}
    \begin{pf}
        This is a local statement so is suffices to prove this convergence on open subsets $\Omega \subseteq M$ that are biholomorphic to balls. In particular, we can treat $\phi_j$ and $P_\lambda \phi_j$ as functions.

        We will conclude this from the "capacity convergence theorem" \cite[Thm.~3.2]{BT} (cf.\ also \cite[Thms.~4.26, 4.21]{GZ}).
        The sequences $\psi_j, \theta_j$ are monotone, hence convergent in capacity \cite[Prop.\ 4.25]{GZ}. By the subadditivity of the MA-capacity, $\psi_j - \theta_j$ is also convergent in capacity.
        To apply the convergence theorem
        we are only left with verifying that the currents $\MA(P_\lambda \phi_j), j=1,\dots,\infty$ are uniformly bounded by capacity, meaning
        there exists a constant $C$ such that for all Borel $B \subseteq \Omega$
        \begin{equation} \label{eq:boundedByCap}
            \int_B \MA(P_\lambda \phi_j) \leq C \Cap(B; \Omega)
            \p
        \end{equation}

        Pick any $\nu > \lambda, \max(\Sigma)$. By \cref{thm:PlLocBdd} the metrics $P_\nu \phi_j $ are bounded,
        so first by the definition of the \MA/ capacity \cite[4.16]{GZ}, and then by monotonicity of $P_\nu \phi_j$ we obtain
        \begin{align}
            \int_B \MA(P_\nu \phi_j) &\leq (\mathop{\operatorname{osc}}_B P_\nu \phi_j) \Cap(B; \Omega) \\
                                &\leq \left(\sup_B P_\nu \phi_1 + \lvert \inf_B P_\nu \phi_\infty \rvert \right) \Cap(B; \Omega)
            \p
        \end{align}
        Since $\MA(P_\lambda \phi_j) \leq \MA(P_\nu \phi_j)$ by \cref{thm:MAlMeasureMonotone}, and $P_\nu \phi_j \in L^\infty_\mathrm{loc}$, we obtain the desired uniform bound.
    \end{pf}

    For the weak convergence claim we need to show that for any test function $\chi$ on $M$,
    \begin{equation} \label{eq:WriteOutWeakConv}
        \left\lvert \int_M \chi \big( (\psi_j - \theta_j) \MAX(\phi_j) - (\psi_\infty - \theta_\infty) \MAX(\phi_\infty) \big) \right\rvert \to 0
        \p
    \end{equation}
    Since $\MA_{g_k}(\phi_j) \xrightarrow{k \to \infty} \MAX(\phi_j)$ uniformly in $j$, and $\psi_j - \theta_j$ is uniformly bounded on $\supp \chi$, we can replace $\MAX$ by $\MA_{g_k}$. Then the statement follows from weak convergence of $\MA(P_\lambda -)$ (and linearity).

    \smallskip
    To prove \cref{eq:psiMAphiConverge}, it suffices to show that for any $\epsilon > 0$ there exists a compact $K \subseteq M$ such that
    \begin{equation} \label{eq:unifMAXboundAtInfty}
        \int_{M \setminus K} \left\lvert (\psi_j - \theta_j) \MAX(\phi_j) \right\rvert < \epsilon
    \end{equation}
    for all $j$. Then we can pick $\chi$ such that $\chi \equiv 1$ on $K$ and use the weak convergence we just proved.

    Since $\phi_1, \phi_\infty \in \Hsing^T$, there are some $a, b > 0$ such that $\phi_j \in \Hsing(a, b)$ for all $j$. We will show that if we pick $K = \bar{B}_{R(\mu)}$ for some $\mu \gg_{a,b} 0$, then \cref{eq:unifMAXboundAtInfty} holds for any $\psi_j, \theta_j, \phi_j \in \Hsing(a, b)^T$.

    First of all, $(\psi_j - \theta_j) \leq (a+b)r^2 + O(1)$.
    Next, recall that by \cref{thm:dr} $\supp \MA(P_\lambda \phi) \subseteq \bar{B}_{R(\lambda)}$ where $R(\lambda) = R_0 \sqrt{\max\{\lambda, 1\}}$ and $R_0$ only depends on $a, b$. Hence
    \[
        r^2 \MA(P_\lambda \phi) \leq R_0^2 \max\{\lambda, 1\} \MA(P_\lambda \phi)
        \p
    \]
    Furthermore \cref{thm:MAlMeasureMonotone} implies
    \[
        (e^{-\lambda_1} - e^{-\lambda_2}) \MA(P_{\lambda_1} \phi) \leq \int_{\lambda_1}^{\lambda_2} \d \lambda e^{-\lambda} \MA(P_\lambda \phi)
        \p
    \]
    Applying this inequality to any of the step functions $g_k(x) \leq e^{-x}$ we obtain for any $\mu > 1$
    \begin{equation} \label{eq:unifMeasureBound}
    \begin{aligned}
        \int_{M \setminus \bar{B}_{R(\mu)}} r^2 \MA_{g_k}(\phi)
        &\leq \int_{\lambda_0}^\infty e^{-\lambda} \int_{M \setminus \bar{B}_{R(\mu)}} r^2 \MA(P_\lambda \phi) \d \lambda \\
        &\leq \int_\mu^\infty e^{-\lambda} R_0^2 \lambda \left( \int_M \MA(P_\lambda \phi)  - \int_M \MA(P_\mu \phi)\right) \d \lambda \\
        &= R_0^2 \int_\mu^\infty e^{-\lambda} \lambda \nu([\mu, \lambda]) \d \lambda \\
        &= R_0^2 \int_\mu^\infty (1 + \lambda) e^{-\lambda} \, \nu(\d \lambda) \\
          &= R_0^2 \int_{\{h^\psi \geq \mu\}} (1 + h^\psi) e^{-h^\psi} (i \partial\bar{\partial} \psi)^n \\
          &= O\left(\mu^{1+n} e^{-\epsilon^2 \mu}\right)
    \end{aligned}
    \end{equation}
    In the second inequality step we have shrunk the domain of the $\lambda$-integration by using that $\supp \MA(P_\lambda \phi) \subseteq \bar{B}_{R(\mu)}$ for $\lambda \leq \mu$.
    The last integrand is bounded by a constant times $(1 + r^2) e^{-\epsilon r^2} (i \partial\bar{\partial} r^2)^n$ with the $\epsilon$ from \cref{thm:dr}.
    To obtain the precise asymptotic we use that $\epsilon^{-1} r^2 \geq h^\psi$.
    We see that $\mu \to \infty$ makes \cref{eq:unifMeasureBound} go to zero.
\end{proof}

\begin{rmk}
    With the same arguments as in \cref{eq:unifMeasureBound} we obtain for any $\phi \in \Hsing(a, b)^T$
    \begin{equation} \label{eq:gaussianpoly}
    \begin{aligned}
        \int_M r^{2k} \MA_X(\phi) &\leq \int_{\lambda_0}^\infty e^{-\lambda} \int_M R(\lambda)^{2k} \MA(P_\lambda \phi) \d \lambda \\
                               &\leq R_0^{2k} \int_{\lambda_0}^\infty e^{-\lambda} (\lambda - \lambda_0 + 1)^k \nu([\lambda_0, \lambda]) \d \lambda \\
                               &\lesssim_k R_0^{2k} \int_{\lambda_0}^\infty e^{-\lambda} (\lambda - \lambda_0 + 1)^k \, \nu(\d \lambda)
                               = R_0^{2k} \int_M (h^\psi - \lambda_0 + 1)^k e^{-h^\psi} (i \partial\bar{\partial} \psi)^n \\
                               &\lesssim_{a,b,k} \int_M (1 + r^2)^k e^{-\epsilon r^2} \left(i \partial\bar{\partial} r^2 + O(r^{-2})\right)^n \lesssim 1
                               \p
    \end{aligned}
    \end{equation}
    where we write "$\lesssim_{A, B}$" to denote "$\leq$" up to a constant that only depends on the quantities $A, B$, and possibly also on background data such as $M$ and $r$. Similarly,
    \begin{equation} \label{eq:boundMminusBall}
        \int_{M \setminus B_R} r^{2k} \MA_X(\phi) = O\left(R^{2k+2n} e^{-\epsilon^2 R^2/ R_0^2}\right)
        \p
    \end{equation}
    We will use these estimates over and over again to justify integrations by parts.
\end{rmk}

\subsection{Definition and convexity of the weighted \MA/ Energy}%

\begin{defn}
    On $\bar{\cH}_2^T$ we define an \demph{Energy} and \demph{$\cF$-}functional by
    \begin{align}
        \cE_X(\phi_1, \phi_0) &= \int_0^1 \int_M (\phi_1 - \phi_0) \MAX(t \phi_1 + (1-t) \phi_0) \d t
        \label{eq:defEX}
        \\
        \cF(\phi) &= -\log \int_M e^{-\phi}
        \label{eq:defF}
        \p
    \end{align}
    The \demph{Ding} functional is defined as
    \[
        D_X(\phi_1, \phi_0) = \vol_X(K_M^{-1})^{-1} \cE_X(\phi_1, \phi_0) - \cF(\phi_1)
        \c
    \]
    where $\vol_X(K_M^{-1}) = \int_{\lambda_0}^\infty e^{-\lambda} \nu(\d \lambda)$ is computed using the \DH/-measure from \cref{thm:DH}.
\end{defn}

From \cref{eq:gaussianpoly} we see that $\exists C = C(M, a, b)$ such that
\begin{equation} \label{eq:EXbound}
    \left| \int_M (\phi_1 - \phi_0) \MAX(t \phi_1 + (1-t) \phi_0) \right| \leq C
\end{equation}
for $\phi_1, \phi_2 \in \Hsing(a, b)^T$. Likewise, $\cF(\phi)$ is uniformly bounded for $\phi \in \Hsing(a, b)$.

\begin{prop} \label{thm:Econt}
    The functional $\cE_X$ is continuous along non-increasing sequences (in both arguments).
\end{prop}
\begin{proof}
    This is an easy consequence of \cref{thm:psiMAphiConverge}: Using \cref{thm:regularization}, pick a decreasing sequence $\phi^\nu_i \in \cH_2^T$ decreasing to $\phi_i$. We know that
    \[
        \int_M (\phi^\nu_1 - \phi^\nu_0) \MAX(t \phi^\nu + (1-t) \phi^\nu_0)
    \]
    converges for every $t$ as $\nu \to \infty$, and by picking $a, b$ such that $\phi_i^\nu \in \Hsing(a, b)$ and applying \cref{eq:EXbound} we see that it is uniformly bounded, so the dominated convergence theorem applies to the $t$-integral.
\end{proof}

In the case of smooth metrics, the weighted \MA/ energy is a primitive of the $X$-\MA/ measure, meaning that we can obtain $\cE_X$ by integrating $\MA_X$ along any curve in $\cH_2$:
\begin{prop}[{\cite[Prop.~10.2]{BoBM}}] \label{thm:dsmooth}
    If $\phi_{t,s}$ is a smooth 2-parameter family in $\cH_2^T(a, b)$ with the property that all derivatives (in $s$-, $t$- and $M$-directions) are uniformly $O(r^2)$ in equivariant charts then
    \begin{equation} \label{eq:dsmooth}
        \frac{d}{d s} \int_M \dot{\phi}_{t,s} \MAX(\phi_{t, s})
            = \int_M \left(
                \frac{\partial^2 \phi_{t,s}}{\partial t \partial s}
                - \left\langle \bar{\partial} \frac{\partial \phi_{t,s}}{\partial t}, \bar{\partial} \frac{\partial \phi_{t,s}}{\partial s} \right\rangle_{i \partial\bar{\partial} \phi_{t,s}}
            \right)
            \MAX(\phi_{t, s})
    \end{equation}
    (dot denotes $t$-derivative).

    If we apply this to $\phi_{t, s} = s(t \phi_1 + (1-t) \phi_0) + (1-s) \phi_t$, where $\phi_t$ is a smooth family with uniformly $O(r^2)$ derivatives, we obtain
    \begin{equation} \label{eq:Eprimitive}
        \cE_X(\phi_1, \phi_0) = \int_0^1 \int_M \dot{\phi}_t \MAX(\phi_t) \d t
        \p
    \end{equation}
\end{prop}
\begin{proof}
    The proof goes by differentiating under the integral sign and integrating by parts. We denote $s$-derivatives by prime and omit the subscripts for $\phi_{t, s}$.
    Note that by \cref{thm:dr} and our assumption of quadratic growth on the derivatives
    \[
        \frac{d}{d s} \left( \dot{\phi} \MA_X(\phi) \right) = \left(\dot{\phi}' i \partial\bar{\partial} \phi - \dot{\phi} (X \phi') i \partial\bar{\partial} \phi + n \dot{\phi} i \partial\bar{\partial} \phi' \right) \wedge (i \partial\bar{\partial} \phi)^{n-1} e^{-h^\phi} = O(r^{2n+4} e^{-\epsilon r^2})
    \]
    uniformly in $s, t$. Hence the $s$-derivative of the integrand in the LHS of \cref{eq:dsmooth} is integrable, justifying differentiation under the integral sign.
    Computing the $s$-derivative of the integrand of the LHS of \cref{eq:dsmooth} is a somewhat tedious computation on forms which is no different than that in \cite[Prop.~10.2]{BoBM}. We nevertheless include it here for future reference.
    \begin{align}
        \frac{d}{d s} \left[ \dot{\phi} (i \partial\bar{\partial} \phi)^n e^{-h^\phi} \right]
            &= \left[ \left(\dot{\phi}' - \dot{\phi} X(\phi')\right) i \partial\bar{\partial} \phi + n \dot{\phi} i \partial\bar{\partial} \phi' \right] \wedge (i \partial\bar{\partial} \phi)^{n-1} e^{-h^\phi} \\
            &= \left[
                \left(
                    \dot{\phi}'
                    \hl{
                        - \dot{\phi} X(\phi')
                    }
                \right)
                i \partial\bar{\partial} \phi
                - n (\partial\dot{\phi}
                    \hl{ - \dot{\phi} \partial h^\phi }
                ) \wedge i \bar{\partial} \phi'
            \right]
            \wedge (i \partial\bar{\partial} \phi)^{n-1} e^{-h^\phi} \\
                &\qquad + n \partial \left( \dot{\phi} i \bar{\partial} \phi' e^{-h^\phi} \right) \wedge (i \partial\bar{\partial} \phi)^{n-1}
    \end{align}
    At this point we use Hamilton's equation \cref{eq:Hamilton}, to replace $\partial h^\phi$ above by $-X^{0,1} \iprod \partial\bar{\partial} \phi$. We see that the \hl{terms involving $\dot{\phi}$ without derivatives} are equal to
    \[
        \dot{\phi} \left[
            - (X^{0,1} \iprod \bar{\partial} \phi') i \partial\bar{\partial} \phi
            + n \bar{\partial} \phi' \wedge (X^{1,0} \iprod i\partial\bar{\partial} \phi)
            \right] \wedge (i\partial\bar{\partial} \phi)^{n-1}
            = -\dot{\phi} X^{0,1} \iprod \left(\bar{\partial} \phi' \wedge (i \partial\bar{\partial} \phi)^n\right)
            = 0
        \p
    \]
    Further note that
    \begin{align}
        i n \partial \dot{\phi} \wedge \bar{\partial} \phi' \wedge (i \partial\bar{\partial} \phi)^{n-1} &= \langle \partial \dot{\phi}, \partial \phi' \rangle_{i \partial\bar{\partial} \phi} (i \partial\bar{\partial} \phi)^n
        \p
    \end{align}
    Continuing the computation from above we arrive at
    \begin{equation} \label{eq:ZhuIBPintegrand}
        \frac{d}{d s} \left[ \dot{\phi} (i \partial\bar{\partial} \phi)^n e^{-h^\phi} \right]
        = \left[ \dot{\phi}' - \langle \partial \dot{\phi}, \partial \phi' \rangle_{i \partial\bar{\partial} \phi}^2\right] \wedge e^{-h^\phi} (i \partial\bar{\partial} \phi)^n
        + n \partial \left( \dot{\phi} i \bar{\partial} \phi' e^{-h^\phi} \right) \wedge (i \partial\bar{\partial} \phi)^{n-1}
        \p
    \end{equation}
    Now \cref{eq:dsmooth} follows once we justify that the boundary term
    \[
        \int_{\partial B_R} e^{-h^\phi} \dot{\phi} i \bar{\partial} \phi' \wedge (i \partial\bar{\partial} \phi)^{n-1}
    \]
    vanishes as $R \to \infty$. Indeed, the integrand is $O(R^{2n+2}e^{- \epsilon R^2})$ (in equivariant coordinates).

    For the second claim, note that by definition of $\cE_X$ and our particular choice of 2-parameter family $\phi_{t,s}$,
    \begin{equation} \label{eq:EX2paths}
        \cE_X(\phi_t, \phi_0)= \int_0^1 \int_M \dot{\phi}_{t, 1} \MAX(\phi_{t, 1}) \d t
        \c
    \end{equation}
    and our goal is to show that we can replace $\dot{\phi}_{t, 1}$ here with $\dot{\phi}_{t, 0}$.
    Using that the RHS of \cref{eq:dsmooth} is symmetric in $s$ and $t$ we can compute
    \begin{align}
        \frac{d}{d s} \int_0^1 \int_M \dot{\phi}_{t, s} \MAX(\phi_{t, s}) \d t
            &= \int_0^1 \frac{d}{d s} \int_M \dot{\phi}_{t, s} \MAX(\phi_{t, s}) \d t\\
            &= \int_0^1 \frac{d}{d t} \int_M \phi'_{t, s} \MAX(\phi_{t, s}) \d t \\
            &= \left. \int_M \phi'_{t, s} \MAX(\phi_{t, s}) \right|_{t=0}^{t=1} = 0
        \p
    \end{align}
    The exchange of derivative and integral in the first step is justified because the RHS of \cref{eq:dsmooth} is uniformly bounded.
    We see that $s = 0$ and $s = 1$ give the same value in \cref{eq:EX2paths}.
\end{proof}

Equation \cref{eq:dsmooth} above essentially provides a formula for $\frac{d^2}{d t^2} \cE_X(\phi_t, \phi_0)$ when $\phi_t$ is a smooth family of metrics, allowing us to conclude convexity/concavity properties of $\cE_X$. We establish these in the next \namecref{thm:EXconvex}, and use \cref{thm:regularization} to generalize them to the non-smooth case. This generalizes \cite[Prop.~2.17]{BWN} to the asymptotically conical case.

\begin{prop} \label{thm:EXconvex}
    Let $\Phi_t \in \Hsing^T(a_0, b_0), 0 \leq t \leq 1$ be a bounded family of metrics. The function $t \mapsto \cE_X(\Phi_t, \Phi_0)$ is
    \begin{thmlist}
        \item convex on $(0, 1)$ if $\Phi$ is a subgeodesic.
        \item affine on $[0, 1]$ if $\Phi$ is a geodesic.
        \item concave on $[0, 1]$ if $\Phi$ is an affine line, i.e.\ $\Phi_t = (1-t) \Phi_0 + t \Phi_1$.
    \end{thmlist}
\end{prop}

\begin{notation} \label{not:Inu}
    We write $I_\nu = [1/\nu, 1 - 1/\nu]$ for $\nu \in \bbN$.
\end{notation}
\begin{proof}
    In each of the three cases we know that there exists some $L > 0$ such that
    \begin{equation} \label{eq:Lip}
        |\Phi_t - \Phi_s| \leq L (1 + r^2)
    \end{equation}
    for all $s, t \in I$. In the first two cases this is part of our definition of (sub)geodesics, while in \point3 it is immediate from the quadratic bounds on $\Phi_0, \Phi_1$.

    \point1.
    Using \cref{thm:regularization} we can approximate $\Phi$ by a decreasing sequence $\Phi^\nu$ of smooth positively curved
    subgeodesics on $I_\nu \times i \bbR \times M$ such that $\Phi^\nu_t \in \cH_2^T(a, b)$ for some $a, b > 0$. Define $E_\nu(t) \coloneqq \cE_X(\Phi^\nu_t, \Phi^\nu_{1/\nu})$, and $E(t) \coloneq \cE_X(\Phi_t, \Phi_0)$.
    Applying \cref{thm:dsmooth} to $\phi_{s, t} = s \Phi^\nu_t + (1-s) \Phi^\nu_0$ we see
    \begin{equation} \label{eq:d2E}
        \frac{d^2}{d t^2} E_\nu(t) =
            \int_M
                \left(
                    \ddot{\Phi}_t^\nu
                    - \lvert \bar{\partial} \dot{\Phi}^\nu_t \rvert^2
                \right)
                e^{-h^{\Phi^\nu_t}} \left( i \partial\bar{\partial} \Phi^\nu_t \right)^n
    \p
    \end{equation}
    Recall that a bit of linear algebra shows
    \begin{equation} \label{eq:detnplus1}
        \left(\iddball \Phi^\nu\right)^{n+1}
        =
            \frac{n+1}{2}
            \left( \ddot{\Phi}^\nu - | \bar{\partial} \dot{\Phi}^\nu |^2 \right)
            (i \partial\bar{\partial} \Phi^\nu)^n
            \wedge \frac{i}{2}\d \tau \wedge \d \bar{\tau}
        \c
    \end{equation}
    so $\Phi^\nu_t$ being a geodesic means $\ddot{\Phi}^n - |{\bar{\partial}} \dot{\Phi}^\nu|^2 \geq 0$
    and we see that $E_\nu$ is convex. Since $\cE_X(\Phi^\nu_t, \Phi^\nu_\delta) = E_\nu(t) - E_\nu(\delta)$ is convex for every $\nu$ and $\delta \geq 1/\nu$ and converges pointwise to $E(t) - E(\delta)$ by \cref{thm:Econt}, we conclude that $E$ is convex on $(0, 1)$.

    \point2.
    As before, we approximate $\Phi$ by smooth $\Phi^\nu \in \cH_2^T(a, b)$, and also write $\Phi^\infty = \Phi$ for notational convenience. Recall that by \cref{thm:dr} there exists $\epsilon = \epsilon(a, b) > 0$ such that $h^{\Phi^\nu_t} \geq \epsilon r^2 - \epsilon^{-1}$.

    By the monotone convergence theorem, $E_\nu - E_\nu(\delta) \searrow E - E(\delta)$ implies
    $E_\nu \weakto E$ as distributions on $(0, 1)$.
    So to conclude that $E$ is affine if suffices to show $E''_\nu \weakto 0$.
    By \cref{thm:dsmooth} this means we must show that for any test function $\chi$ on $I$
    \[
        \int_{I_\nu \times M}
            \chi(t) \left( \ddot{\Phi}^\nu - | \bar{\partial} \dot{\Phi}^\nu |^2 \right) e^{-h^{\Phi^\nu_t}} (i \partial\bar{\partial} \Phi^\nu)^n \d t
        \to 0
        \c
    \]
    where we only consider $\nu$ large enough so that $\supp \chi \subseteq I_\nu$.

    We pick a cutoff function $g_R = \tilde{g}(r - R)$ such that $\supp g_R \subseteq B_{R+1}$ and $g_R \equiv 1$ on $B_R$ to split the integrand above into two parts $I^\nu_R$ and $\II^\nu_R$,
    supported on $B_{R+1}$ and $M \setminus B_R$ respectively:
    \begin{align}
        I^\nu_R &= \int_{I_\nu \times M}
            \chi(t) g_R \left( \ddot{\Phi}^\nu - | \bar{\partial} \dot{\Phi}^\nu |^2 \right) e^{-h^{\Phi^\nu_t}} (i \partial\bar{\partial} \Phi^\nu)^n \d t
        \\
        \II^\nu_R &= \int_{I_\nu \times M}
            \chi(t) (1 - g_R) \left( \ddot{\Phi}^\nu - | \bar{\partial} \dot{\Phi}^\nu |^2 \right) e^{-h^{\Phi^\nu_t}} (i \partial\bar{\partial} \Phi^\nu)^n \d t
    \end{align}
    Next we rewrite $I^\nu_R$ using using \cref{eq:detnplus1}.
    For $(i \partial\bar{\partial})_{M, t} \Phi^\nu$ to make sense we need to work on a complex domain, so we introduce a dummy imaginary direction of period one to $I_\nu$. Thus instead of integrating over $I_\nu$ we integrate over $W_\nu = I_\nu \times i \bbR/\bbZ = \Omega_\nu / i \bbZ$ where $\Omega_\nu$ has complex coordinate $\tau = t + is$, but $\Phi$ is independent of $s$. Then
    \begin{align}
        I^\nu_R &=
            \frac{2}{n+1} \int_{W_\nu \times M}
                g_R \chi(t) e^{-h^{\Phi^\nu_t}}
                \big( \iddball \Phi^\nu_t \big)^{n+1} \\
         &\leq \frac{2}{n+1} \int_{W_\nu \times M}
             g_R \chi(t) e^{\epsilon^{-1} - \epsilon r^2}
                \left( \iddball \Phi^\nu_t \right)^{n+1}
        \p
    \end{align}
    Since $\left(\iddball \Phi^\nu_t\right)^{n+1} \weakto \left(\iddball \Phi_t\right)^{n+1} = 0$
    by \BT/'s convergence theorem,
    we conclude that $I^\nu_R \to 0$ as $\nu \to \infty$ for any $R$.

    We handle $\II^\nu_R$ by first using \cref{eq:ZhuIBPintegrand} and then integrating the first resulting term by parts in $t$ and the second one in $M$, obtaining
    \begin{align}
        \II^\nu_R
        &=
            \int_{I_\nu \times M} \chi(t) (1 - g_R) \left(
                \frac{d}{d t} \left[ \dot{\Phi}_t (i \partial\bar{\partial} \Phi_t)^n e^{h^\Phi_t} \right] - n \partial \left( \dot{\Phi_t} i \bar{\partial} \dot{\Phi_t} e^{-h^\Phi_t} \right) \wedge (i \partial\bar{\partial} \Phi_t)^{n-1}
            \right) \wedge \d t\\
        &=
            - \int_{I_\nu \times M}
                (1 - g_R) \dot{\chi}(t)
                \dot{\Phi}^\nu_t
                e^{-h^{\Phi^\nu_t}}
                \left( i \partial\bar{\partial} \Phi^\nu_t \right)^n \wedge \d t\\
            &\mathrel{\phantom{=}}
                - n
                \int_{I_\nu \times M}
                      \chi(t) \partial g_R
                    \wedge \dot{\Phi}^\nu_t i \bar{\partial} \dot{\Phi}^\nu_t
                    \wedge e^{-h^{\Phi^\nu_t}}
                      (i \partial\bar{\partial} \Phi^\nu_t)^{n-1} \wedge \d t \\
        &\eqcolon
            \III^\nu_R + \IV^\nu_R
        \p
    \end{align}
    To simplify the formulation of some estimates we will assume that $R = R(\mu)$ for some $\mu \geq \lambda_0$.
    Since the integrand in $\III^\nu_R$ is supported in $M \setminus \bar{B}_R$ we can apply the estimate from \cref{eq:unifMeasureBound} together with \cref{eq:Lip} to obtain
    \begin{align}
        \left\lvert \III^\nu_R \right\rvert
            &\leq
                \| \chi \|_{C^1}
                \int_{I_\nu \times (M \setminus B_R)}
                    L (1 + r^2) \MAX(\Phi^\nu_t) \wedge \d t
        \\  &\leq
                L \| \chi \|_{C^1}
                \int_{\{h^\psi \geq \mu\}}
                    (1 + R_0^2 (1 + h^\psi))
                    e^{-h^\psi}
                    (i \partial\bar{\partial} \psi)^n
        \c
    \end{align}
    where $\psi \in (\cH_2^*)^T$ is arbitrary.
    So picking $R = R(\mu)$ large enough we can make $\III^\nu_R$ arbitrarily small, uniformly in $\nu$.

    For the other term, we first integrate by parts and then split it up again:
    \begin{align}
        \IV^\nu_R
            &=
                -\frac{n}{2}
                \int_{I_\nu \times M}
                      \chi(t) \partial g_R
                    \wedge i\bar{\partial} ( \dot{\Phi}^\nu_t )^2
                    \wedge e^{-h^{\Phi^\nu_t}}
                      (i \partial\bar{\partial} \Phi^\nu_t)^{n-1}
                      \wedge \d t
        \\  &=
                +\frac{n}{2}
                \int_{I_\nu \times M}
                      \chi(t) ( \dot{\Phi}^\nu_t )^2 i\partial\bar{\partial} g_R
                    \wedge e^{-h^{\Phi^\nu_t}} (i \partial\bar{\partial} \Phi^\nu_t)^{n-1}
                    \wedge \d t
        \\  &\mathrel{\hphantom{=}}
                -\frac{n}{2}
                \int_{I_\nu \times M}
                      \chi(t) ( \dot{\Phi}^\nu_t )^2 \partial g_R
                    \wedge \bar{\partial} h^{\Phi^\nu_t}
                    \wedge e^{-h^{\Phi^\nu_t}} (i \partial\bar{\partial} \Phi^\nu_t)^{n-1}
                    \wedge \d t
        \\  &=:
            \V^\nu_R + \VI^\nu_R
        \p
    \end{align}
    Then we rewrite the integrand of $\VI^\nu_R$, using that $\bar{\partial} h^{\Phi^\nu_t} = -i X^{1,0} \iprod i \partial\bar{\partial} \Phi^\nu_t$ implies
    \begin{align}
        n \partial g_R\wedge \bar{\partial} h^{\Phi^\nu_t} e^{-h^{\Phi^\nu_t}} (i \partial\bar{\partial} \Phi^\nu_t)^{n-1}
            &= -\partial g_R \wedge i X^{1,0} \iprod (i \partial\bar{\partial} \Phi^\nu_t)^n  \\
            &= i X^{1,0} g_R \wedge (i \partial\bar{\partial} \Phi^\nu_t)^n
        \c
    \end{align}
    to obtain
    \begin{align}
        \lvert \VI^\nu_R \rvert
            &\leq
                \frac{1}{2}
                \int_{I_\nu \times M}
                    \chi(t) L^2 (1+r^2)^2
                    |X g_R|
                    \MA_X(\Phi^\nu_t) \wedge \d t
        \p
    \end{align}
    Since $X g_R = r \partial_r g_R = r \tilde{g}'(r - R) = O(R)$ the $(t,\nu)$-uniform bound
    \[
        R \,\left\lvert \int_{M \setminus B_R} (1 + r^4) \MA_X(\Phi^\nu_t) \right\rvert = O\left(e^{-\epsilon^2 R^2 / R_0^2} R^{5 + 2n}\right)
    \]
    from \cref{eq:boundMminusBall} shows that $| \VI^\nu_R | \to 0$ as $R \to \infty$, uniformly in $\nu$.

    For $\V^\nu_R$ we can estimate, similar to before,
    \[
        i \partial\bar{\partial} g_R = \tilde{g}'(r - R) \partial\bar{\partial} r + \tilde{g}''(r - R) \partial r \wedge \bar{\partial} r \leq \| \tilde{g} \|_{C^2} i \partial\bar{\partial} r^2
        \p
    \]
    Thus for any $\delta > 0$ we have on $B_{R+1}$
    \[
        | i \partial\bar{\partial} g_R \wedge (i \partial\bar{\partial} \Phi^\nu_t)^{n-1} | \leq \| \tilde{g} \|_{C^2} i \partial\bar{\partial} r^2 \wedge (i \partial\bar{\partial} \Phi^\nu_t)^{n-1} \leq \frac{1}{n} \|\tilde{g}\|_{C^2} \delta^{-1} \left(i \partial\bar{\partial} (\Phi^\nu_t + g_{R+1} \delta r^2)\right)^n
        \p
    \]
    We also have on $B_{R+1}$
    \[
        h^{\Phi^\nu_t + \delta g_{R+1} r^2} \equiv h^{\Phi^\nu_t} + \delta r^2
        \p
    \]
    By increasing $b$ if necessary we can assume that $\Phi^\nu_t + g_{R+1} \delta r^2 \in \cH_2(a, b)^T$ for all $\nu, t$ and $\delta \leq 1$.
    Writing $R \coloneq R(\mu)$, and applying \cref{eq:boundMminusBall} we then obtain the $\nu$-uniform bound
    \begin{align}
        \frac{2}{n} | \V^\nu_R | &\leq \|\chi\|_{\infty}  \| \tilde{g} \|_{C^2} \delta^{-1}
                                  \int_{I_\nu \times (B_{R+1} \setminus B_R)} L^2 (1 + r^2)^2 e^{-h^{\Phi^\nu_t}}
                                      \MA_X\left(\Phi^\nu_t + g_{R+1} \delta r^2/2 \right) \\
                   &\leq \|\chi\|_{\infty}  \| \tilde{g} \|_{C^2} \delta^{-1} e^{\delta (R+1)^2} O\left(R^{4+2n} e^{-\epsilon^2 R^2 / R_0^2}\right)
       \p
    \end{align}
    If we pick $\delta$ sufficiently small this vanishes as $R \to \infty$.

    So by first picking $R \gg 1$ such that $\III^\nu_R + \V^\nu_R + \VI^\nu_R$ is as small as necessary and then letting $\nu \to \infty$ we see that $\langle E''_\nu, \chi \rangle \to 0$. 

    As remarked before, the $t$-integrand in the definition of $\cE_X$ in \cref{eq:defEX} is uniformly bounded along a curve in $\Hsing(a, b)$. Hence $E(t)$ is continuous in $t$, so convexity on $(0, 1)$ implies convexity on $[0, 1]$.

    \point3.
    Concavity along affine lines is very easy to see. If $\Phi_t = t \Phi_1 + (1-t) \Phi_0$ then $\Phi^\nu_t = t \Phi^\nu_1 + (1-t) \Phi^\nu_0$ (this follows from the proof of \cref{thm:regularization} but we could also take this as a definition). Then when we use \cref{eq:dsmooth} to compute $\frac{d^2}{d t^2} \cE(\Phi^\nu_t, \Phi^\nu_0)$ we obtain a negative integrand on the RHS since $\frac{\partial^2 \Phi^\nu_t}{\partial t^2} = 0$.
    This means that $E_\nu'' \leq 0$ and arguing as in \point2 we conclude $E'' \leq 0$.
\end{proof}

\section{Uniqueness along geodesics}%
\label{sec:uniq_along}

We now consider two shrinkers $(M, X, \omega_i, f_i), i = 0, 1$ on $(M, X)$, with corresponding metrics $\phi_i = \log(\omega_i^n) + f_i$ on $K_M^{-1}$.
By \cref{thm:geod_exist} there exists a geodesic $\Phi$ connecting $\phi_0$ to $\phi_1$ with uniform bounds $\Phi_t \in \Hsing(a, b)$ and $|\dot{\Phi}_t| \leq L (1 + r^2)$.

In \cref{sub:uniform_l2,sub:Brunn--Mink} we will show that since both endpoints of the geodesic $\Phi$ are shrinkers, the Ding functional $D_X(\Phi_t, \Phi_0)$ along $\Phi$ must be an affine function of $t$.
Then we will use this in \cref{sub:producing_vf} to define a time-dependent vector field whose flow pulls back $\phi_1$ to $\phi_0$, establishing that shrinkers on $(M, X)$ are unique up to biholomorphism.

Since many of our tools only work for smooth metrics, we will use the approximations $\Phi^\nu$ to the geodesic $\Phi$ furnished by \cref{thm:regularization} throughout this section.
The arguments in this section mostly follow those in \cite[\S 2.4, \S 4]{BoBM}.
The main difficulty lies in obtaining uniform estimates in \cref{thm:unifL2}, since, as opposed to the compact case, the metrics $\phi \in \Hsing(a, b)^T$ are not all uniformly equivalent.

\subsection{Uniform $L^2$-estimate}%
\label{sub:uniform_l2}
\newcommand{\dbst}{\bar{\partial}^*_{\phi, \omega}}

For $\phi \in \cH_2^T$, we define operators $\partial^\phi \coloneqq e^\phi \partial e^{-\phi}: \cA^{n-1, 0}(M, K_M^{-1}) \to \cA^{n, 0}(M, K_M^{-1})$.
Formally (e.g.\ for compactly supported $v$), we have the relation
\[
    \dbst (v \wedge \omega) = i \partial^\phi v
    \c
\]
for any Kähler form $\omega$ on $M$ (locally this is equivalent to the Kähler identity $[\bar{\partial}^*, L] = i \partial$).

We pick, and fix for the rest of this section, some (complete) Kähler metric $\omega$ which is the curvature of a smooth metric $\chi \in \cH_2^T$ satisfying $\chi \geq \Phi^\nu_t$ for all $\nu, t$.
We can obtain such a metric by a similar trick as used in Step 1.2 of \cref{thm:geod_exist}:
Pick a convex increasing function $p: [0, \infty) \to [0, \infty)$ such that $p \equiv 1$ near $0$ and $p(r) = r^2/2$ for $r > 2$, and then add a sufficiently large multiple of $p(r)$ to any of the metrics $\Phi^1_t$. Since $i \partial\bar{\partial} \chi \gtrsim \frac{i}{2} \partial\bar{\partial} r^2$ for $r > 2$ this metric will be complete.

The holomorphic line bundle $\Omega^n(K_M^{-1})$ is trivialized by the \demph{tautological section} $u$ defined by
\[
    u(v_1 \wedge \cdots \wedge v_n) = v_1 \wedge \cdots \wedge v_n
\]
for any $v_i \in T^{1, 0}_p M$.

\begin{defn}
    Given a metric $\phi \in \cH_2^T$, let $L^2_{\omega, \phi} \cA^{p, q}(M; K_M^{-1})$ be the space of $K_M^{-1}$-valued $(n, 0)$-forms which are square-integrable with respect to the inner product
    \[
        (\alpha, \beta)_{\omega, \phi} = \int_M \langle \alpha, \beta \rangle_{\phi, \omega} \frac{\omega^n}{n!}
        \p
    \]
    In the special case $(p, q) = (n, 0)$,
    \[
        (\alpha, \beta)_{\omega, \phi} = \int_M e^{-\phi} \alpha \wedge \bar{\beta}
    \]
    is independent of $\omega$, so we just write $L^2_\phi \cA^{n, 0}(M; K_M^{-1})$.
    We define a projection operator $\pi_\perp^\phi$ on $L^2_\phi \cA^{n, 0}(M; K_M^{-1})$ by
    \begin{equation} \label{eq:defpiperp}
        \pi_\perp^\phi (f u) = \left[ f - \left(\int_M e^{-\phi}\right)^{-1} \int_M e^{-\phi} f i^{n^2} u \wedge \bar{u} \right] u
        \c
    \end{equation}
    where $\int_M e^{-\phi} \coloneq \int_M e^{-\phi} u \wedge \bar{u} = \|u\|^2_{\phi}$.

    We write $\pi_\perp^{\nu, t} = \pi_\perp^{\Phi^\nu_t}$ and $\pi_\perp^t = \pi_\perp^{\Phi_t}$.
    On the $L^2_{\Phi^\nu}$-sections of the bundle $\pr_M^* (\cA^{n, 0} \otimes K_M^{-1}) \to I_\nu \times M$ we then obtain operators $\pi_\perp^\nu$ and $\pi_\perp$, which for every $t$ project onto the $\Phi^\nu_t$-orthogonal complement of $u$.
\end{defn}

\begin{lem} \label{thm:Liouville}
    The only $T$-invariant holomorphic sections of $\cA^{n, 0}\left(K_M^{-1}\right)$ are constant multiples of $u$.
    Thus on $T$-invariant sections, the operator $\pi_\perp^\phi$ is the same as the projection onto the orthogonal complement of the subspace $L^2 \Omega^n(K_M^{-1})$ of square-integrable holomorphic sections.
\end{lem}
\begin{proof}
    The tautological section provides an isomorphism $\Omega^n(K_M^{-1})^T \cong \cO_M(M)^T$. Since $JX \in T$, any $f \in \cO_M(M)^T$ must be $X$-invariant. Every $X$-orbit intersects $\bar{B}_1 \subseteq M$, so $\sup_M f = \sup_{\bar{B}_1} f$, hence $f$ achieves its maximum and is therefore constant by the maximum principle.
\end{proof}

\begin{prop} \label{thm:unifL2}
    Let $\phi \in \cH_2(a, b)^T$ be such that $\phi \leq \chi$ and $i \partial\bar{\partial} \phi > \kappa \omega$ for some $\kappa > 0$ and let $f \in C^\infty(M; \bbR)^T$ satisfy $f \leq L(1 + r^2)$ for some $L > 0$. Then the equation
    \begin{equation}
        \partial^\phi v = \pi_\perp^\phi (f u)
    \end{equation}
    has a unique solution $v \in \cA^{n-1, 0}(M; K_M^{-1})^T$ such that $\bar{\partial} v$ is primitive,
    with a bound
    \[
        \left\lVert v \right\rVert_{\omega, \phi} \leq C(a, b, L)
        \p
    \]
\end{prop}
\begin{rmk} \label{rmk:unifL2Phi}
    We want to apply this to the metrics $\Phi^\nu_t$ and $f = \dot{\Phi}^\nu_t$. We obtain $v^\nu_t$, defined for $t \in I_\nu$ (recall $I_\nu = [1/\nu, 1 - 1/\nu]$), that solve
    \begin{equation} \label{eq:defvnu}
        \partial^{\Phi^\nu_t} v^\nu_t = \pi_\perp^{\nu, t}\left(\dot{\Phi}^\nu_t u\right)
        \p
    \end{equation}
    Strict positivity $i \partial\bar{\partial} \Phi^\nu_t > \kappa \omega$ for some $\kappa = \kappa(\nu, t) > 0$ is guaranteed by \cref{thm:regularization}.
    We extend $v^\nu_t$ to $I = [0, 1]$ by zero,
    obtaining $v^\nu \in L^2_{\omega, \Phi^\nu} \cA^{n-1, 0}(I \times M; K_M^{-1})^T$ such that $C \geq \| v^\nu \|_{\omega, \phi^\nu} \geq \| v^\nu \|_{\omega, \chi}$.
    By Alaoglu's theorem there exists a $L^2$-weak convergent subsequence $v^{\nu_k} \weakto v \in L^2_{\omega,\chi}\cA^{n-1, 0}(I \times M; K_M^{-1})^T$.
\end{rmk}

\begin{proof}
    To simplify notation, we shall write $V \coloneqq L^2_{\phi, \omega} \cA^{n, 0}(M; K_M^{-1})^T$, $S \coloneqq \cA^{n, 0}_c(M; K_M^{-1})^T$ (the latter space consisting of compactly supported forms) and we will omit the subscript $\omega$ from norms.

    First, recall that the mere existence and uniqueness of such a $v$ is easily obtained by standard "$L^2$-theory": Since $i \partial\bar{\partial} \phi > \kappa \omega$ by assumption and $\omega$ is complete, $i \partial\bar{\partial} \phi$ is a complete Kähler metric.
    By the standard $L^2$-estimate on complete
    manifolds \cite[Thm.~VIII.4.5]{agbook} we have for any $\alpha \in V$
    \begin{equation} \label{eq:L2}
        \| \pi_\perp^\phi \alpha \|_{\phi}^2 \leq \int_M \langle [i \partial\bar{\partial} \phi, \Lambda_\omega]^{-1} \bar{\partial} \alpha, \bar{\partial} \alpha \rangle_{\phi, \omega} \omega^n
            \leq c(\phi, \omega)^{-1} \lVert \bar{\partial} \alpha \rVert_\phi^2
        \c
    \end{equation}
    where $c(\phi, \omega) > \kappa$ is the infimum of the eigenvalues of $\omega^{-1} i \partial\bar{\partial} \phi$ over all of $M$.

    The vector $\beta \coloneqq \pi_\perp^\phi(f u)$ gives a bounded linear functional on $(\Dom \bar{\partial})^T \subseteq V$ satisfying
    \begin{nicealign} \label{eq:L2basic}
        \lvert \langle \alpha, \beta \rangle_\phi \rvert &= \lvert \langle \pi_\perp^\phi \alpha, \beta \rangle_\phi \rvert
        \leq \lVert \beta \rVert_\phi \lVert \pi_\perp^\phi \alpha \rVert_\phi \\
                         &\leq c(\phi, \omega)^{-1/2} \| \beta \|_\phi \lVert \bar{\partial} \alpha \rVert_\phi
        \c
    \end{nicealign}
    and applying Hahn--Banach we can extend the linear functional $\bar{\partial} \alpha \mapsto \langle \alpha, \beta \rangle_\phi$ from $( \Im \bar{\partial} )^T \subseteq L^2_{\phi, \omega} \cA^{n, 1}(M; K_M^{-1})$ to its $L^2_\phi$ closure, obtaining $g_\phi \in \overline{\Im \bar{\partial}}^T$ such that $\beta = \dbst g_\phi$ in the sense of unbounded operators.
    Equivalently we can define $v_\phi = \Lambda_\omega g_\phi$ so that $g_\phi = \omega \wedge v_\phi$ and $\beta = i \partial^\phi v_\phi$. Since $g_\phi$ is $\bar{\partial}$-closed, $\bar{\partial} v_\phi \wedge \omega = 0$, i.e.\ $\bar{\partial} v_\phi$ is an $\omega$-primitive form.
    If there are two $\bar{\partial}$-closed solutions $g_1, g_2$ to $\dbst g_i = f$ then
    \[
        g_1 - g_2 \in \ker \dbst \cap \ker \bar{\partial} = \scH^{n, 1}_{\phi, \omega}(M; K_M^{-1})
        \c
    \]
    where $\scH^{n,1}_{\phi, \omega}$ denotes $\Delta_{\phi, \omega}$-harmonic forms. But since $\phi$ is positively curved, the standard $L^2$ estimate for $\bar{\partial}$ shows that this space is zero. Hence we have uniqueness for $g_\phi$ and therefore also for $v_\phi$.
    Furthermore, observe that $\Delta_{\phi, \omega} g_\phi = 2 \bar{\partial} \bar{\partial}_{\phi, \omega}^* g = \bar{\partial} f \wedge u$ is smooth, so $g_\phi$ and $v_\phi$ are too, by elliptic regularity.

    From \cref{eq:L2basic} we can see that
    \begin{align} \label{eq:estv}
        \| v_\phi \|_\phi &= \| g_\phi \|_\phi = \sup_{\alpha \in S} \frac{|\langle \omega \wedge v_\phi, \bar{\partial} \alpha \rangle_\phi|}{\| \bar{\partial} \alpha \|_\phi}
        = \sup_{\alpha \in S} \frac{|\langle \beta, \alpha \rangle_\phi|}{\| \bar{\partial} \alpha \|_\phi}  \\
                            &\leq c(\phi, \omega)^{-1/2} = \kappa^{-1/2}
    \end{align}
    which does not provide the uniform bound we claimed.
    However note that for any $\alpha \in S$ and $\epsilon > 0$ we have
    \begin{align}
        \left\lVert \bar{\partial} \alpha \right\rVert^2_\phi
            \geq \left\lVert \bar{\partial} \alpha \right\rVert^2_{(1-\epsilon) \phi + \epsilon \chi}
            &\geq c((1-\epsilon) \phi + \epsilon \chi, \omega) \left\lVert \pi^{(1-\epsilon) \phi + \epsilon \chi}_\perp \alpha \right\rVert^2_{(1-\epsilon) \phi + \epsilon \chi} \\
            &\geq \epsilon \left\lVert \pi^{(1-\epsilon) \phi + \epsilon \chi}_\perp \alpha \right\rVert^2_{(1-\epsilon) \phi + \epsilon \chi}
    \end{align}
    where the first inequality is a consequence of $\chi \geq \phi$ and the second one follows from \cref{eq:L2}.
    Plugging this into the first line of \cref{eq:estv} instead of immediately using \cref{eq:L2basic} we see that
    \begin{align}
        \| v_\phi \|_\phi &\leq \epsilon^{-1/2} \sup_{\alpha \in S} \frac{|\langle \beta, \alpha \rangle_\phi|}{\| \pi^{(1-\epsilon) \phi + \epsilon \chi}_\perp \alpha \|_{(1-\epsilon) \phi + \epsilon \chi}} \\
                  &= \epsilon^{-1/2} \sup_{\alpha \in S} \frac{|\langle \beta, \alpha \rangle_\phi|}{\| \alpha \|_{(1-\epsilon) \phi + \epsilon \chi}} \\
                  &= \epsilon^{-1/2} \sup_{\alpha \in S} \frac{|\langle e^{\epsilon(\chi - \phi)} \beta, \alpha \rangle_{(1-\epsilon) \phi + \epsilon \chi}|}{\| \alpha \|_{(1-\epsilon) \phi + \epsilon \chi}} \\
                  &\leq \epsilon^{-1/2} \left\lVert e^{\epsilon (\chi - \phi)} \beta \right\rVert_{(1-\epsilon) \phi + \epsilon \chi}
        \label{eq:unifL2main}
    \end{align}
    In the first equality we use that, since $\beta \perp u$, the numerator does not change when adding a multiple of $u$ to $\alpha$, so it is sufficient to optimize over $\alpha$ such that $\alpha = \pi^{(1-\epsilon)\phi + \epsilon \chi}_\perp \alpha$.
    Plugging the definition of $\beta$ into the RHS we get
    \begin{equation} \label{eq:est_proj}
        \lVert v_\phi \rVert_\phi
            \leq \epsilon^{-1/2} \left\lVert e^{\epsilon (\chi - \phi)} \left(f u - \left( \int_M e^{-\phi} \right)^{-1} \left(\int_M e^{-\phi} f\right) u \right) \right\rVert_{(1 - \epsilon) \phi + \epsilon \chi}
        \p
    \end{equation}
    Now we pick $\epsilon \ll 1$ such that $\epsilon(\chi - \phi) \leq \frac{a}{2} r^2 + C$ for some $C > 0$. Since we have a lower bound $\phi \geq a r^2 - a^{-1}$ in equivariant charts, this can be done uniformly in $\phi$.
    The square of the RHS can be estimated, independent of $\phi$, by observing that
    \begin{align}
        \left\lVert e^{\epsilon (\chi - \phi)} f u \right\rVert^2_{(1-\epsilon) \phi + \epsilon \chi}
            &\leq L^2 \int_M e^{-\phi} e^{C -a r^2 / 2} (1 + r^2)^2 \\
        \left\lVert e^{\epsilon (\chi - \phi)} u \right\rVert^2_{(1-\epsilon) \phi + \epsilon \chi}
            &\leq \int_M e^{- \phi} e^{C - ar^2 / 2} \\
        \left\lvert \int_M e^{-\phi} f u \right\rvert
            &\leq L \int_M e^{-\phi} (1 + r^2)
        \p
    \end{align}
    and bounding the integrals using that $\phi \geq a r^2 - a^{-1}$ in equivariant charts. We also use that $\int e^{-\phi} \geq C'$ uniformly in $\phi$ since $\phi \leq b r^2 + b$ in equivariant charts.
\end{proof}

\subsection{Berndtsson's Brunn--Minkowski inequality}%
\label{sub:Brunn--Mink}

The next theorem is completely analogous to \cite[Thm.~3.1]{BoBM}, although we phrase it in slightly different terms.
\begin{restatable}{thm}{bernBM}
    \label{thm:BerndtssonFormula}
    Let $\omega$ be as at the beginning of \cref{sub:uniform_l2}.
    Let $\phi_t \in \cH_2(a, b)$ be a smooth family of metrics whose derivatives are all $O(r^2)$ in equivariant charts and satisfy $i \partial\bar{\partial} \phi_t > \kappa(t) \omega$ for some $\kappa(t) > 0$.
    Let $v_t$ be the solution to
    \begin{equation} \label{eq:defv}
        -i \bar{\partial}^*_{\omega, \phi_t} ( \omega \wedge v_t ) = \partial^{\phi_t} v_t = \pi_{\perp}( \dot{\phi}_t u )
    \end{equation}
    for which $\bar{\partial} v_t$ is primitive (which exists by \cref{thm:unifL2}). As before, we extend $\phi_t$ to $\tau = t + i s \in (0, 1) \times i \bbR$ by declaring it to be independent of $s$.
    Then
    \begin{align} \label{eq:Berndtsson}
        \left( \int_M e^{-\phi_t} \right) \frac{d^2}{d t^2} \cF \left( \phi_t \right)
            &= 4 i \frac{\partial}{\partial \tau} \iprod \frac{\partial}{\partial \bar{\tau}} \iprod \int_M \iddball \phi_t \wedge i^{n^2} \hat{u} \wedge \overline{\hat{u}} e^{-\phi_t} \\
            &\qquad+ \int_M | \bar{\partial} v_t |_{\omega, \phi_t}^2 \frac{\omega^n}{n!}
    \end{align}
    where $\hat{u} = u - \frac12 \d \tau \wedge v_t$ (we are also asserting that the integrals on the RHS converge).
\end{restatable}

\begin{proof}[Proof sketch]
    A proof in the compact case can be found in \cite[\S 4]{BerCurv}. All integrands in that proof are of the form $P(\dot{\phi}_t, D \phi_t, D^2 \phi_t) e^{-\phi_t}$ for some polynomial $P$, where $D \phi_t = (\partial_t \phi_t, \d \phi_t)$. Thus in equivariant coordinates all integrands and their derivatives can be bounded by $e^{-a r^2}$ times a polynomial in $r$, which justifies all integrations by parts and differentiations under the integral sign.
    The details of the proof can be found in \cref{apx:bernBM}.
\end{proof}

When we apply this to the family of subgeodesics $\Phi^\nu$ we have $\iddball \Phi^\nu_t \geq 0$, showing that $t \mapsto \cF(\Phi^\nu_t)$ is convex. Since $\cF$ is continuous along decreasing sequences by the monotone convergence theorem, we conclude
\thminflow
\begin{cor} \label{thm:Fconvex}
    $\cF$ is convex along geodesics.
\end{cor}
\subsection{Producing the vector field}%
\label{sub:producing_vf}

By \cref{thm:EXconvex,thm:Fconvex} the Ding functional $D_X = \vol(K_X^{-1})^{-1} \cE_X - \cF$ is concave along the geodesic $\Phi$ between the shrinkers $\phi_0$ and $\phi_1$.
Next we want to compute its derivatives at the endpoints of the interval $[0,1]$, which we do separately for $\cE_X$ and $\cF$.
By \cref{eq:Eprimitive} we would expect that
\[
    \ddat{t}{0} \cE_X(\Phi_t, \Phi_0) = \int_M \dot{\Phi}_0^+ \MAX(\Phi_0)
    \p
\]
We cannot actually show this directly in the non-smooth case, but the following inequality is good enough:
\begin{lem}["Euler--Lagrange inequality"]
    The endpoint derivatives of $\cE_X$ along $\Phi$ satisfy
    \[
        \ddat{t}{0}^+ \cE_X(\Phi_t, \Phi_0) \leq \int_M \dot{\Phi}_0^+ \MAX(\Phi_0)
        \qtext{and}
        \ddat{t}{1}^- \cE_X(\Phi_t, \Phi_0) \geq \int_M \dot{\Phi}_1^- \MAX(\Phi_1)
        \p
    \]
\end{lem}
\begin{proof}
    Consider $\psi_s = (1-s) \psi_0 + s \psi_1$ for any two (smooth!) $\psi_i \in \cH_2^T$.
    By the concavity result in \cref{thm:EXconvex} and \cref{eq:Eprimitive} we conclude
    \[
        \cE_X(\psi_1, \psi_0) \leq \ddat{s}{0}^+ \cE_X(\psi_s, \psi_0) = \int_M (\psi_1 - \psi_0) \MAX(\psi_0)
        \p
    \]
We can extend this inequality to $\psi_i \in \Hsing^T$ by approximating $\psi_i$ by smooth $\psi^\nu_i \in \cH_2^T$ from above and invoking \cref{thm:Econt}, respectively \cref{thm:psiMAphiConverge}, to see that the left, respectively right, hand side converge.
    Applying this inequality to $\psi_0 = \Phi_0$ and $\psi_1 = \Phi_t$ we obtain
    \begin{align}
        \ddat{t}{0}^+ \cE_X(\Phi_t, \Phi_0) \coloneqq \lim_{t \searrow 0} \frac{\cE(\Phi_t, \Phi_0)}{t}
                        &\leq \lim_{t \searrow 0} \int_M \frac{\Phi_t - \Phi_0}{t} \MAX(\Phi_0) \\
                        &= \int_M \dot{\Phi}_0^+ \MAX(\Phi_0)
        \c
    \end{align}
    where we use that $\Phi_t$ is convex as a function of $t$ to apply monotone convergence in the last step.
    The bound at $t = 1$ is analogous.
\end{proof}

Being convex, $\cF$ also has one-sided derivatives at the endpoints:%
\[
    \ddat{t}{0}^+ \cF(\Phi_t) = \frac{\int_M \dot{\Phi}_0^+ e^{-\Phi_0}}{\int_M e^{-\Phi_0}}
    \qtext{and similarly}
    \ddat{t}{1}^- \cF(\Phi_t) = \frac{\int_M \dot{\Phi}_1^- e^{-\Phi_1}}{\int_M e^{-\Phi_1}}
    \p
\]
    They can be computed by the dominated convergence theorem, using that
    \(
    h^{-1} \left\lvert e^{-\Phi_{t+h}} - e^{-\Phi_t} \right\rvert
        \leq e^{- \Phi_t + L h (1 + r^2)} L (1 + r^2)
    \)
    is integrable when $L h < a$ (recall that $a r^2$ is the asymptotic lower bound on $\Phi_t$ in equivariant charts).

The shrinker equation for $\phi_{0,1}$ can be rewritten as
\begin{equation}
    e^{-f_i} (i \partial\bar{\partial} \phi_i)^n = e^{-\phi_i}
\end{equation}
and $f_i = h^{\phi_i} + \mathrm{const.}$ since the $i \partial\bar{\partial} \phi_i$-gradient of both functions is the soliton vector field $X$.
Hence the soliton equation is equivalent to the condition that $\MAX(\phi_i)$ and $e^{-\phi_i}$ are equal up to multiplication by a constant.
Therefore
\begin{align}
    \frac{d}{d t}\bigg|_0^+ D_X(\Phi_t, \Phi_0) &\leq \int_M \dot{\Phi}^+ \left[ \vol(K_M^{-1})^{-1} \MAX(\phi_0) - \left(\int_M e^{-\phi_0}\right)^{-1} e^{-\phi_0} \right] = 0
    \\
    \frac{d}{d t}\bigg|_1^- D_X(\Phi_t, \Phi_0) &\geq \cdots = 0
    \c
\end{align}
since the measures in the square brackets are proportional and both are normalized to unit mass.
At the same time, $D_X$ is concave, so it must be constant along $\Phi$. Since we saw that $\cE_X$ is affine along $\Phi$ in \cref{thm:EXconvex}\point2 we conclude that $\cF$ is affine along $\Phi$.

At this point we will briefly explain the rest of the argument, ignoring issues of regularity and boundedness:
Since $\cF(\Phi_t)$ is affine, applying \cref{thm:BerndtssonFormula} to $\Phi_t$ we see that both terms on the RHS of \cref{eq:Berndtsson} must vanish -- implying $\bar{\partial} v_t = 0$.
Then taking $\bar{\partial}$ of \cref{eq:defv} we obtain
\[
    \partial\bar{\partial} \Phi_t \wedge v_t = \dot{\Phi}_t u
    \qquad
    \text{or equivalently}
    \qquad
    \cL_{V_t} \bar{\partial} \Phi_t = \frac{d}{d t} \bar{\partial} \dot{\Phi}
    \c
\]
in terms of the vector field $V_t$ satisfying $V_t \iprod u = v_t$.
If $F_t$ is the time-dependent flow of the vector field $\Re V_t$ the above formula shows that $F_t^* i \partial\bar{\partial} \Phi_t$ is time-independent, and $\psi = F_1$ is the biholomorphism whose existence we claimed in \cref{thm:uniq}.

\subsection{Approximation arguments}%
\label{sub:approximation_arguments}

If we define
\[
    \cF_\nu(t) = -\log \int_M e^{-\Phi^\nu_t}
    \c
\]
then by monotone convergence, $\cF_\nu \weakto \cF$ and hence $\cF_\nu'' \weakto \cF'' = 0$.
As in \cref{rmk:unifL2Phi} we have $v^\nu \in \cA^{n-1, 0}(I_\nu \times M; K_M^{-1})^T$ satisfying Berndtsson's formula \cref{eq:Berndtsson} for $\phi = \Phi^\nu$. Since $\cF_\nu'' \weakto 0$ and both terms on the RHS of \cref{eq:Berndtsson} are non-negative, we can pair \cref{eq:Berndtsson} with a bump function on $(\epsilon, 1-\epsilon)$ to conclude that
\[
    \| \bar{\partial} v^\nu \|_{L^2_\chi([\epsilon, 1-\epsilon] \times M)} \leq 
    \| \bar{\partial} v^\nu \|_{L^2_{\phi^\nu}([\epsilon, 1-\epsilon] \times M)} \to 0
\]
for any $\epsilon > 0$.
Consequently $\bar{\partial} v^\nu \weakto 0$ as distributions. We already know that $v^\nu$ has a $L^2_\chi(I \times M)$-weak limit $v$, so we conclude $\bar{\partial} v = 0$ as distributions.

The equation $\bar{\partial} v = 0$ is "elliptic in $M$-directions" so $v$ is smooth in those directions. More precisely:
\begin{lem}
    If $\bar{\partial}_M v = 0$ for some $L^2$ section $v$ of a hermitian line bundle $(\pi^* L, h)$ on $I \times M$ then $v \in L^2(I; L^2 H^0(M, L))$.
\end{lem}

\begin{proof}[Proof sketch]
    We have $v \in L^2(I; L^2_h(M; L))$ meaning $v_t$ is an $L^2$ section of $(L, h)$ for every
    $t$
    (we do not need to say "almost every $t$" because we could just re-define $v_t = 0$ for problematic values of $t$).
    The distributional equation $\bar{\partial} v = 0$ implies that for every ($L^2$-valued) test form $\alpha$ and every $\psi \in \cD(I)$:
    \[
        0 = \langle \bar{\partial} v, \psi \boxtimes \alpha \rangle = \int_0^1 \psi(t) \int_M v_t \wedge \bar{\partial} \alpha \d t
        \p
    \]
    Since this is true for any $\psi$ we see that $\bar{\partial} v_t = 0$ weakly for (almost) every $t$.
    The claim now follows from elliptic regularity.
\end{proof}

Our next goal is to take the limit of \cref{eq:defvnu} to obtain $\partial^{\Phi_t} v_t = \pi^t_\perp (\dot{\Phi}_t u)$.
Integrating by parts, we know from \cref{eq:defvnu} that for any $w \in \cD^{n, 0}(I \times M)$
\begin{equation} \label{eq:weakdelPhi}
    \int_{I \times M} \d t \wedge v^\nu \wedge \overline{\bar{\partial} w} e^{-\Phi^\nu}
        = \int_{I \times M} \d t \wedge \pi_\perp^{\nu, t} (\dot{\Phi}^\nu_t u) \wedge \bar{w} e^{-\Phi^\nu}
    \p
\end{equation}
To take the limit of the LHS we first use that $e^{\chi-\Phi^\nu} \bar{\partial} w \to e^{\chi-\Phi} \bar{\partial} w $ in $L^2_\chi$ by monotone convergence.
We pair this with the $L^2_\chi$-weak convergence of $v^\nu$ established in \cref{rmk:unifL2Phi} to obtain
\[
    \int_{I \times M} \d t \wedge v^\nu \wedge \overline{\bar{\partial} w} e^{-\Phi^\nu} \to \int_{I \times M} \d t \wedge v \wedge \overline{\bar{\partial} w} e^{-\Phi}
    \p
\]
Meanwhile the RHS of \cref{eq:weakdelPhi} can be rewritten as
\[
    \int_{I \times M} \d t \wedge \pi_\perp^{\nu, t} (\dot{\Phi}^\nu_t u) \wedge \bar{w} e^{-\Phi^\nu}
        = \int_{I \times M} \d t \wedge \dot{\Phi}^\nu u \wedge \overline{\pi_\perp^{\nu, t}(w)} e^{-\Phi^\nu}
    \p
\]
Since
\[
    \big\lVert \dot{\Phi}^\nu \big\rVert^2_\Phi \leq L^2 \int_{I \times M} \d t \wedge (1 + r^2)^2 e^{-\Phi_t}
\]
is uniformly bounded, and $\dot{\Phi}^\nu \weakto \dot{\Phi}$ as distributions, we can conclude that $\dot{\Phi}^\nu \to \dot{\Phi}$ weakly in $L^2_\Phi$.
\begin{claim*} \label{thm:piperpConv}
    For any test function $f \in C_c((0, 1) \times M)$ we have $\pi_\perp^\nu (f u) \to \pi_\perp (f u)$ in $L^2_\Phi$.
\end{claim*}
\vspace{-0.7\topsep}
\begin{proof}
    W.l.o.g.\ we may assume $f \geq 0$. Then by monotone convergence
    \[
        c_\nu(t) \coloneqq \left(\int_M e^{-\Phi^\nu_t}\right)^{-1} \int_M e^{-\Phi^\nu_t} f_t \to \left(\int_M e^{-\Phi_t}\right)^{-1} \int_M e^{-\Phi_t} f_t
    \]
    for every $t$, and all $c_\nu$ are uniformly bounded by a function of $a, b$. Since $\pi^\nu_\perp (f u) - \pi_\perp (f u) = (c_\nu(t) - c_\infty(t)) u$ and $\| u \|_\Phi < \infty$, the claim follows from the dominated convergence theorem.
\end{proof}

Furthermore,
\begin{align}
    \left\lVert e^{\Phi - \Phi_\nu} \pi_\perp^\nu(w) - \pi_\perp(w) \right\rVert_\Phi
        &\leq
            \left\lVert e^{\Phi - \Phi_\nu} \big(\pi_\perp^\nu(w) - \pi_\perp(w)\big) \right\rVert_\Phi +
            \left\lVert (1 - e^{\Phi - \Phi_\nu} ) \pi_\perp(w) \right\rVert_\Phi
    \\  &\leq
            \left\lVert \pi_\perp^\nu(w) - \pi_\perp(w) \right\rVert_\Phi +
            \left\lVert \left( 1 - e^{\Phi - \Phi_\nu} \right) \pi_\perp(w) \right\rVert_\Phi
    \p
\end{align}
So using the claim for the first summand, and monotone convergence for the second, we also obtain $e^{\Phi - \Phi_\nu} \pi_\perp^\nu(w) \to \pi_\perp(w)$ in $L^2_\Phi$.
Therefore the RHS of \cref{eq:weakdelPhi} also converges to the obvious limit, meaning that the limit of \cref{eq:weakdelPhi} is
\[
    \int_I \int_M v_t \wedge \overline{\bar{\partial} w_t} e^{-\Phi_t} \d t
        = \int_I \int_{M} \pi^t_\perp(\dot{\Phi}_t) u \wedge \overline{w_t} e^{-\Phi_t} \d t
    \c
\]
and this holds for any test form $w$ on $I \times M$.
We can integrate the LHS by parts on $M$ and use that the resulting integrands must be equal for almost every $t$ to see that
\[
    \partial\left(e^{-\Phi_t} v_t\right) = e^{-\Phi_t} \pi_\perp^t(\dot{\Phi}_t) u
\]
as distributions, for a.e.\ $t$. Recalling that PSh functions have $W^{1,1}_\mathrm{loc}$-regularity \cite[Thm.~1.48]{GZ} we see that $\partial (e^{-\Phi_t} v_t) \in L^1_\mathrm{loc}$, so its product with the $L^\infty_\mathrm{loc}$ function $e^{\Phi_t}$ is a well-defined distribution. Hence
\[
    \pi_\perp^t(\dot{\Phi}_t) u = e^{\Phi_t} \partial \left(e^{-\Phi_t} v_t\right) = e^{\Phi_t} \left( e^{-\Phi_t} \partial v_t - e^{-\Phi_t} \partial \Phi_t \wedge v_t \right) = \partial v_t - \partial \Phi_t \wedge v_t
    \c
\]
where every term is a $L^1_\mathrm{loc}$ function.
Since $v_t$ is holomorphic, taking $\bar{\partial}$ gives us
\[
    \partial\bar{\partial} \Phi_t \wedge v_t = \bar{\partial} \dot{\Phi}_t \wedge u
\]
for a.e.\ $t$. Writing $v = V \iprod u$ for some $V \in L^2\left(I; L^2_{\Phi_t} H^0(M, \cT_M)^T\right)$ we can equivalently say
\begin{equation}
    V_t \iprod \partial\bar{\partial} \Phi_t = \bar{\partial} \dot{\Phi}_t
    \qtext{or}
    ( \cL_{V_t} - \partial_t )\bar{\partial} \Phi_t = 0
\end{equation}
or, in terms of the real-holomorphic vector field $\Re V_t$ (cf.\ \cref{eq:Hamilton}),
\begin{align} \label{eq:main}
    \Re V_t \iprod \ddc \Phi_t = \d^c \dot{\Phi}_t
    \qtext{or}
    ( \cL_{\Re V_t} - \partial_t ) \d^c \Phi_t = 0
    \p
\end{align}
\subsection{Flow of the vector field}%
\label{sub:flow_of_the_vector_field}
We have $\Re V \in L^2(I; \aut(M, X))$, where $\aut(M, X)$ is the Lie algebra of real-holomorphic vector fields commuting with $X$. Even though $\Re V$ is not a continuous time-dependent vector field, this regularity is sufficient to define its flow. The relevant results about ODEs can be found in \cref{apx:odes}.

Since $\aut(M, X)$ is finite dimensional \cite[Prop.\ 5.9]{CDS}, the norm of any $Z \in \aut(M, X)$ controls its local Lipschitz constants, so
$\Re V|_{I \times U} \in L^2(I; C^{0,1}(U)) \subseteq L^1(I; C^{0,1}(U))$ for any open bounded $U \subseteq M$. Applying \cref{thm:PicLind} we obtain a continuous map $F: [0, T] \times U \to U$  such that the time derivative $\dot{F}$ exists almost everywhere and solves the ODE
\[
    \dot{F}_t = -\Re V_t \circ F_t
    \p
\]
\begin{claim*}
    $F_t$ is holomorphic for \emph{every} $t$.
\end{claim*}
\vspace{-0.7\topsep}
\begin{proof}
    $\Re V_t$ is real-holomorphic, meaning that $[J, \d (\Re V_t)] = 0$ in holomorphic coordinates. Thus
    \[
        \frac{d}{d t} [J, \d_x F_t] = - \d (\Re V_t) \circ [J, \d_x F_t]
        \c
    \]
    where we have computed $\frac{d}{d t} \d_x F_t$ using \cref{thm:odeDiff}.
    Now \cref{thm:PicLind} guarantees that this ODE has a unique solution for the initial value $[J, \d_x F_0] = 0$, showing $[J, \d_x F_t] = 0$ for all $x, t$.
\end{proof}

Since $\Re V_t \in \aut(M, X)$ we have $(\gamma_s)_* \Re V_t = \Re V_t$ for all $s, t$ and hence also
\[
    F_t \circ \gamma_{s} = \gamma_{s} \circ F_t
    \p
\]
Since we have a lower bound on the maximum existence interval on any compact subset of $M$, this implies that $F_t: M \to M$ is defined for all $t$ by the argument in \cite[Lem.\ 2.34]{CDS}.

Pick some chart $U$ of $M$ so that we can treat $\Phi$ as a function in this chart and consider an arbitrary test form (of appropriate degree) $\alpha$ such that ${F_t}_* \alpha$ is supported in that chart for $t$ in some range $(t_1, t_2)$.
Using that $\Phi_t$ is Lipschitz in $t$ and ${F_t}_* i \partial\bar{\partial} \alpha$ is compactly supported, we can differentiate under the integral sign obtaining
\[
    \frac{d}{d t} \int_U \Phi_t {F_t}_* \ddc \alpha = \int_U \left[ \dot{\Phi}_t {F_t}_* \ddc \alpha - \Phi_t \cL_{-\Re V_t} {F_t}_* \ddc \alpha \right]
        = \int_U \left[ \d^c \dot{\Phi}_t - \cL_{\Re V_t} \d^c \Phi_t \right] \wedge {F_t}_* \d \alpha = 0
\]
for $t \in (t_1, t_2)$%
.
By \cref{eq:main} the RHS vanishes, and since $\alpha$ was arbitrary we conclude that the current $F_t^* (i \partial\bar{\partial} \Phi_t|_U) = \frac12 F_t^* (\ddc \Phi_t|_U)$ is independent of $t \in (t_0, t_1)$ in the chart. Since the integral on the LHS is continuous in $t$ we even conclude equality for $t \in [t_0, t_1]$, and since $U, t_0, t_1$ were arbitrary we see that $F_t^* i \partial\bar{\partial} \Phi_t$ is independent of $t \in [0, 1]$ globally,
implying $i \partial\bar{\partial} \Phi_0 = F_1^* i \partial\bar{\partial} \Phi_1$.
Since $i \partial\bar{\partial} \Phi_0$ and $i \partial\bar{\partial} \Phi_1$ were picked to be two arbitrary shrinkers with the same soliton vector field at the beginning of \cref{sec:uniq_along},
we have proven \cref{thm:uniq} in the case $X_1 = X_2$ (with $\psi = F_1$ as our biholomorphism).

Our uniqueness theorem also allows us to re-prove the Matsushima theorem \cite[Thm.~5.1]{CDS} for asymptotically conical shrinkers. Recall that all vector fields commuting with $X$ are complete, so the Lie algebra of $\Aut(M, X)$ is $\aut(M, X)$.

\begin{cor}[Matsushima theorem] \label{thm:Matsushima}
    Let $(M, \omega, X)$ be a shrinker with Ricci curvature tending to zero at infinity.
    Then the Lie algebra $\aut(M, X)$ is the complexification of the Lie algebra $\fg_X = \aut(M, \omega, X)$ of holomorphic Killing fields commuting with $X$, and thus reductive. Equivalently, $\Aut(M, X)$ is reductive, with maximal compact subgroup $G_X \coloneq \Aut(M, \omega, X)$.
\end{cor}
\begin{proof}
    The group $G_X$ is compact since it acts faithfully by isometries on the Link of the asymptotic cone of $M$.
    Therefore the Lie group $G_X^\bbC \subseteq \Aut(M, X)$ with Lie algebra $\fg_X^\bbC \coloneq \fg_X \otimes \bbC \subseteq \aut(M, X)$ is reductive with maximal compact subgroup $G_X$.

    Let $\psi \in \Aut(M, X)$. Then $(M, \psi_* \omega, X)$ is another shrinker with Ricci curvature tending to zero at infinity. Applying the arguments above to $\omega_0 = \omega$ and $\omega_1 = \psi_* \omega$ we obtain the continuous family of biholomorphisms $F: I \times M \to M$ such that $F_t^* i\partial\bar{\partial} \Phi_t = \omega_0$ is constant. Then pulling back \cref{eq:main} by $F_t$ we obtain
    \begin{equation}
        (F_t^* \Re V_t) \iprod \omega = \frac12 \d^c F_t^* \dot{\Phi}_t
        \p
    \end{equation}
    The exact form of the RHS is not important, we only want to see that $J F_t^* \Re V_t$ is Hamiltonian (hence Killing), or in other words
    $u_t \coloneq - F_t^* \Re V_t \in J \fg_X$.
    Now consider the time-dependent left-invariant vector field $U_t$ on $G_X^\bbC$ corresponding to $u_t$. Its integral curve $a$ starting at the identity (which exists by \cref{thm:PicLind}) satisfies the ODE
    \[
        \dot{a}(t) = U_t (a(t)) = \d L_{a(t)} u_t
        \p
    \]
    This is essentially the same ODE as
    \[
        \dot{F}_t = - \d F_t (F_t^* \Re V_t) = \d F_t (u_t)
        \c
    \]
    therefore $F_t(p) = a(t)(p)$ by the uniqueness part of \cref{thm:PicLind}. In particular, $F_1 \in G_X^\bbC$ and since $F_1^* \psi^* \omega_0 = \omega_0$ we have $F_1 \circ \psi \in \Isom(M, \omega, X)$, so $\psi \in G_X^\bbC$. Since $\psi$ was arbitrary we see that $G_X^\bbC = \Aut(M, X)$ and $\fg_X^\bbC = \aut(M, X)$.
\end{proof}

\section{Vector fields in different maximal tori}%
\label{sec:vf_uniq}

\begin{notation}
    A lowercase Fraktur letter always denotes the Lie algebra of the group denoted by the same uppercase letter.
\end{notation}

Finally we show that two shrinkers $(\omega_1, X_1)$ and $(\omega_2, X_2)$ with Ricci curvature decay%
\footnote{recall that by \cite[Cor.~2.3]{MWconical} and \cite[Thm.~A]{CDS} this implies the metrics are asymptotically conical}
on the same complex manifold $M$ are equal up to automorphism.
The proof idea is simple: We pick maximal tori $T_1, T_2$ such that $X_i \in \ft_i$ and find $\psi \in \Aut(M)$ such that $\psi T_1 \psi^{-1} = T_2$. Then $\psi_* X_1 \in \ft_2$ and since $(\psi_* \omega_1, \psi_* X_1)$ is a shrinker with bounded Ricci tensor we can conclude that $X_2 = \psi_* X_1$ by the volume minimization principle in \cite{CDS}. The results of the previous sections the allow us to conclude $\psi_* \omega_1 = \omega_2$.

The difficulty lies in showing the existence of $\psi$. This is easy in the compact case \cite[Thm.~3.2]{TianZhu02} because $\Aut(M)$ is a Lie group, so all maximal tori are conjugate.
Unfortunately, no results are known to us which would guarantee that maximal tori in the automorphism group (holomorphic or algebraic) of a quasiprojective manifold are conjugate.
Therefore our strategy is to restrict the action of the tori to the \emph{compact} exceptional set $E$ and conjugate the tori in the finite-dimensional automorphism group of $E$.
This shows that there exists a biholomorphism $\hat{g}$ defined around $E$ such that $\hat{g}^* X_2 = X_1$ up to terms vanishing on $E$ to arbitrarily high order. Then \cref{thm:AutGermsv2}, the main technical ingredient for this section, allows us to produce another biholomorphism $\eta$ of a neighborhood of $E$ such that $X_1 = \eta_* \hat{g}^* X_2$, so we can take $\psi = \hat{g} \eta^{-1}$.

This method can be applied more generally to polarized Fano fibrations (see the definitions below), which include shrinkers without any curvature decay, as well as Kähler cones.

\newcommand{\IE}{\scI}

\subsection{Polarized Fano fibrations}%
\label{sub:polarized_fano_fibrations}

We now generalize our setting to polarized Fano fibrations.
This is not directly relevant to our main result about the uniqueness of asymptotically conical shrinkers. However, the method we will use to prove the uniqueness of the soliton vector field naturally fits into the context of polarized Fano fibrations, and furthermore we can use it to prove that polarized Fano fibrations are isomorphic (as algebraic varieties) if{}f they are biholomorphic.
A reader who is only interested in the case of Kähler--Ricci shrinkers can just skip to \cref{sub:conjugating} as we will keep our notation consistent.

\begin{defn}[{\cite[Def~2.1]{SunZhang}}]
    \label{def:ff}
    A \demph{Fano fibration} consists of a surjective projective morphism $\pi: M \to Y$, where
    \begin{enumerate}
        \item $M$ is a normal variety with klt singularities and $-K_M$ is a $\pi$-ample $\bbQ$-Cartier divisor
        \item $Y$ is a normal affine variety
        \item $\pi_* \Oalg_M = \Oalg_Y$.
    \end{enumerate}
\end{defn}

Recall that the algebraic and analytic definitions of normality and properness are equivalent \cite[Prop.\ XII.2.1, 3.2]{SGA1}.

As the following argument shows, $Y$ and $\pi$ are determined by the variety $M$, so the data of a Fano fibration does not need to include $Y$ or $\pi$. Thus we can simply call the variety $M$ a Fano fibration.
\begin{lem}
    If $\pi: M \to Y$ is Fano fibration then $Y = \Spec \Oalg_M(M)$ and $\pi$ is the canonical map $M \to \Spec \Oalg_M(M)$.
\end{lem}
\begin{proof}
    First recall that there is a canonical map $\pi: M \to \Spec \Oalg_M(M)$ which is initial among all maps form $M$ to affine schemes. Let $\pi': M \to \Spec A$ be another Fano fibration.
    Then $\pi' = f \circ \pi$ for some $f: \Spec \Oalg_M(M) \to \Spec A$, but $\cO_{\Spec A} = \pi'_* \Oalg_M = f_* \pi_* \Oalg_M$, so $f$ is an isomorphism of affine schemes.
\end{proof}

In fact, the \emph{biholomorphism} class of $Y$ is also determined by the biholomorphism class of $M$:
\begin{lem} \label{thm:coneAnalytic}
    Let $\pi_1: M_1 \to Y_1$ and $\pi_2: M_2 \to Y_2$ be two Fano fibrations with the same underlying analytic space $M$.
    Then there is an isomorphism of analytic spaces $Y_1 \cong Y_2$ commuting with the maps $\pi_i$.
\end{lem}
\begin{proof}
    Since $\pi_i$ are closed continuous surjections (in the euclidean topology), they are topological quotient maps.
    The variety $Y$ cannot have any compact analytic subset besides points, so by the proper mapping theorem connected compact analytic subsets of $M$ must be contained in a fiber of $\pi$. Thus the fibers of $\pi$ are finite unions of maximal connected compact analytic subsets of $M$.
    Since $\pi_i$ is proper, $({\pi_i}_* \Oalg_{M_i})^\mathrm{an} \cong {\pi_i}_* \cO_M$ by \cite[Thm.\ XII.4.2]{SGA1}. Thus ${\pi_i}_* \cO_M = \cO_{Y_i}$, which also shows that the fibers of $\pi_i$ are connected, so they are precisely the maximal connected compact analytic subsets of $M$.
    Therefore both quotient maps have the same fibers, and we conclude that $Y_1$ and $Y_2$ are homeomorphic.
    Finally the analytic structure sheaves of $Y_1$ and $Y_2$ are both given by the pushforward of $\cO_M$, so they agree.
\end{proof}
This lemma also shows that whenever we have a biholomorphism $\phi: M \to M$ of the analytic space underlying a Fano fibration, it will induce a biholomorphism $\bar{\phi}$ of $Y$ determined by $\bar{\phi} \circ \pi = \pi \circ \phi$.

The definition that is relevant for this work is:
\begin{defn}[{\cite[Def~2.5]{SunZhang}, \cite[Def~2.4]{ColSzeStab}}]
    \label{def:pff}
    Let $\pi: M \to Y$ be a Fano fibration equipped with the faithful algebraic action of a torus $T$ on $M$ and $Y$ for which $\pi$ is equivariant, and such that $T$ has exactly one fixed point $o \in Y$. The algebra $R = \Oalg_M(M) = \Oalg_Y(Y)$ has a weight
    decomposition
    \[
        R = \bigoplus_{\alpha \in \ft^*} R_\alpha
        \c
    \]
    meaning that $\cL_\xi f = i \alpha(\xi) f$ for $f \in R_\alpha$ and $\xi \in \ft$.
    The \demph{Reeb cone} of $M$ consists of $\xi \in \ft$ with the property that $\langle \alpha, \xi \rangle > 0$ for all weights $\alpha \neq 0$ such that $R_\alpha \neq \{0\}$.

    A \demph{polarized} Fano fibration is a triple $(\pi: M \to Y, T, \xi)$ as above with $\xi$ in the Reeb cone of $M$. Since the action of $\ft$ on $M$ is faithful we can consider $\xi$ to be a (real-)holomorphic vector field on $M$. In the case $M = Y$ the triple $(Y, T, \xi)$ is called a \demph{polarized affine cone}.
\end{defn}

Sun and Zhang have proved that if $(M, \omega, X)$ is a shrinker then $M$ is quasiprojective and $(\pi: M \to Y, T, JX)$ is a polarized Fano fibration, where $T$ is the torus generated by $JX$ and $Y = \Spec \Oalg_M(M)$ \cite[Thm.~3.1]{SunZhang}.

\begin{notation} \label{not:Exp}
    We write $\Exp(V): U \to M$ for the time 1 flow of a (time-independent) vector field $V \in \Re H^0(M, \cT_M)$. Here $U \subseteq M$ is the maximal open subset on which the flow can be defined.
\end{notation}
We briefly recall how to define the flow of a holomorphic vector field on an analytic space $M$.
By a holomorphic vector field we mean (the real part of) a section $X$ of the sheaf $\cT_M = (\Omega^1_M)^\vee$.
We can pick embeddings $j: W \to V$ of an open $W \subseteq M$ into an open $V \subseteq \bbC^m$ and extend $X|_W$ to a holomorphic vector field $\tilde{X}$ on $V$. The flow of $\tilde{X}$ will preserve $j(W)$, and any two extensions will have the same flow on $j(W)$. Thus the local flows glue to a flow on the whole variety $M$. As in the smooth case, the flow has a nonzero existence time on any compact set.
We can also define the flow of a time-dependent vector field in the same way.

\begin{rmk} \label{rmk:complex_torus}
    The $T$-action on a polarized Fano fibration $M$ complexifies to a holomorphic action $T_\bbC \times M \to M$ of the complexified torus $T_\bbC \cong (\bbC^*)^m$.
    At the Lie algebra level we just let the vector $i V \in i \ft$ act by $J V$. For $V, W \in \ft$ we have $[V, J W] = J [V, W] = 0$ since $V$ is real-holomorphic, so this defines a infinitesimal action of Lie algebra $\ft \otimes  \bbC$.
    This action is faithful (because the complexification of the faithful action of $\ft$ on the vector space $R$ has to be faithful), so we will view $\ft \otimes \bbC \cong \ft \oplus J \ft$ as a subalgebra of $\Re H^0(M, \cT_M)$.
    Then from the weight decomposition of $R$ it is easy to see that all the vector fields in $J \ft$ are complete on $Y = \Spec R$ and thus complete on $M$ since $\pi$ is proper.
\end{rmk}

From the requirement in \cref{def:pff} that there exists only one $T$-fixed point on $Y$ we obtain
\begin{lem}
    For a polarized Fano fibration $R_0 = \bbC$, i.e.\ no non-constant $f \in R = \Oalg_Y(Y)$ can have zero $T$-weight.
\end{lem}
\begin{proof}
    As outlined in \cite[Lem.~2.5]{ColSzeStab} we can pick weight vectors $f_1, \dots, f_N \in R \setminus \bbC$ that generate $R$, thus obtaining a closed equivariant embedding $Y \to \bbC^N$ of varieties. If the weight space $R_0$ is nontrivial then w.l.o.g.\ the first $k$ coordinate functions on $\bbC^N$, will have $T$-weight zero. Thus for any $p \in Y \subseteq \bbC^N$, the first $k$ coordinate functions will be constant along $\Exp(t J \xi) (p)$. But $\lim_{t \to \infty} \Exp(t J \xi)(p) \in \bbC^k$ is a $T$-fixed point. Since there can be only one fixed point we conclude that $f_1, \dots, f_k$ must be constant on $Y$, a contradiction.
\end{proof}

The next lemma shows that the condition of being a Reeb vector field only depends on the vector field, and not on the choice of torus $T$ containing it in its Lie algebra.
\begin{lem} \label{thm:Reebflow}
    The Reeb cone $\Lambda$ of $(M, T)$ can be characterized as
    \begin{equation} \label{eq:Reebv3}
        \Lambda = \ft \cap \left\{\xi \in \Re H^0(M, \cT_M) \:\middle|\: \forall p \in M: \lim_{t \to \infty} \Exp(t J\xi)(p) \in \pi^{-1}(o) \right\}
        \p
    \end{equation}
    (The existence of the limits is part of the condition).
\end{lem}
\begin{proof}
    Let $R = \Oalg_M(M) \cong \Oalg_Y(Y)$. If $f \in R \setminus \bbC$ has $T$-weight $\alpha: \ft \to \bbR$ then for any $\xi \in \Lambda$
    \[
        \Exp(t J \xi)^* f = e^{-\alpha(\xi) t} f \to 0
    \]
    as $t \to \infty$ because $\alpha(\xi) > 0$. Since the cone $Y$ can be embedded into $\bbC^N$ by weight vectors $f_1, \dots, f_N \in R \setminus \bbC$ \cite[Lem.\ 2.5]{ColSzeStab}, we see that
    \[
        \lim_{t \to \infty} \Exp(t J\xi)(p) \in (f_1, \dots, f_N)^{-1}(0) = \pi^{-1}(o)
        \p
    \]
    Conversely, suppose $\xi$ is in the RHS of \cref{eq:Reebv3} and there exists $f \in R$ of $T$-weight $\alpha$ such that $\alpha(\xi) \leq 0$. Then for any $t \in \bbR$
    \[
        |f(p)| = \left| e^{\alpha(\xi) t} f(\Exp(t J\xi)(p)) \right| = \lim_{t \to \infty} \left| e^{\alpha(\xi) t} f(\Exp(t J\xi)(p)) \right| \leq \sup |f(\pi^{-1}(o))|
        \p
    \]
    So by Liouville's theorem (applied to the affine variety $Y$), $f$ must be a constant.
\end{proof}

\begin{rmk}
    If $Y$ is smooth away from the vertex $o \in Y$ (the unique $T$-fixed point) then the Kähler metric induced by the embedding mentioned above makes $Y \setminus \{o\}$ into a Kähler cone.
\end{rmk}

As it turns out, $\xi$ (viewed as a holomorphic vector field on $M$) is sufficient to determine the algebraic structure of a polarized Fano fibration.

\begin{prop} \label{thm:vfDeterminesAlg}
    Let $M_1$ and $M_2$ be polarized Fano fibrations with the same underlying analytic space $M$ and the same vector field $X \in H^0(M, \cT_M)$. Then $M_1 \cong M_2$ as algebraic varieties.
\end{prop}
\begin{proof}
    We can replace the tori $T_i$ by the sub-torus generated by $\xi$ and thus assume $T_1 = T_2 \eqcolon T$.
    By \cref{thm:coneAnalytic} the cones $Y_i \coloneq \Spec \Oalg(M_i)$ for $i = 1,2$ have the same underlying analytic space $Y$. First we will show that the (global) regular functions on $Y_i$ are generated by the set of weight vectors in the space $\cO_Y(Y)$ of holomorphic functions on $Y$, showing that the algebraic structure on $Y$ is determined by the analytic structure on $Y$ and the $T$-action.

        So let $f \in \cO_Y(Y)$ be a holomorphic function that is a weight vector for the action of $T$, i.e.\ $(t \cdot f)(p) \coloneq f(t \cdot p) = \alpha(t) f(p)$ for some character $\alpha : T \to U(1)$ and all $p \in M$.
        Consider a $T$-equivariant algebraic embedding $Y_1 \inj \bbC^N$ as in \cite[Lem.~2.5]{ColSzeStab}, where $T$ acts on $\bbC^N$ linearly with positive weights. In a neighborhood of $o \in Y_1 \subseteq \bbC^N$, we can extend $f$ to a germ of a holomorphic function $\tilde{f}$ on $(\bbC^N, 0)$. Replacing $\tilde{f}$ by its projection on the weight-$\alpha$ space
        $ \int_T \alpha(t)^{-1} t \cdot \tilde{f} \d t$,
        we can assume that $\tilde{f}$ has weight $\alpha$.
        Now $\tilde{f}$ can be expanded as a power series around $0$, but since each coordinate function $z_i$ on $\bbC^N$ has positive weight $w_i > 0$, we see that the expansion can only have monomials of total degree $\leq \lambda / \min_i (w_i)$. Thus $f$ is algebraic on $Y_1$. Since all monomials on $\bbC^N$ are weight vectors, it is clear that the weight vectors in $\cO_Y(Y)$ generate the regular functions on $Y_1$.

    By the definition of a Fano fibration, the canonical divisors $K_{M_1}$ and $K_{M_2}$ are ($\bbQ$-Cartier and) $\pi$-ample, so for sufficiently divisible $m_i$ there are closed embeddings $j_i: M \to Y \times \bbP^{N_i}$ given by
    \[
        p \mapsto (\pi(p), [s_0(p) : \dots : s_{N_i}(p)])
    \]
    for some (algebraic) $s^{(i)}_j \in H^0(M, -m_i K_{M_1})$. Since the action of $T$ on $M_i$ is algebraic, every $s^{(i)}_j$ can be written as a sum of weight vectors in the $T$-representation $H^0(M, -m_i K_{M_i})$. Thus replacing $s^{(i)}_j$ by their weight components, we can assume that the $j_i$ are $T$-equivariant w.r.t.\ some linear diagonal $T$-action on $\bbP^{N_i}$.
    In particular, we have a $T$-equivariant algebraic isomorphism $j_1^* \pr_2^* \cO(1) \iso -m_1 K_{M_1}$. This also induces an isomorphism of the analytifications, so
    \[
        H^0(M, j_1^* \pr_2^* \cO(m_2)) \cong H^0(M, -m_1 m_2 K_M) \cong H^0(M, j_2^* \pr_2^* \cO(m_1))
    \]
    $T$-equivariantly ($H^0$ denotes analytic sections).
    \begin{claim*}
        For sufficiently large $k$, any weight vector in $H^0(M, j_1^* \pr_2^* \cO(k))$ is algebraic on $M_1$.
    \end{claim*}
    \begin{pf}
        For  sufficiently large $k$, the sequence of sheaves on $\bbC^N \times \bbP^{N_1}$
        \[
            0 \to \pr_2^* \cO(k) \otimes \cI_M \to \pr_2^* \cO(k) \to {j_1}_* j_1^* \pr_2^* \cO(k) \to 0
        \]
        will be exact on global sections, see \cref{thm:_vanish} below. Thus a weight vector $f \in H^0(M, -k m_1 m_2 K_M)$ has a preimage
        \[
            \tilde{f} \in  H^0(Y \times \bbP^N, \pr_2^* \cO(k m_2)) \cong \cO_Y(Y) \otimes H^0(\bbP^N, \cO(k m_2))
        \]
        which we can take to be a weight vector as well. Then we can write it as a sum of tensor products of weight vectors in $\cO_Y(Y)$ and $H^0(\bbP^N, \cO(k m_2))$. Since the weight vectors in $\cO_Y(Y)$ are algebraic, this shows that $\tilde{f}$ is algebraic and thus so is $f$.
    \end{pf}
    Hence when $m_2$ is sufficiently large, the embedding $j_2: M \inj Y \times \bbP^{N_2}$ is defined by (algebraic) sections of $-m_1 m_2 K_{M_1}$ and is thus algebraic on $M_1$. Since the algebraic structure of $M_2$ is induced by this embedding, we conclude $M_1 \cong M_2$.
\end{proof}

\begin{lem} \label{thm:_vanish}
    Let $\cF$ be a coherent algebraic sheaf on $\bbC^N \times \bbP^{N_1}$. Then the analytic cohomology group $H^1(\bbC^N \times \bbP^{N_1}; \cF \otimes \pr_2^* \cO(k))$
    vanishes for $k \gg 0$.
\end{lem}
\begin{proof}
    Since $\pr_2^* \Oalg(1)$ is an ample algebraic sheaf on $\bbP^{N_1}_{\bbA^N}$, the algebraic coherent sheaf $\cF \otimes \pr_2^* \Oalg(m_2)$ is generated by global (algebraic) sections. So there exists a coherent algebraic $\cF_1$ such that
    \[
        0 \to \cF_1 \to \cO^{\oplus m_1} \to \cF(n_1) \to 0
    \]
    is exact.
    Applying the same reasoning recursively to the kernel we obtain exact sequences,
    \begin{align}
        0 \to \cF_2 \to \cO^{\oplus m_2} \to \cF_1(n_2) \to 0 \\
        0 \to \cF_3 \to \cO^{\oplus m_3} \to \cF_2(n_3) \to 0 \p
        \\[-0.4em]
        \vdots \phantom{\to \cF_2(n_3) \to 0}
    \end{align}
    Hence for sufficiently large twists $k$, the long exact sequence in cohomology shows
    \begin{align}
        H^1\left(\bbC^N \times \bbP^{N_1}; \cF(k)\right) &= H^2\left(\bbC^N \times \bbP^{N_1}; \cF_1(k - n_1)\right)
            = \cdots \\
            &= H^{N + N_1 + 1}\left(\bbC^N \times \bbP^{N_1}; \cF_{N + N_1}(k - n_1 - n_2 -\cdots)\right) = 0
        \c
    \end{align}
    where the last cohomology group vanishes for dimension reasons.
\end{proof}

Finally we show that, for topological reasons, the exceptional fiber $E = \pi^{-1}(o)$ of a polarized Fano fibration $(\pi: M \to Y, T, \xi)$ is determined, up to bioholomorphism, by the analytic structure. Note that for any smooth AC shrinker except $\bbC^n$, this follows from the fact that $E$ is the unique connected compact analytic subset of $M$ of positive dimension.

\begin{prop} \label{thm:Euniq}
    Let $(M_i, T_i, \xi_i), i = 1, 2$ be polarized Fano fibrations with the same underlying analytic space $M$. At the level of analytic spaces, the fibration map $\pi: M \to Y$ agrees for both varieties (see \cref{thm:coneAnalytic}).
    Let $o_i$ be the unique fixed point of the action induced by $T_i$ on $Y$. Then there exists a biholomorphism $\psi: M \to M$ such that $\psi(\pi^{-1}(o_1)) = \pi^{-1}(o_2)$. In fact, $\psi$ lies in the subgroup of $\Aut(M)$ generated by the images of the action maps $T_i \to \Aut(M)$.
\end{prop}
\begin{proof}
    Since the Reeb cones of $(M, T_1)$ and $(M, T_2)$ are open, we can pick rational vectors $\tilde{\xi}_i \in \ft_i$, meaning that $\exp(\bbR \tilde{\xi}_i) \subseteq T_i$ is an $S^1$. This way we obtain two holomorphic actions of $\bbC^*$ on $M$ and therefore also on $Y$. We denote the actions on $Y$ by $\rho_i: \bbC^* \to \Aut(Y)$.
    We will show by contradiction that the subgroup $G \coloneq \langle \rho_1(\bbC^*), \rho_2(\bbC^*) \rangle \leq \Aut(Y)$ must contain an element $\phi$ such that $\phi(o_1) = o_2$.
    So we assume that the orbit $A = G \cdot o_1$ does not contain $o_2$.

    If we pick a holomorphic embedding $j: Y \inj \bbC^{N_1}$, equivariant for $\rho_1$ and some linear action on $\bbC^{N_1}$ with weights $w_i$,
    we see that $j(A)$ is contracted by the nullhomotopy
    \[
        j(A) \times I \to j(A), \quad (z, t) \mapsto (t^{w_1} z_1, \dots, t^{w_{N_1}} z_{N_1})
        \p
    \]
    Since $j$ is a homeomorphism onto $j(A)$, this shows $A$ is contractible.
    
    Now consider another embedding $i: Y \inj \bbC^N$, equivariant for $\rho_2$ and some linear action with weights $m_1, \dots, m_N > 0$. The stabilizer of a point $z \in \bbC^N$ under this action is
    \[
        \bigcap_{k \text{ s.t.\ } z_k \neq 0} \bbZ_{m_k} \leq U(1)
        \c
    \]
    where $\bbZ_q$ denotes the cyclic group of order $q$, viewed as a subgroup of $U(1) \leq \bbC^*$.
    Since $o_2 = i^{-1}(0) \not\in A$, every $z \in i(A)$ has a $\rho_2$-stabilizer contained in $\bbZ_M \leq U(1)$, for $M \coloneq \operatorname{lcm}(m_1, \dots, m_n)$. Pick $p \in \bbN$ coprime to $N$. Then $\bbZ_p \cap \bbZ_M = \{1\}$, so the action of $\bbZ_p \leq \bbC^*$ on $A$ via $\rho_2$ is free, making the quotient map $A \to A / \bbZ_p$ a principal $\bbZ_p$-bundle with contractible total space. In fact, we have a commuting diagram of principal $\bbZ_p$-bundles
    \begin{cd} \label{eq:diagBZp}
        A \arrow[r, "i"] \arrow[d] & \bbC^N \setminus \{0\} \arrow[r] \arrow[d]  & \bbC^\infty \setminus \{0\} \arrow[d] \simeq S^\infty \\
        A/\bbZ_p \arrow[r, ""] & (\bbC^N \setminus \{0\}) / \bbZ_p \arrow[r] & (\bbC^\infty \setminus \{0\}) / \bbZ_p \simeq B \bbZ_p
    \end{cd}
    where the $\bbZ_p$-action on $\bbC^\infty$ is taken to have weights $m_1, \dots, m_n, 1, 1, \dots$, so that the top vertical arrows are $\bbZ_p$-equivariant.

    Since $A$ is contractible, the long exact sequence in homotopy for a Serre fibration \cite[Thm.~4.41]{Hatcher} now shows that $\pi_2$ and all higher homotopy groups of $A / \bbZ_p$ vanish. For $\pi_1$, we obtain from the functoriality of the long exact sequence, that
    \begin{cd*}
        \pi_1(A/\bbZ_p) \arrow[d, "\sim" {anchor=south, rotate=90}] \arrow[r]  & \pi_1(B \bbZ_p) \arrow[d, "\sim" {anchor=south, rotate=90}] \\
        \pi_0(\bbZ_p) \arrow[r, equal]  & \pi_0(\bbZ_p)
    \end{cd*}
    commutes. Therefore the inclusion of $A/\bbZ_p$ in $B\bbZ_p$, which factors as
    \[
        A/\bbZ_p \to (\bbC^N \setminus \{0\})/\bbZ_p \to B \bbZ_p
        \c
    \]
    is a weak homotopy equivalence and therefore an isomorphism on cohomology \cite[Prop.~4.21]{Hatcher}. Now on the one hand the action of $\bbZ_p$ on $\bbC^N \setminus \{0\}$ commutes with the scaling action of $\bbR_{>0}$, so $(\bbC^N \setminus \{0\}) / \bbZ_p \simeq S^{2n-1} / \bbZ_p$ is (homotopy equivalent to) a compact manifold and thus has finite-dimensional cohomology.
    On the other hand, it is well-known
    (for example by the Leray--Hirsch theorem,
    or comparing with group cohomology)
    that the total cohomology group $H^*(B \bbZ_p; \bbQ)$ cannot be finite-dimensional. Thus the composite map
    \[
        H^*(A / \bbZ_p; \bbQ) \to H^*((\bbC^N \setminus \{0\})/\bbZ_p; \bbQ) \to H^*(B \bbZ_p; \bbQ)
    \]
    cannot be surjective, a contradiction.

    To conclude the proof, not that
    since $\rho_i$ factors through the "restriction map" $\Aut(M) \to \Aut(Y)$, the biholomorphism $\phi$ has a preimage $\psi \in \Aut(M)$ satisfying $\psi(\pi^{-1}(o_1)) = \pi^{-1}(o_2)$.
\end{proof}

\subsection{Conjugating tori in $\Aut(E^N)$}%
\label{sub:conjugating}

We consider a polarized Fano fibration $(\pi: M \to C_0, T, \xi)$. In particular $M$ could be the complex manifold underlying a Kähler--Ricci shrinker, $C_0$ its asymptotic cone and $T$ the torus generated by $\xi = JX$ where $X$ is the soliton vector field, cf.\ \cref{sub:defs}. Regardless of the origin of the polarized Fano fibration, we will define $X = - J \xi$.

We denote the vertex of the cone $C_0$ (the unique fixed point of $T \acts C_0$) by $o$ and call $E \coloneq \pi^{-1}(\{o\})$ the \demph{exceptional fiber}, with inclusion map $i: E \to M$.
Let $\IE \leq \cO_M$ be the sheaf of holomorphic functions vanishing on $E$ and let $E^N$ be the analytic scheme associated to $\IE^N$.

\begin{defn}
    We define the group of \demph{germs of automorphisms of $M$ at $E$} as
    \[
        \Aut(M)_E = \{f : U \to M \text{ biholomorphic onto its image} \mid
                        U \supset E \text{ open},~
                        f(E) = E
                    \} \,/ \sim
        \c
    \]
    where $f \sim g$ if{}f they agree in a neighborhood of $E$. There exists a restriction homomorphism $\Aut(M, E) \to \Aut(M)_E$ from the automorphisms of $M$ preserving $E$ and it is injective by the identity theorem.
\end{defn}

We can push forward germs of vector fields $V \in \Re H^0(E, i^{-1}\cT_M)$ along germs of automorphisms $\eta \in \Aut(M)_E$ by picking an open $U \supset E$ on which $\eta$ and $\eta^{-1}$ are defined and on which $V$ can be represented by $\hat{V} \in \Re H^0(U, \cT_M)$ and then setting
\[
    \eta_* V = \left[ \d \eta \hat{V} \circ \eta^{-1} \right]
    \p
\]

By \cite[Satz~7]{Kaup}, the group $H^N = \Aut(E^N)$ of holomorphic automorphisms of the \emph{compact} analytic scheme $E^N$ is a Lie group with Lie algebra $\fh^N = \Hom(\Omega^1_{E^N}, \cO_{E^N}) \cong \operatorname{Der}(\cO_{E^N}, \cO_{E^N})$.

Consider the subsheaf $\cG \leq \Re (\cT_M)$ consisting of real-holomorphic vector fields that infinitesimally preserve $E$, or equivalently preserve $\IE$. Concretely this means that $V \in \cG(U)$ if{}f for every holomorphic function $f$ defined on some open $V \subseteq U$
\[
    f \in \IE(V) \implies V^{1,0} f \in \IE(V)
    \p
\]
Since these vector fields also preserve $\IE^N$ they act as derivations on $\cO_{E^N}$, giving us the restriction map $\rho_*: \cG(U) \to \fh^N$, at least when $U \supseteq E$. The map $\rho_*$ should be thought of as the map induced on Lie algebras by the restriction homomorphism
\[
    \rho: \Aut(M)_E \to \Aut(E^N)
    \p
\]
The kernel of $\rho_*$ consists of the derivations that are zero on $\cO/\IE^N$, i.e.
\begin{equation} \label{eq:kerrho}
    \cK \coloneq \ker \rho_* = \sHom(\Omega^1_M, \IE^N)\leq \sHom(\Omega^1, \cO_M) = \cT_M
    \c
\end{equation}
or in other words, the vector fields that restrict to zero on $E^N$.

Note that $\cG(M) = \aut(M, E)$ consists of all global real-holomorphic vector fields on $M$ that preserve $E$.
We define $G^N \leq \Aut(E^N)$ to be the immersed Lie subgroup with Lie algebra
\[
    \fg^N \coloneqq \Im\left(\rho_*: \aut(M, E) \to \fh^N\right)
    \p
\]
The next lemma gives another, more concrete, characterization of $G^N$.

\begin{lem} \label{thm:GNimage}
    Every $g \in G^N$ is the restriction to $E^N$ of some $\hat{g} \in \Aut(M)_E$.
\end{lem}
\hyperref[not:Exp]{Recall} that we write $\Exp(V)$ for the time 1 flow of a vector field $V$, as opposed to lowercase $\exp$, which will always denote the exponential map of a Lie group
\begin{proof}
    A Lie group is generated by the image of its Lie algebra under the exponential map, so $g = \exp(\rho_* V_1) \cdots \exp(\rho_* V_n)$ for some vector fields $V_i \in \aut(M, E)$. While the time 1 flow of the vector fields $V_i$ may not be defined globally, there exists a neighborhood of $E$ on which $\hat{g} \coloneq \Exp(V_1) \circ \cdots \circ \Exp(V_n)$ is defined, so $\hat{g} \in \Aut(M)_E$. 
    Finally by \cite[Satz~2(c)]{Kaup} we have $\rho(\Exp(V_i)) = \exp(\rho_* V_i)$, showing that $\rho(\hat{g}) = g$.
\end{proof}

We set $\fg_T \coloneq \aut(M, E)^T$ to be the algebra of $T$-invariant holomorphic vector fields on $M$ that preserve $E$. It contains the Lie algebra $\ft$ of $T \leq \Aut(M, E) \leq \Aut(M)_E$.

\begin{lem} \label{thm:liftT}
    Let $T' \leq G^N$ be a torus subgroup whose Lie algebra contains $\rho_* \ft$. Then $\ft' \subseteq \rho_* \fg_T$.
\end{lem}
\begin{proof}
    Let $y \in \ft' \leq \fg^N$. By definition of $\fg^N$ there exists $\hat{Y} \in \aut(M, E)$ such that $\rho_* \hat{Y} = y$. Using the (volume 1) Haar measure on $T$ we can average $\hat{Y}$ to
    \[
        \bar{Y} = \int_T \left(t_* \hat{Y}\right) \d t \in \fg_T
        \p
    \]
    Since $y$ is invariant for the adjoint action of $\rho(T)$, we have $\rho_* (t_* \hat{Y}) = \Ad(\rho(t)) \rho_*(\hat{Y}) = y$. Thus $\rho_* \bar{Y} = y$, proving the claim.
\end{proof}

The key step in the proof of the uniqueness (up to conjugation) of soliton vector fields is the following result stating that soliton vector fields are stable under perturbations of sufficiently high order:

\begin{prop} \label{thm:AutGermsv2}
    Let $(\pi: M \to C_0, T, JX)$ be a polarized Fano fibration with exceptional fiber $E \subseteq M$.
    Then there exists $N > 0$, independent of $X \in J\ft$, with the following property:
    For every vector field $H$ on some neighborhood $U$ of $E$ that vanishes to $N$-th order on $E$, meaning $H \in H^0(U, \cK^N)$, the germs of the vector fields $X$ and $X - H$ along $E$ are conjugate by some $\eta \in \Aut(M)_E$, i.e.\ $\eta_* (X - H) = X$.
\end{prop}

\begin{proof}
    The main fact about polarized Fano fibrations that we exploit is simply that there exists a nice Lyapunov function for $X$, i.e.\ a continuous function $r: M \to \bbR_{\geq 0}$ that is strictly increasing along the flow of $X$. To see this observe that by \cite[Lem.~2.5]{ColSzeStab} there is an embedding of cones
    \[
        f: (C_0, X) \inj \left( \bbC^N, V_a \right)
        \qtext{where}
        V_a = 2 \Re \left(\sum_i a_i z_i \frac{\partial}{\partial z_i}\right)
    \]
    and $a = (a_1, \dots, a_N) \in \bbR^N$ are the $X$-weights of some generators $f_1, \dots, f_N$ of $\Oalg(C_0)$.
    The continuous function $\hat{r} = \sum_{i = 1}^N |z_i|^{1/a_i}$ has weight 1 under the flow of $V_a$. Thus
    \begin{equation} \label{eq:rAlgebraic2}
        r = \hat{r} \circ f = \sum_{i = 1}^N |f_i|^{1/a_i}
    \end{equation}
    has the desired property $r \circ \gamma_{2t} = e^{t} r$, where we use the notation $\gamma_{2 t} \coloneq \Exp(t X)$ from \cref{sec:functionals}.
    We also still use the notation $B_s \coloneq \{r < s\}$.
    Since $E = r^{-1}(0)$ is compact, we can replace $U$ by $B_\delta$ for some $\delta > 0$.
    Consider the family
    \[
        \eta_t \coloneqq \Exp\left(t X\right) \circ \Exp\left(-t (X - H)\right)
        \c
    \]
    of biholomorphisms of the analytic space $M$,
    which are defined on an open neighborhood of $E$ on which the flow of $-(X - H)$ exists up to time $t$ (the vector field $X$ on the other hand is complete since it comes from a torus action). Note that $\eta_t|_E = \id_E$.
    We will show that $\eta = \lim_{t \to \infty} \eta_t$ exists in the compact-open topology on some neighborhood of $E$.

    The time derivative of $\eta_t$ is%
    \footnote{
        We can compute the time derivative of $\eta_{t_0}$ at a point $p$ by writing
        \[
            \eta_{t_0 + s} = \Exp(s X) \circ \Exp(t_0 X) \circ \Exp(-s (X - H)) \circ \Exp(-t_0 (X - H))
            \p
        \]
        Then we find neighborhoods $W_1, W_2$ of $\Exp(-t_0(X-H))(p)$ and $\eta_{t_0}(p)$ that embed into open sets $V_1, V_2 \subseteq \bbC^m$ such that there are extensions of $-(X - H)$ and $X$ to $V_1$ and $V_2$ respectively. Shrinking $W_1$ and $W_2$, we can also extend $\Exp(t_0 X)$ to a map from $V_1$ to $V_2$. At this point the left three maps are smooth maps between manifolds and the claimed formula follows from the chain rule.
    }
    \begin{align}
        \dot{\eta}_t
            &= X \circ \eta_t - \d\! \Exp(t X) (X - H) \circ \Exp(-t (X - H)) \\
            &= \d\! \Exp(t X) H \circ \Exp(-t X) \circ \eta_t
            \p
    \end{align}
    In other words, $\eta_t$ is the flow of the time-dependent vector field
    \[
        H_t \coloneqq (\gamma_{2t})_* H = \d \gamma_{2t} H \circ \gamma_{-2t}
        \c
    \]
    where $\gamma_{2t}$ still denotes the time $t$ flow of the vector field $X$.
    $H_t$ is defined on the open neighborhood $\{r < e^t \delta\}$ of $E$.

    We fix an embedding $j: M \to \bbP^m$ of the quasi-projective variety $M$ and use the restriction of some metric $g$ on $\bbP^m$ as a background metric.
    Around every $p \in M \subseteq \bbP^m$ we can find holomorphic extensions of $X$ to a neighborhood in $\bbP^m$. Combining them with a smooth partition of unity, we obtain a smooth extension $\tilde{X}$ of $X$ to a neighborhood of $E$ in $\bbP^m$.
    By compactness of $E$, there exists $C > 0$ such that $\cL_{\tilde{X}} g \leq C g$ on $E$, and hence also on some neighborhood $W$ of $E$ in $\bbP^m$.
    Shrinking $\delta$, we can assume that $B_\delta \subseteq M \cap W$, so that
    $\gamma_{2t}^* g \leq e^{2 C t} g$ on $B_{e^{2t} \delta}$.  Then for $z \in B_\delta$
    \begin{equation} \label{eq:compareMetricv2}
        |H_t(z)|_g \leq e^{C t} |H(\gamma_{-2t}(z))|_g
        \p
    \end{equation}

    We have three ways of quantifying "degree of vanishing at $E$" at our disposal: Being contained in a power of $\IE$, being a multiple of a power of $r$ and being a multiple of a power of the distance function $d_E$ induced by the metric $g$. The next lemma lets us compare these degrees:
    \begin{lem} \label{thm:compareVanish2}
        Let $C^0$ denote the sheaf of continuous functions on $M$.
        \begin{thmlist}
            \item There exists $k \in \bbN$ such that $\IE^k \subseteq r C^0$. As a consequence, for any $a \in \bbN$, a section of $\cK^{a k} \subseteq j^* \cT_{\bbP^m}$ can be written as $r^a$ times a continuous section of $j^* T \bbP^m$.
            \item There exists $K > \bbN$ such that $r^K \leq d_E$ in a neighborhood of $E$, where $d_E$ is the distance to $j(E)$ in $\bbP^m$, measured (extrinsic to $M$) with the Fubini--Study metric.
        \end{thmlist}
    \end{lem}
    \begin{proof}
        \point1
        Pick $\bbZ \ni \ell > \max_i \left(a_i^{-1}\right)$, where $a_i$ are the numbers from \cref{eq:rAlgebraic2}. Then for any multi-index $\alpha$ of order $|\alpha| = \ell$ we have that
        \[
            |f^\alpha| \leq \Big(\max_i |f_i|\Big)^\ell \leq \sum_{i = 1}^N |f_i|^\ell = o(r)
            \c
        \]
        so $f^\alpha C^0 \subseteq r C^0$.
        Since
        $\sqrt{(f_1, \dots, f_N)} = \IE$, there exists $m > 0$ such that 
        $\IE^m \subseteq (f_1, \dots, f_N)$.
        Hence
        \[
            \IE^{m \ell} \subseteq (f_1, \dots, f_N)^\ell \subseteq r C^0
            \c
        \]
        and we can pick $k = m \ell$.

        From the surjection $j^* \Omega^1_{\bbP^m} \surj \Omega^1_M$ and the left-exactness of $\sHom$ we obtain
        \[
            \cK^{a k} = \sHom(\Omega^1_M, \IE^{a k}) \subseteq \sHom(j^* \Omega^1_{\bbP^m}, \IE^{a k}) \subseteq \sHom_{\cO_M}(j^* \Omega^1_{\bbP^m}, r^a \cdot C^0) = r^a C^0(-; j^* T \bbP^m)
            \c
        \]
        where $C^0(-; j^* T \bbP^m)$ denotes the sheaf of continuous sections of $j^* T \bbP^m$. For the last step we can use that $j^* \Omega^1_{\bbP^m}$ is locally free.

        \point2
        Starting from \cref{eq:rAlgebraic2} we see that if we pick an integer $K \geq \sum_i a_i$, then by pigeonhole
        \[
            r^K \in (|f_1|, \dots, |f_N|)_{C^0}
            \c
        \]
        where the RHS is the subsheaf of $C^0$ generated by the functions $|f_i|$.
        Now we pick finitely many open subsets $W_\alpha \subseteq \bbP^m$ covering $E$, to which the $f_i$ can be extended holomorphically. Shrinking $W_\alpha$, we can assume that the $|f_i|$ are Lipschitz with some constant $L \in \bbR$. Since we can locally write $r^K = \sum_i g_i |f_i|$ for some continuous functions $g$, we conclude
        \[
            r^K|_{W_\alpha} \leq L \left(\sup_{V_i} |g_i|\right) d_E
            \p
        \]
        After possibly shrinking $W_\alpha$ we obtain that $r^K \leq B d_E$ near $E$ for some $B > 0$. Replacing $K$ by $K+1$ we can take $B = 1$ in a neighborhood of $E$.
    \end{proof}
    Back to the proof of \cref{thm:AutGermsv2}, we will assume w.l.o.g.\ that $C > K$.

    Now we require $N > (C+1) k$ so that by \point1 of the Lemma
    \[
        H \in \Re H^0\left(B_\delta, \cK^N \right) \subseteq r^{C+1} C^0(B_\delta; j^* T \bbP^m)
    \]
    i.e.\ $H = r^{C+1} F$ for some continuous vector field $F$ on $B_\delta$. Continuing from \cref{eq:compareMetricv2} we have that on $B_\delta$
    \begin{align}
        | H_t |_g
            &\leq e^{C t} | H \circ \gamma_{-2t} |_g \\
            &\leq e^{C t} | F \circ \gamma_{-2t} |_g e^{-(C+1) t} r^{C+1}
             = e^{-t} r^{C + 1} | F \circ \gamma_{-2t} |_g \\
            &\leq e^{-t} r d_E \| F \|_{g, L^\infty(B_\delta)} \\
            &< e^{-t} d_E
            \p
    \end{align}
    For the last step we assume $\delta < \| F \|_{g, L^\infty(B_\delta)}^{-1}$.
    Now consider a flow line $x: [0, T) \to B_\delta \subseteq \bbP^m$ of $H_t$ and let $s(t) = d(E, x(t))$. Then
    \[
        \dot{s}(t) \leq |H_t(x(t))|_g < e^{-t} s(t)
        \c
    \]
    so by ODE comparison $s(t) < s(0) \exp(1 - e^{-t}) < e^1 s(0)$.
    This shows that flow lines of $H_t$ starting in $\bar{B}_{e^{-1} \delta}$ stay in $B_\delta$ and thus exist for all times. Furthermore the length of their derivatives is bounded by $\delta e^{-t}$, showing that
    \[
        \eta_t|_{\bar{B}_{e^{-1} \delta}}: \bar{B}_{e^{-1} \delta} \to B_\delta
    \]
    is Cauchy in the compact-open topology and therefore has a limit
    \[
        \eta = \lim_{t \to \infty} \Exp(t X) \circ \Exp(- t(X-H))|_{\bar{B}_{e^{-1} \delta}}
        \p
    \]
    As the compact-open limit of holomorphic maps, $\eta$ is holomorphic \cite[Thm.~H.15]{GunII}.
    By construction
    \begin{align}
        \eta \circ \Exp(s(X - H)) &= \lim_{t \to \infty} \Exp(t X) \circ \Exp((s-t)(X-H)) \\
                           &= \lim_{t \to \infty} \Exp((t + s) X) \circ \Exp(t (X-H)) \\
                           &= \Exp(s X) \circ \eta
    \end{align}
    on $B_{e^{-1} \delta}$. Taking $s$-derivatives we see that $\eta_* (X - H) = X$. Since all $\eta_t$ fix $E$ the same is true of $\eta$, so $\eta$ represents an element of $\Aut(M)_E$.
\end{proof}

The concept of a maximal torus makes sense in the automorphism group of a \emph{noncompact} analytic space:
We say a faithful torus action $T \acts M$ is \emph{maximal} if there exists no torus $T' \geq T$ to which the action can be extended faithfully. A faithful action of a torus $\bbT^m$ on an $n$-manifold always factors through the isometry group of a $\bbT^m$-invariant Riemannian metric, so $m \leq n^2$. The same bound holds for analytic spaces, since a faithful action restricts to a faithful action on the regular locus. Hence there are no unbounded chains of tori in the automorphism group of an analytic space, so a torus action $\bbT \acts M$ can always be extended to a maximal torus action.

Note that on a polarized affine variety $(T \acts Y, \xi)$ with vertex $o$ we have $\lim_{t \to \infty} \Exp( t J \xi )(p) = o$ for any $p \in Y$.
Thus if $T' \geq T$ is a bigger torus acting faithfully on a polarized Fano fibration $(M, T, JX)$, the induced action on $Y$ still fixes $o$, and the triple $(M, T', JX)$ still is a polarized Fano fibration.

Using \cref{thm:AutGermsv2} we can prove Iwasawa's theorem on conjugate tori for the automorphism group of a polarized Fano fibration, at least for tori that contain a Reeb vector field:

\begin{prop} \label{thm:pffTori}
    Let $(\pi: M \to C_0, T_1, X_1)$ and $(\pi: M \to C_0, T_2, X_2)$ be two polarized Fano fibrations with maximal tori $T_i$ and with the same underlying complex analytic space $M$.
    Then there exists a biholomorphism $\psi: M \to M$ such that $\psi T_1 \psi^{-1} = T_2$.
\end{prop}
\begin{proof}
    Replacing $X_1$ by a general vector in the Reeb cone we can assume that $J X_1$ generates $T_1$. Then it suffices to show that $\psi_* X_1 \in J\ft_2$.
    Let $E = \pi^{-1}(\{o\})$ be the exceptional fiber of $\pi$ over the vertex of $C_0$, which by \cref{thm:Euniq} we can assume is the same for both polarized Fano fibrations. Observe that the definition of the group $G^N$ only depends on the analytic structure of $M$ and on the subvariety $E$. Thus we have two tori $\rho(T_1)$ and $\rho(T_2)$ in $G^N$, with $\rho_* J X_i \in \Lie \rho(T_i)$. First we extend these to maximal tori $\rho(T_i) \leq T_i' \leq G^N$.
    Maximal tori in any connected Lie group are conjugate (because maximal compact subgroups of a connected Lie group are conjugate \cite{IwasawaCpt}, and maximal tori in compact groups are conjugate), so there exists some $g \in G^N$ such that $gT_1'g^{-1} = T_2'$.

    By \cref{thm:GNimage}, $g$ is the restriction of some $\hat{g} \in \Aut(M)_E$ defined on some open neighborhood $U$ of $E$, hence $\Ad(g) (\rho_* J X_1) = \rho_* (\hat{g}_* J X_1)$.
    On the other hand, \cref{thm:liftT} shows that there exists a vector field $J Y \in \fg_{T_2}$ such that $\Ad(g) (\rho_* J X_1) = \rho_* J Y$.
    By \cref{eq:kerrho} this means that
    \[
        X_1 - \hat{g}^* Y \in H^0\left(U, \cK^N\right)
        \p
    \]
    Then by \cref{thm:AutGermsv2}, there further exists $\eta \in \Aut(M)_E$ such that $\eta_* \hat{g}^* Y = X_1$. Pick $\epsilon > 0$ small enough so that $\hat{g} \eta^{-1}$ is defined on $B_\epsilon$ (here $B_\epsilon = \pi^{-1}(\{r_1 < \epsilon\})$ for some radial function $r_1$ on the cone $T_1 \acts C_0, X_1$).

    By the same argument as in \cite[Lem.~2.34]{CDS}, the vector field $Y \in \Re H^0(M, \cT_M)$ is complete since it commutes with $X_2$.
    Thus $\hat{g} \eta^{-1}$ can be extended to a $\psi \in \Aut(M)$ such that $\psi_* X_1 = Y$ on all of $M$ by setting
    \[
        \psi = \lim_{t \to \infty} \Exp(t Y) \circ (\hat{g} \eta^{-1}) \circ \Exp(-t X_1)
        \p
    \]
    Note that this limit is actually constant in $t$ where defined and its domain contains $B_{e^t \epsilon}$, thus $\psi \in \Aut(M)$.

    Since $J X_1$ generates $T_1$, the flow of $Y$ generates $\psi T_1 \psi^{-1}$. This torus commutes with $T_2$ since $Y$ does. Since both tori are maximal, we must have $\psi T_1 \psi^{-1} = T_2$.
\end{proof}
\subsection{Applications}%
\label{sub:applications}

\Cref{thm:pffTori} lets us apply the volume minimization theorem for shrinkers \cite[Thm.~5.19]{CDS} to any two vector fields, even if they do not commute. Thus we conclude that the soliton vector field of a shrinker with Ricci curvature decay is unique up to automorphism, which finishes the proof of our main theorem:
\uniq*

\begin{proof}
    Let $T_i$ denote the torus generated by $J X_i$ and let $T'_i \geq T_i$ be a maximal torus.
    Applying \cref{thm:pffTori} we obtain a shrinker $(\psi_* X_1, \psi_* \omega_1)$ whose vector field lies in $\ft'_2$.
    By the volume minimization principle \cite[Thm.~5.19]{CDS} we conclude that $\psi_* X_1 = X_2$. Finally, the uniqueness of the metric for a given soliton vector field (under these curvature decay hypotheses) was established in \crefrange{sec:geodesics}{sec:uniq_along}.
\end{proof}

Since \cite[Thm.~5.19]{CDS} only requires the Ricci curvature to be bounded, the argument above proves the uniqueness of the soliton vector field (but not the metric) under slightly weaker assumptions:

\uniqVF*

With the same arguments we also obtain a uniqueness theorem for Calabi--Yau cone metrics, using the volume minimization principle of \cite{MSYSasaki}.

\uniqCY*
\begin{proof}
    Let $(L, g_2)$ be the link of $(Y, \omega_2)$. We may consider $\Isom(L, g_2)$ as a subgroup of the isometry group of $\omega_2$.
    Let $T_2$ be a maximal torus in $K \coloneq \Isom(L, g_2) \cap \Aut(Y)$. By \cite[Thm.~1]{HeIso}, $K$ is a maximal compact subgroup of $\Aut(Y)$.
    Furthermore $\Aut(Y, X_2)$ is a Lie group by \cite[Prop.~4.9]{DonSunII}, and it contains $K$. Since maximal compact subgroups are conjugate in Lie groups this shows that $T_2$ is a maximal torus in $\Aut(Y, X_2)$. And since $X_2 \in J \ft_2$, any torus in $\Aut(Y)$ containing $T_2$ lies in $\Aut(Y, X_2)$, so $T_2$ is a maximal torus in $\Aut(Y)$.
    Hence by \cref{thm:pffTori} we can assume $X_1 \in \ft_2$. Then by \cite[\S 3.3]{MSYSasaki}, $X_1 = X_2$ and by \cite[Prop.~4.8]{DonSunII} $\omega_1$ and $\omega_2$ agree up to a biholomorphism preserving $X_2$.
\end{proof}

Using \cref{thm:pffTori} we can also show that the algebraic structure of a polarized Fano fibration is uniquely determined by the complex analytic space $M$.
\begin{cor} \label{thm:pffalg}
    Let $(\pi_1: M_1 \to Y_1, T_1, X_1)$ and $(\pi_2: M_2 \to Y_2, T_2, X_2)$ be polarized Fano fibrations with the same underlying analytic space $M$.
    Then there exists an isomorphism (of algebraic varieties) $\psi: M_1 \to M_2$.
\end{cor}
\begin{proof}
    W.l.o.g.\ we may assume that both tori $T_i$ are maximal. 
    By \cref{thm:coneAnalytic} the maps $\pi_1$ and $\pi_2$ are the same map $\pi: M \to Y$ at the level of analytic spaces. Applying \cref{thm:pffTori} we obtain a biholomorphism $\psi \in \Aut(M)$ such that $\psi_* X_1 \in J \ft_2$. Since the characterization of Reeb vector fields in \cref{thm:Reebflow} does not rely on the algebraic structure (and is biholomorphism-invariant), we see that $\psi_* X_1$ is in the Reeb cone of $(M_2, T_2)$, so $(M_2, T_2, \psi_* X_1)$ is a polarized Fano fibration. Then \cref{thm:vfDeterminesAlg} shows that $\Oalg_{M_2} = \psi_* \Oalg_{M_1}$.
\end{proof}

\begin{appendix}
\section{A regularization theorem}%
\label{apx:a_regularization_theorem}

\newcommand{\BandK}{B{\l}ocki--Ko{\l}odziej{}}
We adapt \BandK's elementary regularization technique \cite{BK} to our setting.
We write $\Omega = (0, 1) \times i \bbR, \Omega_\nu = (1/\nu, 1 - 1/\nu) \times i \bbR$ and $\pr_M$ for the projection $\Omega \times M \to M$.
\begin{prop} \label{thm:regularization}
    Assume $(\cH_2^*)^T \neq \emptyset$.
    Let $\Phi$ be a locally bounded $i \bbR$-invariant, PSH metric on $\pr_M^* K_M^{-1}$ over $\bar{\Omega} \times M$ such that $\Phi_t \in \Hsing(a,b)^T$ for every $t$.
    Then there exists a sequence $\Phi_\nu$ of smooth $i \bbR$-invariant PSh metrics on $\pr_M^* K_M^{-1}$ over $\overline{\Omega_\nu} \times M$ non-increasing pointwise to $\Phi$, such that for every $t$:
    \[
        \Phi^\nu_t \in \cH_2(a', b')^T
        \qtext{and}
		i \partial\bar{\partial} \Phi^\nu_t > \epsilon(\nu) \pr_M^* \omega_C
    \]
	for some cone metric $\omega_C$ and $a' = a'(a, b), b' = b'(a, b)$.
    Furthermore in equivariant charts, all (spatial \emph{and} temporal) derivatives of $\Phi^\nu$ are $O(r^2) = O(e^{2x})$, although not uniformly in $\nu$.

    If furthermore $\Phi$ is locally Lipschitz in $t$ with Lipschitz "constant" $L (1 + r^2)$ then so is $\Phi^\nu$.
\end{prop}

\begin{proof}[Proof sketch]
    We will first construct $\Phi^\nu$ satisfying all the claimed properties except $T$-invariance. Then at the end we explain why replacing $\Phi^\nu$ by its $T$-average preserves all claimed properties.

    We pick a background metric $\omega_0 = i \partial\bar{\partial} \phi_0$ with $\phi_0 \in (\cH_2^*)^T$. We will pick $r^2/2$ to be the potential of the cone metric that $i \partial\bar{\partial} \phi_0$ is asymptotic to.

    We claim that applying the regularization method described in the proof of \cite[Thm.~2]{BK} on the manifold $\Omega \times M$, with $\gamma = \pr_M^* \omega_0$
    to the $\gamma$-PSH functions $(t, p) \mapsto \Phi(t, p) - \phi_0(p)$,
    using the charts $\{U_\alpha\}_\alpha = \cA_M$
    will result in functions $\Phi^\nu$ (called $\phi_j$ in \cite{BK}) satisfying the claimed properties.

    We will now go through the steps of the proof in \cite{BK} that need to adapted for our version. Unless otherwise noted we use the same notations and definitions as \cite{BK}.
    We will pretend that $\cA$ only consists of equivariant charts, since the proof does not need to be modified for the precompact charts (recall these are the charts that we use to cover the exceptional set $E$). Furthermore, we will only consider the subset $(1, \infty) \times i \bbR \times U'_\alpha$ of every equivariant chart, since we can assume that the precompact charts, cover the remaining part of $U_\alpha$.
    For the subsets $V_\alpha \subseteq U_\alpha$ of \cite{BK} we pick $\overline{V'_\alpha} \subseteq U_\alpha'$ and set $V_\alpha = (1, \infty) \times i(-1, 1) \times V_\alpha'$.

    In the first step of the proof of \cite[Thm.~2]{BK} we can pick $f_\alpha$ to be equal to the PSH function corresponding to the background metric $\phi_0$ in that chart, so in particular $\epsilon = 0$ and $u_\alpha$ is just the function corresponding to the metric $\Phi$ in the chart $U_\alpha$.
    Thus the quadratic bounds for the metric $\Phi$ (from the assumption $\Phi_t \in \Hsing(a, b)^T$) translate into
    \[
        a e^{2x_\alpha} - a^{-1} \leq u_\alpha \leq b e^{2x_\alpha} + b
    \]
    in every equivariant chart $U_\alpha$, where we just write $x_\alpha$ for the first (real) coordinate in the chart $U_\alpha$.
    We might as well simplify these bounds to
    \[
        a e^{2x_\alpha} \leq u_\alpha \leq b e^{2 x_\alpha}
    \]
    since we can change $b$, and add constants to $u_\alpha$ by rescaling $\phi_0$.

    Performing convolutions with a smooth test function in every chart we obtain $u_{\alpha, \delta} \in C^\infty(V_\alpha \times \Omega_\nu)$, as long as $\delta$ is sufficiently small so that $\delta$-balls around points in $V_\alpha \times \Omega_\nu$ fit inside $U_\alpha \times \Omega$.

    \proofstep{Step 1: Adapting \cite[Lem.\ 4,~5]{BK}}
    The crucial point for us is that we can modify the proof of \cite[Lem.\ 4, 5]{BK} to obtain an analogue of \cite[Eq.~(3)]{BK}, namely that
    \begin{equation} \label{eq:unifConvr2}
        \epsilon(\delta) \coloneq \max_{\alpha, \beta} \sup_{V_\alpha \cap V_\beta} e^{-2x_\alpha}|(u_{\alpha, \delta} - u_{\beta, \delta}) - (f_\alpha - f_\beta)| \to 0
        \p
    \end{equation}
    where $x_\alpha$ is the first coordinate of the equivariant chart $U_\alpha$.
    Changing charts changes $x$ by at most adding a constant, so we will not keep track of which $x_\alpha$ we are using.
    We call this kind of convergence "$e^{2x}$-uniform".
    There are two places in the proof of \cite[Lem.\ 4]{BK} where the authors obtain locally uniform convergence. We will show that in both cases we can get $e^{2x}$-uniform convergence on $V_\alpha \cap V_\beta$.

    \textbf{First place.\quad}%
    \BandK{} use the zero Lelong number assumption and Dini's theorem
    to conclude that for any $r > 0$
    \[
        \frac{\hat{u}_r - \hat{u}_\delta}{\log r - \log \delta} \to 0
    \]
    locally uniformly as $\delta \to 0$. Instead of invoking Dini's theorem we will capitalize on the fact that our functions don't just have zero Lelong number, but are actually "$e^{2x}$-bounded", meaning that $|u| \leq e^{2x}$ for some constant $A$. Therefore
    \[
        \frac{\hat{u}_r - \hat{u}_\delta}{\log r - \log \delta} \leq 2b e^{2x} \frac{1}{\log r - \log \delta}
    \]
    so the RHS goes to zero $e^{2x}$-uniformly in $\delta$. A very similar argument is also used in \cite[Lem.\ 5]{BK}, but for the regularization $\tilde{u}_\delta$ instead of $\hat{u}_\delta$. We can use the same argument to get $e^{2x}$-uniform convergence there.

    \textbf{Second place.\quad}%
    \BandK{} use that the transition functions between charts $V_\alpha$ and $V_\beta$ only change the volume of small balls by at most some fixed factor.
    This is still true in our setting because the chart transition function $F = F_{\alpha, \beta}$ is equivariant with respect to
    the vector field $X = \partial_{x_\alpha} = \partial_{x_\beta}$ and the sets $L \cap V_\alpha$ on the link $L$ are precompact.
    Hence there is a finite upper bound $A$ on the Jacobian of $F$ and $F^{-1}$ and we still conclude
    \[
        \bar{B}(F(z), \delta) \subseteq F(\bar{B}(z, A\delta))
        \qtext{and}
        F(\bar{B}(z, \delta)) \subseteq \bar{B}(F(z), A \delta)
    \]
    for $z \in V_\alpha$ and $\delta$ sufficiently small as in the end of the proof \cite[Lem.\ 4]{BK}.

    \proofstep{Step 2: Adapting \cite[Thm.~2]{BK}}

    Once we have established \cref{eq:unifConvr2},
    we pick functions $\eta_\alpha \in C^\infty\left(U_\alpha \cap L; [-1, 0]\right)$ on the link $L$ such that $\eta_\alpha \equiv 0$ on $V_\alpha \cap L$ while $\eta_\alpha \equiv -1$ outside of some compact subset of $L \cap U_\alpha$. We then extend $\eta_\alpha$ to all of $U_\alpha$ by requiring $X(\eta_\alpha) = 2 \eta_\alpha$. We thus still have $i \partial\bar{\partial} \eta_\alpha|_{L \cap U_\alpha} > -C \omega_0$ for some $C > 0$ and since $i \partial\bar{\partial} \eta_\alpha$ also has weight 2 we obtain $| \partial\bar{\partial} \eta_\alpha |_{\omega_0} > -C $. We may assume $C > 3$

    Continuing as in \cite{BK} we finally set
    \begin{align}
        \label{eq:BKpatchRegularization}
        m_\delta &\coloneq \widetilde{\max_\alpha}^\kappa \left( u_{\alpha, \delta} - f_\alpha + \epsilon(\delta) \eta_\alpha /C \right) \\
        \label{eq:Phidelta}
        \Phi_\delta &\coloneq \phi_0 + \frac{1}{1 + 3 \epsilon(\delta)} m_\delta
    \end{align}
    using a regularized max function $\tilde{\max}^\kappa$ (cf.\ \cite[0.5.18]{agbook}).
    From the definition of $u_\alpha$ and $f_\alpha$ it is clear that $\Phi_\delta \to \Phi$ pointwise as $\delta \to 0$. It is also clear that for $\delta$ small enough $\frac{a}{2} e^{2 x_\alpha} \leq u_{\alpha, \delta} \leq 2 b e^{2 x_\alpha}$, and since each $f_\alpha$ also satisfies a bound of the form $f_\alpha \geq a_0 e^{2 x_\alpha} - a_0^{-1}$ we can easily conclude that $\Phi_\delta \in \cH_2(a', \infty)$ for some $a'$ depending on $a$ and $b$.

    By the definition of $\epsilon(\delta)$ in \cref{eq:unifConvr2}, the sets
    where $\eta_\alpha = -e^{2x_\alpha}$ (for $\alpha \leq N$) do not contribute to the regularized max when the regularization parameter $\kappa = \lVert \max - \tilde{\max}^\kappa \rVert_{\infty}$ is less than $\epsilon(\delta)/3$ (here we use that $-e^{2x_\alpha} \leq -1$ and $C > 3$).

    We obtain bounds on curvature
    \[
        \iddball (u_{\alpha, \delta} - f_\alpha + \epsilon \eta_\alpha / C) \geq 0 - (1 + \epsilon(\delta)) \pr_M^* \omega_0
    \]
    and conclude $\iddball \Phi_\delta \geq 0$ and $i \partial\bar{\partial} \Phi_\delta > \epsilon(\delta) \pi^* \omega_0$ (at least when $\epsilon(\delta) < 1/3$).

    To obtain our sequence $\Phi^\nu$ we pick a sequence $\delta_\nu \to 0$, set $\kappa = \epsilon(\delta_\nu)/4$, but are left with the task of fixing the fact that $\Phi_{\delta_\nu}$ might not be monotone decreasing. To simplify things assume that $\epsilon(\delta)$ is decreasing as $\delta \to 0$. Then we can say, for $\epsilon(\delta)$ small
    \begin{align}
        \Phi_{\delta'} - \Phi_\delta &\leq \left(\frac{1}{1 + 3 \epsilon(\delta')} - \frac{1}{1 + 3 \epsilon(\delta)}\right) \max\{|m_\delta|, |m_{\delta'}|\} \\
            &\leq 4(\epsilon(\delta) - \epsilon(\delta')) \max\{|m_\delta|, |m_{\delta'}|\} \\
            &\leq 4 (\epsilon(\delta) - \epsilon(\delta')) \left(2 b e^{2x} \right)
        \p
    \end{align}
    This shows that for some large $C'$
    \[
        \Phi^\nu \coloneq \Phi_{\delta_\nu} + C' \epsilon(\delta_\nu) (1 + r^2)
    \]
    will be decreasing.

    Note that
    \begin{equation}\label{eq:bigO}
        u_{\alpha, \delta}, f_\alpha, \eta_\alpha, 1 + r^2 = O(r^2) = O(e^{2x})
    \end{equation}
    and the implicit constant only depends on $b$. Thus we see that $\Phi^\nu \in \Hsing(a', b')$.

    The claim about the local Lipschitz constant follows from the fact that convolutions preserve Lipschitz constants.

    \proofstep{Step 3: Quadratic bounds on the derivatives}

    Actually \cref{eq:bigO} remains true if we replace any of the functions on the LHS with their derivative (in one of the equivariant charts), for the following reasons:
    \begin{enumerate}
        \item $u_{\alpha, \delta}$: In a chart $V_\alpha$ we have $\partial^\beta u_{\alpha, \delta} = (\partial^\beta \psi_\delta) * u_\alpha = \| \partial^\beta \psi_\delta \| O(e^{2x})$ (this bound is of course not uniform in $\nu$).
        \item $f_\alpha$: We can write $f_\alpha = r^2 + \rho_\alpha$ (see the proof of \cref{thm:KRSasymptotics}): The claim is true for $r^2$ because it has $X$-weight 2, and it was established in \cref{rmk:higherAsymp} that all derivatives of $\rho_\alpha$ are $O(x) = O(\log r)$
        \item $\eta_\alpha, 1 + r^2$: Because these functions (and hence also their derivatives) have $X$-weight 2.
    \end{enumerate}

    \proofstep{Step 4: $T$-averaging}

    Finally we show that we still satisfy all the claimed properties if we replace $\Phi^\nu$ by its $T$-average, computed as
    \[
        \langle \Phi^\nu \rangle_T \coloneq \langle \Phi^\nu - \phi \rangle_T + \phi
        \c
    \]
    where $\phi$ could be any $T$-invariant metric, for example $\Phi$
    (here $\langle - \rangle_T$ on the RHS denotes the $T$-average of a function).
    \begin{enumerate}
        \item PSh: Preserved because $T$ acts holomorphically.
        \item $i \partial\bar{\partial} \Phi^\nu_t > \epsilon(\nu) \omega_C$: Preserved since $\omega_C$ is $T$-invariant.
        \item $\Phi^\nu_t \in \Hsing(a', b')$: By \cref{eq:bigOalt} we can formulate the condition $\Phi^\nu \in \Hsing(a', b')$ as
            \[
                a' r^2 < (\Phi^\nu - \Phi) - \log \left( \frac{e^{-\Phi}}{r^{-2n} \omega_C^n} \right) < b' r^2
            \]
            for $r \gg 1$. Since the left and right of this inequality are $T$-invariant, the inequalities are preserved when passing to the $T$-average, showing $\langle \Phi^\nu \rangle_T \in \Hsing(a', b')$.
        \item nonincreasing: trivial

        \item $\Phi^\nu \to \Phi$ pointwise: While $\Phi^\nu - \Phi$ is discontinuous in general, it is continuous along any individual orbit of $T$ because $\Phi$ is $T$-invariant. Thus by Dini's theorem the convergence is uniform on every orbit, hence the average $\langle \Phi^\nu - \Phi \rangle_T$ converges uniformly on orbits, which implies pointwise convergence.

        \item quadratic bounds on derivatives: Observe that the relative potential of $r^{-2n} \omega_C^n$ and $\bigwedge_j i\d z_j^\alpha \wedge \d \bar{z}_j^\alpha$, and all of its coordinate derivatives, are bounded for every equivariant chart, because both volume forms are $X$-invariant. Hence our claim of quadratic growth of the coordinate derivatives of $\Phi^\nu$ is equivalent to claiming that their potentials w.r.t.\ $r^{-2n} \omega_C^n$, i.e.\ the functions
        \[
            f^\nu = - \log \left( \frac{e^{-\Phi^\nu}}{r^{-2n} \omega_C^n} \right)
        \]
        and their (coordinate) derivatives are $O(r^2)$ in equivariant charts. As explained in \cref{rmk:higherAsymp} this is equivalent to 
        \begin{equation} \label{eq:Or2alt}
            \left\lvert (\nabla^{\omega_C})^k f^\nu \right\rvert_{\omega_C} = O(r^{2-k})
            \p
        \end{equation}
        Since $T$ acts isometrically on $(M \setminus E, \omega_C)$, this condition is preserved if we replace $f^\nu$ by its $T$-average $\langle f^\nu \rangle_T$, which corresponds to replacing $\Phi^\nu$ by $\langle \Phi^\nu \rangle_T$.
    \end{enumerate}
\end{proof}

\section{ODEs with discontinuous right hand sides}%
\label{apx:odes}

The standard results from ODE theory (and their proofs) still hold when the RHS of an (non-autonomous) ODE only has $L^1$ regularity in time.

\begin{thm}[Picard--Lindelöf] \label{thm:PicLind}
    Let $U \subseteq \bbR$ be open, $V \in L^1(I; C^{0,1}(U; \bbR^n))$ and $K \subseteq U$ compact.
    \begin{thmlist}
        \item There exists $T > 0$ such that for every $x_0 \in K$ there exists a unique continuous function $u: [0, T] \to U$ such that
            \begin{equation} \label{eq:IVP}
                u(0) = x_0
                \qtext{and}
                \dot{u}(t) = V(t, u(t))
                \quad\text{for a.e.\ $t$}
                \p
            \end{equation}

\item The resulting function $F: [0, T] \times K \to U$ is continuous.
    \end{thmlist}
\end{thm}
\begin{rmk}
    $V(t, u(t))$ is an $L^1$ function of $t$ since
    \[
        \int_I | V(t, u(t)) | \d t \leq \int_I \|V(t)\|_{C^0} \d t \leq \| V \|
        \p
    \]
\end{rmk}
\begin{proof}
    \point1 \textit{Existence and uniqueness.}
    We apply the contraction mapping principle to the map $\Gamma$ on $C^0([0, T]; U)$ given by
    \[
        \Gamma(\phi)(t) = x_0 + \int_0^t V(t, \phi(s)) \d s
        \p
    \]
    The integrand is in $L^1$, so the RHS is indeed an (absolutely) continuous function of $t$. In addition,
    \[
        \Gamma(\phi)(t) - x_0 \leq \| V \|_{L^1([0, T])} \to 0
    \]
    as $T \to 0$, so if we pick $T \ll 1$ we indeed have $\Gamma(\phi)(t) \in U$.
    Furthermore for any $\phi, \phi' \in C^0([0, T]; U)$
    \[
        \int_0^T |V(s, \phi(s)) - V(s, \phi'(s))| \d s \leq \int_0^T \| V(s) \|_{C^{0, 1}} |\phi(s) - \phi'(s)| \leq \| \phi - \phi' \|_{C^0} \| V \|_{L^1([0, T]; C^{0,1})}
        \p
    \]
    Since the RHS goes to zero as $T \to 0$ we see that $\Gamma_T$ is a contraction for $T \ll 1$. Thus by the contraction mapping principle $T$ has a unique fixed point $u$, meaning $u$ satisfies
    \[
        u(t) = x_0 + \int_0^t V(s, u(s)) \d s
        \p
    \]
    By the Lebesgue fundamental theorem of calculus this is equivalent to $u$ being differentiable a.e.\ with $u'(t) = V(t, u(t))$.

    \medskip
    \point2 \textit{Continuous dependence on the initial value.}
    Let $x, y \in K$. Then the function $\delta(t) \coloneq | u(x, t) - u(y, t) |$ satisfies the differential inequality
    \[
        \dot{\delta}(t) \leq | \dot{u}(x, t) - \dot{u}(y, t) | \leq \| V_t \|_{C^{0,1}} \cdot \delta(t)
        \p
    \]
    Define
    \[
        \Xi(t) = \int_0^t \| V_s \|_{C^{0,1}} \d s
        \p
    \]
    Applying the product rule for almost every $t$ we see that
    \[
        \frac{d}{d t} \left[ e^{-\Xi(t)} \delta(t) \right] = e^{-\Xi(t)} \left( \delta'(t) -\Xi(t) \delta(t)\right) \leq 0
    \]
    almost everywhere. Integrating yields $\delta(t) \leq e^{\Xi(t)} \delta(0)$ so we can bound
    \[
        |u(x, t) - u(x', t')| \leq |u(x, t) - u(x', t)| + |u(x', t) - u(x', t')| \leq e^{\Xi(t)} | x - x' | + \int_t^{t'} \| V_s \|_{C^0} \d s
        \c
    \]
    showing that $u$ is continuous.
\end{proof}

Since the set of points where the partial derivative of a function exists is measurable, $\partial_t F$ exists almost everywhere in the setting of the theorem above by Tonelli's theorem.

\begin{thm}[Differentiable dependence on initial value]
    \label{thm:odeDiff}
    If in the setting of the theorem above $V \in L^1(I; C^1(U; \bbR^n))$ then the map $F: [0, T] \times K \to V$ we obtained is differentiable in the $K$-directions, and its differential satisfies the ODE
    \begin{equation} \label{eq:ODEvar}
        \frac{d}{d t} \d_x F_t = \d_{F(x, t)} V_t \circ d_x F_t
        \p
    \end{equation}
    In particular for any $x \in U$ the differential $\d_x F: [0, T] \to \bbR^{n \times n}$ is absolutely continuous.
\end{thm}
\begin{proof}
    We may assume that $U$ is convex, so that $\hat{K} \coloneq \operatorname{conv}(F(K \times [0, T])) \subseteq U$.
    Let $\delta_{\nabla V}, \delta_F: \bbR_{\geq0} \to \bbR_{\geq 0}$ be the moduli of continuity of $\nabla V$ on $\hat{K}$ and of $F$ on $K \times [0, T]$. Write $\delta = \delta_{\nabla V} \circ \delta_F$.
    Let $A_{t, x}$ be the solution of the variational equation \cref{eq:ODEvar} and set $w(t) = F_t(y) - F_t(x) - A_{t, x}(y - x)$. Then
    \begin{align} \label{eq:ODEvarerr}
        \dot{w}(t) &= V_t(F_t(y)) - V_t(F_t(x)) - \d_{F_t(x)} V_t (A_{t, x} (y - x))
        \p
    \end{align}
    We now apply the mean value theorem to every component $V^k$ of $V$ to obtain points $\xi_k$ on the line connecting $F_t(x)$ to $F_t(y)$ such that
    \begin{align}
        V^k_t(F_t(y)) - V^k_t(F_t(x)) &= \d_{\xi_k} V^k_t ( F_t(y) - F_t(x) ) \\
                                          &= \d_{F_t(x)} V^k_t ( F_t(y) - F_t(x) ) + \delta_{\nabla V}(|F_t(y) - F_t(x)|) |F_t(y) - F_t(x)|
        \p
    \end{align}
    Setting $C = \delta_{\nabla V}(\diam \hat{K})$ we observe
    \begin{align}
        \delta_{\nabla V}(|F_t(y) - F_t(x)|) |F_t(y) - F_t(x)| &\leq \delta_{\nabla V}(\delta_F(|y - x|)) | w(t) + A_{t, x}(y - x) |  \\
                                                     &\leq C |w(t)| + o(|y - x|)
    \end{align}
    because $\delta_{\nabla V} \circ \delta_F (|y - x|) = o(1)$. Plugging back into \cref{eq:ODEvarerr} we see that
    \[
        |\dot{w}(t)| \leq C | w(t) | + o(|y - x|)
        \p
    \]
    Defining $\Xi(t) \coloneq \int_0^t \|V_s\|_{C^1} \d s + C t$ we obtain
    \[
        \frac{d}{d t} e^{-\Xi(t)} |w(t)| \leq e^{-\Xi(t)} o(|y-x|)
        \c
    \]
    and integrating this equation yields
    \begin{align}
        |w(t)| \leq o(|y - x|) \cdot e^{\Xi(t)} \int_0^t e^{-\Xi(s)} \d s = o(|y - x|)
        \c
    \end{align}
    showing that $A_t$ is the derivative of $F_t$.
\end{proof}

\begin{lem} \label{thm:dtPushfwd}
     Let $F$ be as above and $\beta \in \cA^p(I \times U)$. Then $\frac{d}{d t} F_t^* \beta = F_t^* \cL_{V_t} \beta_t$.
\end{lem}
\begin{proof}
    By definition $(F_t^* \beta)(p) = \beta\left(t, F_t(p)\right) \circ \d_p F_t$. This is differentiable for almost every $t$ so we can apply product and chain rule almost everywhere to conclude the claim.
\end{proof}
\section{Proof of Berndtsson's Brunn--Minkowski inequality}%
\label{apx:bernBM}

\bernBM*

\begin{rmk}
    Our definition of $v_t$ differs from \cite{BoBM} by a factor of 2, but our definition of $\hat{u}$ coincides.
\end{rmk}

\begin{proof}
    For convenience of notation, we omit the subscripts in $\phi_t$ and $v_t$, as well as the subscript $\omega$ in norms. We also use Berndtsson's notation in which $\int_M e^{-\phi} = \int_M e^{-\phi} i^{n^2} u \wedge \bar{u}$.

    We can directly compute the LHS of \cref{eq:Berndtsson}, obtaining
    \[
        \int_M e^{-\phi} \left( \ddot{\phi} - \dot{\phi}^2 \right) + \left( \int_M e^{-\phi} \right)^{-1} \left( \int e^{-\phi} \dot{\phi} \right)^2 = \int_M e^{-\phi} \ddot{\phi} - \left( \dot{\phi} u, \pi^\phi_\perp (\dot{\phi} u) \right)_\phi
        \p
    \]
    This involves differentiating under the integral sign, which is justified by the DCT since $\dot{\phi}, \ddot{\phi}$, etc.\ are integrable with respect to $e^{-\phi}$.

    Meanwhile we can split $\iddball \phi$ into $tt$-, $MM$- and mixed parts to write out the first term on the RHS as
    \begin{align}
          &\int_M e^{-\phi} \ddot{\phi} + (-1)^{n+1} i^{n^2} \int_M e^{-\phi} \partial\bar{\partial} \phi \wedge v \wedge \overline{v} + (-1)^{n+1} 2 \Re i^{n^2} \int_M e^{-\phi} i \bar{\partial} \dot{\phi} \wedge u \wedge \overline{v}
        \label{eq:RHS1}
        \p
    \end{align}
    The first integral is clearly finite since $\ddot{\phi} = O(r^2)$.
    Recall that by our construction of $\omega$ we have $\omega \geq c i \partial\bar{\partial} r^2$ for some $c > 0$. Since $i \partial\bar{\partial} \phi = O(r^2)$ in equivariant charts this implies $i \partial\bar{\partial} \phi \leq C \omega$. This shows that the second integral is bounded up to a constant by $\| v \|_{\omega, \phi}^2$. Finally the third integral is bounded since
    \[
        \left\lVert \bar{\partial} \dot{\phi} \wedge u \right\rVert^2_\phi \leq \left\lVert |\partial \dot{\phi}|_\omega u \right\rVert^2_\phi = \int e^{-\phi} u \wedge \bar{u} \,|{\bar{\partial}} \dot{\phi}|^2_\omega
    \]
    and $|\partial \dot{\phi}|_\omega = O(r)$, so $\bar{\partial} \dot{\phi} \wedge u \in L^2_\phi \cA^{n, 1}(M; K_M^{-1})$.
    Now the last term in \cref{eq:RHS1} can be rewritten as
    \[
        2 \Re i\left(\bar{\partial} (\dot{\phi} u), \omega \wedge v \right)_\phi
            = - 2 \Re \left( \dot{\phi} u, -i \bar{\partial}^*_{\omega, \phi} (\omega \wedge v) \right)
            = - 2 \Re \left(\dot{\phi} u, \pi_\perp^\phi (\dot{\phi} u)\right)_\phi
            = - 2 \left(\dot{\phi} u, \pi_\perp^\phi (\dot{\phi} u)\right)_\phi
        \p
    \]
    To see that the second term on the RHS of \cref{eq:RHS1} is finite, recall from the proof of \cref{thm:unifL2} that $\omega \wedge v \in \Dom \Delta_{\bar{\partial}}$, where $\Delta''_\phi = \bar{\partial} \bar{\partial}^*_{\phi, \omega} + \bar{\partial}^*_{\phi, \omega} \bar{\partial}$. By \cite[Thm.~3.2]{agbook} this implies $v \in \Dom \bar{\partial}$.
    We can further use the Akizuki--Nakano identity \cite[Cor.~1.3]{agbook} to rewrite
    \[
        \| \bar{\partial} v \|^2_\phi = \left(v, \Delta''_\phi v\right) = \left(v, \Delta'_\phi v\right)_\phi + \left( v, [i \partial\bar{\partial} \phi, \Lambda] v \right)_\phi
        \c
    \]
    where $\Delta'_\phi = (\partial^\phi)^*_\phi \partial^\phi + \partial^\phi (\partial^\phi)^*_\phi$.
    Since $\Lambda \sim r^{-2}$ and $i \partial\bar{\partial} \phi = O(r^2)$ both terms on the RHS are finite.
    Since $(\partial^\phi)^*_\phi v = [\Lambda, \bar{\partial}] v = 0$ by Kähler identities, degree, and the hypothesis that $\bar{\partial} v$ is primitive, we further have
    \begin{equation} \label{eq:RHS2}
        \left(v, \Delta'_\phi v\right)_\phi = \left\| \partial^\phi v \right\|_\phi^2 = \left\| \pi^\phi_\perp (\dot{\phi} u) \right\|_\phi^2
        \p
    \end{equation}
    Also note that
    \begin{align} \label{eq:side_comp}
        \left([i \partial\bar{\partial} \phi, \Lambda ] v, v \right)_\phi &= -\left( i \partial\bar{\partial} \phi \wedge v, \omega \wedge v \right)_\phi = -\int e^{-\phi} i \partial\bar{\partial} \phi \wedge v \wedge \star \left(\omega \wedge \bar{v}\right) \\
                                                &= -\int e^{-\phi} i \partial\bar{\partial} \phi \wedge v \wedge \star (\omega \wedge \bar{v}) = (-1)^{\binom{n}{2}} (-i)^{n} \int e^{-\phi} i \partial\bar{\partial} \phi \wedge v \wedge \bar{v} \\
                                                &= i^{n^2} (-1)^{n} \int e^{-\phi} \partial\bar{\partial} \phi \wedge v \wedge \bar{v}
        \c
    \end{align}
    which cancels the second term in the RHS of \cref{eq:RHS2}.
    Combining \cref{eq:RHS1,eq:RHS2,eq:side_comp} we get that the RHS of \cref{eq:Berndtsson} is equal to
    \begin{align}
        &\int_M e^{-\phi} \ddot{\phi} - 2 (\dot{\phi} u, \pi^\phi_\perp (\dot{\phi} u))_\phi + (\pi_\perp (\dot{\phi} u), \pi_\perp (\dot{\phi} u))_\phi
            = \int_M e^{-\phi} \ddot{\phi} - \left(\dot{\phi} u, \pi^\phi_\perp (\dot{\phi} u)\right)_\phi
    \end{align}
    which is also the LHS.
\end{proof}

\end{appendix}

\linespread{1.00}
\bibliographystyle{amsalpha}
\bibliography{complex_bib.bib}

\end{document}